\documentclass[11pt]{amsart}
\usepackage[margin=0.85in]{geometry}
\usepackage{amsmath,amssymb,amsthm,mathrsfs}
\newcommand{\C}{\mathbb{C}}
\newcommand{\R}{\mathbb{R}}
\newcommand{\Z}{\mathbb{Z}}
\newcommand{\Q}{\mathbb{Q}}
\newcommand{\T}{\mathbb{T}}
\newcommand{\PP}{\mathbb{P}}
\newcommand{\E}{\mathbb{E}}
\newcommand{\A}{\mathscr{A}}
\newcommand{\Hm}{\mathcal{H}}
\newcommand{\G}{\mathcal G}
\newcommand{\Log}{\operatorname{Log}}
\newcommand{\Arg}{\operatorname{Arg}}
\newcommand{\Vol}{\operatorname{Vol}}
\newcommand{\MV}{\operatorname{MV}}
\newcommand{\Area}{\operatorname{Area}}
\newcommand{\Newt}{\operatorname{Newt}}
\newcommand{\conv}{\operatorname{conv}}
\newcommand{\dist}{\operatorname{dist}}

\newtheorem{theorem}{Theorem}[section]
\newtheorem{lemma}[theorem]{Lemma}
\newtheorem{proposition}[theorem]{Proposition}
\newtheorem{corollary}[theorem]{Corollary}
\newtheorem{conjecture}[theorem]{Conjecture}
\newtheorem*{theoremA}{Theorem A}
\newtheorem*{theoremB}{Theorem B}
\newtheorem*{theoremC}{Theorem C}
\newtheorem*{theoremD}{Theorem D}
\newtheorem*{theoremE}{Theorem E}
\newtheorem*{theoremF}{Theorem F}
\theoremstyle{definition}
\newtheorem{definition}[theorem]{Definition}
\theoremstyle{remark}
\newtheorem{remark}[theorem]{Remark}
\newtheorem{example}[theorem]{Example}

\newcounter{meq}[section]
\renewcommand{\themeq}{\thesection\alph{meq}}
\makeatletter
\newcommand{\eqtag}[1]{\addtocounter{meq}{1}\protected@edef\@currentlabel{\themeq}\label{#1}\eqno(\themeq)}
\makeatother

\title{On the volume of amoebas of complete intersections}
\author{Mounir Nisse}
\address{Department of Mathematics, Xiamen University Malaysia, Jalan Sunsuria, Bandar Sunsuria, 43900 Sepang, Selangor, Malaysia}
\email{mounir.nisse@gmail.com, mounir.nisse@xmu.edu.my}
\date{November 16, 2026}
\subjclass[2020]{Primary 14T90, 14M25; Secondary 32A60, 52A39, 52B20, 60D05, 28A75}
\keywords{Amoeba, coamoeba, Newton polytope, mixed volume, Bernstein theorem, complete intersection, twisted conjugation, Hausdorff measure, random lines}

\begin{document}

\begin{abstract}
We study the $2k$-dimensional volume of the amoeba of a smooth complete intersection $V$ of dimension $k$ in $(\C^*)^n$, $n\ge2k$, with Newton polytopes $\Delta_1,\dots,\Delta_{n-k}$. For $n=2k$ we prove $\Vol(\A(V))\ge\mathcal M(V)/\MV(\Delta_1,\dots,\Delta_k,-\Delta_1,\dots,-\Delta_k)$, where $\mathcal M(V)$ is the coamoeba volume counted with multiplicity, and we show $\mathcal M(V)\ge\pi^{2k}$ for plane curves ($k=1$, $n=2$) whose Newton polygon is not centrally symmetric, for nondegenerate varieties with odd $\MV(\Delta_1,\dots,\Delta_k,\Delta_1,\dots,\Delta_k)$, and for smooth towers of curve fibrations. We determine the equality case for curves and towers, and prove that $k$-planes in general position never attain the bound. For $n>2k$ we prove two-sided bounds for $\Hm^{2k}(\A(V))$ through coordinate projections, and show that for generic coefficients almost every fibre of the amoeba map is an orbit of a finite group. For lines we prove that in several random models the area of the amoeba of a line in $(\C^*)^n$ is asymptotic to $\frac\pi2n\log n$, and we compute its limit under degenerations of the configuration. Finally we characterise the varieties attaining the extremal volumes. For $n=2k$ every regular fibre of the amoeba map has an even number of points, divisible by $4$ for real varieties with $k\ge2$; the volume-maximal varieties are those with two-to-one amoeba maps and maximal rolled coamoebas, and they are never real when $k\ge2$; a volume-minimal plane curve is characterised by a group of symmetries of maximal order, and its amoeba has no holes. For $n>2k$ the extremal varieties lie in translates of $2k$-dimensional subtori.
\end{abstract}

\maketitle

\section{Introduction}\label{L:s1}\label{L:sec:intro}

\subsection{Amoebas and coamoebas}\label{L:s1.1}

Let $\C^*=\C\setminus\{0\}$ and let $V\subset(\C^*)^n$ be an algebraic variety, that is, the common zero set of finitely many Laurent polynomials in $n$ variables. Every point $z=(z_1,\dots,z_n)$ of $(\C^*)^n$ has a modulus part and an argument part, $z_j=|z_j|e^{i\theta_j}$. Separating them gives two maps,
$$\Log(z)=(\log|z_1|,\dots,\log|z_n|)\in\R^n,\qquad\Arg(z)=(\theta_1,\dots,\theta_n)\in\T^n=(\R/2\pi\Z)^n .$$
The image $\A(V)=\Log(V)$ is called the amoeba of $V$, and $\Arg(V)$ is its coamoeba; see [Jo], [Ni], [NiP], [NS], [PS] for coamoebas. Amoebas were introduced by Gelfand, Kapranov and Zelevinsky [GKZ, Chapter 6] in their study of discriminants, and their name comes from their shape. They are closed sets with thin tentacles going off to infinity in finitely many directions, and holes in the middle. They have become a central object in tropical geometry, real algebraic geometry and the theory of hypergeometric functions; see for instance [Mi04, PR, PT].

The simplest example is the line $L=\{1+z+w=0\}\subset(\C^*)^2$. Its amoeba is the set of $(x,y)\in\R^2$ such that the three numbers $1$, $e^x$ and $e^y$ are the side lengths of a possibly degenerate triangle, because $1+z+w=0$ means that the three complex numbers $1,z,w$ add up to zero. So it is the region bounded by the three curves $e^x+e^y=1$, $1+e^x=e^y$ and $1+e^y=e^x$, which has three tentacles, in the directions $(-1,0)$, $(0,-1)$ and $(1,1)$, and area $\pi^2/2$ [Pa]. The coamoeba of $L$ is the set of $(\theta_1,\theta_2)$ for which the unit vectors $1$, $e^{i\theta_1}$, $e^{i\theta_2}$ have positive combinations adding up to zero. It is the union of the two open triangles with vertices $(0,\pi),(\pi,\pi),(\pi,2\pi)$ and $(\pi,0),(\pi,\pi),(2\pi,\pi)$, together with parts of their boundaries. Its area is $\pi^2$, and every point of the open triangles has exactly one preimage, because the positive weights are determined by the three directions up to a common factor, which is fixed by the coefficient $1$. So the area of the amoeba is half the area of the coamoeba. This is no coincidence, as we explain below.

The combinatorics of an amoeba is governed by the Newton polytopes of the defining polynomials. If $f=\sum_\alpha c_\alpha z^\alpha$ is a Laurent polynomial, with $z^\alpha=z_1^{\alpha_1}\cdots z_n^{\alpha_n}$, its Newton polytope $\Newt(f)$ is the convex hull of the exponents $\alpha$ with $c_\alpha\neq0$. For a hypersurface $\{f=0\}$, the tentacles of the amoeba point in the directions of the outer normals of the facets of $\Newt(f)$, and the bounded components of its complement correspond to some of the lattice points inside $\Newt(f)$ [FPT].

\subsection{How large and how small can an amoeba be?}\label{L:s1.2}

For a curve $C=\{f=0\}\subset(\C^*)^2$, Passare and Rullg\aa rd [PR] proved that the area of the amoeba is at most $\pi^2\Area(\Newt(f))$, and Mikhalkin and Rullg\aa rd [MR] showed that the maximum is attained exactly, up to multiplicative translation, by the simple Harnack curves, a class of real curves with extremal topology. In higher dimension, the amoeba of a hypersurface in $(\C^*)^n$, $n\ge3$, has infinite volume, since its tentacles are thick. The natural higher-dimensional setting is that of half-dimensional varieties, $\dim_\C V=k$ in $(\C^*)^{2k}$. For these, Mikhalkin [Mik, Theorem 1 and Proposition 2] showed that the amoeba map covers the amoeba at most $\MV(\Delta_1,\dots,\Delta_k,-\Delta_1,\dots,-\Delta_k)$ times and the rolled coamoeba at most $\MV(\Delta_1,\dots,\Delta_k,\Delta_1,\dots,\Delta_k)$ times, that the multivolumes of the amoeba and of the rolled coamoeba coincide [Mik, Corollary 5 and Lemma 6], and that $\Vol(\A(V))\le\frac{\pi^{2k}}2\MV(\Delta_1,\dots,\Delta_k,\Delta_1,\dots,\Delta_k)$, which for curves is the bound of Passare and Rullg\aa rd. For $k$-planes in general position he deduced $\pi^{2k}/\binom{2k}k\le\Vol(\A(V))\le\pi^{2k}/2$ [Mik, (3.1)]\footnote{In [Mik] the lower bound of (3.1) is misprinted, with denominator $(2n!)^2$. The intended value, which follows from [Mik, Theorem 1 and Proposition 8], is $\pi^{2n}(n!)^2/(2n)!=\pi^{2n}/\binom{2n}n$, in Mikhalkin's notation $n=k$.}, using the value $\pi^{2k}$ of the multivolume of the coamoeba of a generic plane, due to Nisse and Passare [NP, Theorem 1.2] and, independently, Johansson [Jo] (see [Mik, Proposition 8]), and he asked for the minimal and maximal volumes [Mik, Problem 3.1]. As Mikhalkin notes [Mik, Remark 1.6], the finiteness of the volume had been observed earlier by Madani and the author [MN1], and the upper bounds of [Mik, Theorem 1 and Proposition 2] had been conjectured by the author. Madani and the author proved in [MN1, Theorem 1.1] that, more generally, the amoeba of a $k$-dimensional algebraic variety in $(\C^*)^n$ with $n\ge2k$ has finite $2k$-dimensional volume as soon as no irreducible component of the variety lies in a subtorus of dimension less than $2k$, the converse, that the volume is infinite otherwise, is observed in the remark following that theorem, and the proof in [MN1] is written for generic varieties. For a rational curve parametrised by $f=(f_1,\dots,f_n)$ they also gave the explicit upper bound $p\,\Area(\A_f)\le\pi^2\sum_{j_1<j_2}n_{j_1}n_{j_2}$ [MN1, Theorem 4.1], where $n_j$ is the number of zeros and poles of $f_j$ in $\C$, counted with multiplicity, and $p$ is the minimal number of points in the generic fibres of $\Log\circ f$, with equality on the left if and only if $\Log\circ f$ is a $p$-sheeted covering. In [MN2, Theorem 1.1] they showed that, for generic $k$-dimensional analytic varieties in $(\C^*)^n$ with $n\ge2k$, the finiteness of the volume of the amoeba characterises the algebraic varieties, and they proved a comparison theorem between the volumes of the amoeba and of the coamoeba, $p\,\Vol(\mathrm{co}\A)\le\Vol(\A)\le P\,\Vol(\mathrm{co}\A)$, in which $p$ and $P$ are ratios of the extreme cardinalities of the fibres of the argument map and of the logarithmic map, under the assumption that $\Log|_V$ and $\Arg|_V$ are locally finite coverings [MN2, Theorem 3.1]. They also computed, from this comparison theorem and the results of [NP], the amoeba volume $\pi^{2k}/2^k$ of a generic real $k$-plane in $(\C^*)^{2k}$ [MN2, Corollary 7.1], see also [NP, Theorem 1.2], and proved that the argument map is injective on generic linear spaces [MN2, Theorem 7.1]. These results concern finiteness, upper bounds, and the comparison of amoebas with coamoebas.

The present paper is concerned mainly with lower bounds, and, when $n>2k$, with two-sided bounds expressed through Newton polytopes. Apart from the bound of Mikhalkin for $k$-planes, we know of no earlier lower bound for the volume of amoebas of curves or of half-dimensional varieties. A lower bound in terms of the Newton polytopes cannot be expected in general. For instance, the curves $(1+z)(1+w)=\varepsilon$ have the unit square as Newton polygon, but their amoebas have area tending to zero as $\varepsilon\to0$ (Proposition \ref{L:prop:zonotope}). We shall see that this phenomenon is tied to central symmetry, and that it disappears under natural combinatorial hypotheses.

\subsection{The basic inequality}\label{L:s1.3}

Throughout, $f_1,\dots,f_k$ are Laurent polynomials in $2k$ variables with Newton polytopes $\Delta_1,\dots,\Delta_k\subset\R^{2k}$, and $V\subset(\C^*)^{2k}$ is their common zero set. We say that $V$ is a smooth complete intersection of dimension $k$ if $df_1\wedge\dots\wedge df_k\neq0$ at every point of $V$. We use the normalised mixed volume $\MV$ of $2k$ convex bodies in $\R^{2k}$, normalised so that $\MV(P,\dots,P)=(2k)!\Vol(P)$ (Section \ref{L:sec:mv}), and we put
$$M_-(V)=\MV(\Delta_1,\dots,\Delta_k,-\Delta_1,\dots,-\Delta_k),\qquad M_+(V)=\MV(\Delta_1,\dots,\Delta_k,\Delta_1,\dots,\Delta_k).$$
For $k=1$ these are $\MV(\Delta,-\Delta)=\Area(\Delta-\Delta)-2\Area(\Delta)$ and $\MV(\Delta,\Delta)=2\Area(\Delta)$.

For $x\in\R^{2k}$ and $\theta\in\T^{2k}$ let $N(x)$ be the number of points of $V$ over $x$, and $N_{\Arg}(\theta)$ the number of points of $V$ with argument $\theta$. The multivolume of the coamoeba is $\mathcal M(V)=\int_{\T^{2k}}N_{\Arg}(\theta)\,d\theta$, the volume of the coamoeba counted with multiplicity.

\begin{theoremA}
Let $V$ be a smooth complete intersection of dimension $k$ in $(\C^*)^{2k}$.

\emph{(a) (Mikhalkin [Mik])} One has $N\le M_-(V)$ almost everywhere, and $\int_{\R^{2k}}N(x)\,dx=\mathcal M(V)\le\pi^{2k}M_+(V)$.

\emph{(b)} One has
$\Vol(\A(V))\ \ge\ \frac{\mathcal M(V)}{M_-(V)} .$
Here the right side is read as $0$ when $M_-(V)=0$, in which case $N=0$ almost everywhere and $\mathcal M(V)=0$. In particular, if $\mathcal M(V)\ge\pi^{2k}$, then $\Vol(\A(V))\ge\pi^{2k}/M_-(V)$.
\end{theoremA}

Part (a) of Theorem A is due to Mikhalkin: the fibre bound $N\le M_-(V)$ and the bound for the rolled coamoeba are [Mik, Theorem 1 and Proposition 2], the identity $\int N=\mathcal M(V)$ is [Mik, Corollary 5], and the bound $\mathcal M(V)\le\pi^{2k}M_+(V)$ is [Mik, Lemma 6]. The identity holds because $V$ has exactly half the dimension of the torus, so that the two maps $\Log|_V$ and $\Arg|_V$ have Jacobians of the same absolute value [MN2, Proposition 3.1]. The fibre bound comes from the fact that a fibre of $\Log|_V$ is fixed by an antiholomorphic involution and lies in the intersection of $V$ with its image, whose points are counted by Bernstein's theorem. We reproduce the proofs in Section~\ref{L:sec:log} for completeness, since the rest of the paper uses the details. Part (b) of Theorem A is new. It follows from part (a) by integrating the fibre bound over the regular values of $\Log|_V$, and it holds for every smooth complete intersection of dimension $k$ in $(\C^*)^{2k}$. For $k$-planes in general position, where $\mathcal M(V)=\pi^{2k}$ and $M_-(V)=\binom{2k}k$, it reduces to the lower bound of [Mik, (3.1)]; Mikhalkin states no lower bound for other varieties. Part (b) also makes explicit, through the Newton polytopes, the constant in the lower half of the comparison theorem [MN2, Theorem 3.1]. Theorem A(b) reduces lower bounds for the volume to lower bounds for $\mathcal M(V)$, and the new contributions of Part~\ref{part:one} concern this coamoeba bound: the bound $\mathcal M(V)\ge\pi^{2k}$ for curves whose Newton polygon is not centrally symmetric (Theorem B), for nondegenerate varieties with $M_+(V)$ odd (Theorem C), and for smooth towers (Theorem D), together with the equality cases for curves and towers (Theorems B, E and F) and the strict inequality for $k$-planes in general position (Corollary~\ref{K2:cor:strict}). The rest of the paper is devoted to these questions. For the line $L$ above, $\mathcal M(L)=\pi^2$, $M_-(L)=2$ and every regular fibre has two points, so that $\Area(\A(L))=\pi^2/2$ is an equality case.

\subsection{Curves}\label{L:s1.4}

A lattice polygon $\Delta$ is centrally symmetric if $\Delta=2c-\Delta$ for some point $c$.

\begin{theoremB}
Let $\Delta\subset\R^2$ be a two-dimensional lattice polygon. If $\Delta$ is not centrally symmetric, then every reduced Laurent polynomial $f$ with Newton polygon $\Delta$ satisfies $\mathcal M(\{f=0\})\ge\pi^2$, and if moreover $C=\{f=0\}$ is smooth,
$$\Area(\A(C))\ \ge\ \frac{\pi^2}{\MV(\Delta,-\Delta)}\ \ge\ \frac{\pi^2}{4\Area(\Delta)},$$
with equality in the first inequality if and only if $\Delta$ is a triangle and $f$ is a trinomial, that is, $f$ has exactly three monomials, at the vertices of $\Delta$. If $\Delta$ is centrally symmetric, there are smooth curves with Newton polygon $\Delta$ for which $\mathcal M(C)$ and $\Area(\A(C))$ are arbitrarily small.
\end{theoremB}

The coamoeba bound is proved by following the arguments of the roots of $f(z,\cdot)$ along rays $z=re^{it}$, $0<r<\infty$, and using Vieta's formula. The non-symmetry of $\Delta$ enters through a pair of opposite edges of different lattice lengths, which forces the total variation of the arguments to be at least that of a sawtooth function. The equality case is more delicate. We show that an extremal curve is invariant under a finite group of symmetries containing twisted complex conjugations $z\mapsto\theta\bar z$, that the fibres of its amoeba map are orbits of this group, and that its quotient by the holomorphic part of the group is a trinomial curve on a unimodular triangle, that is, a line after a monomial change of coordinates and a multiplicative translation.

\subsection{Higher dimension}\label{L:s1.5}

\begin{theoremC}
Let $V$ be a nondegenerate smooth complete intersection of dimension $k$ in $(\C^*)^{2k}$ (Definition \ref{L:def:nondeg}). If $M_+(V)$ is odd, then for almost every $\theta\in\T^{2k}$ the number of points of $V$ whose arguments are congruent to $\theta$ modulo $\pi$ is odd. Hence $\mathcal M(V)\ge\pi^{2k}$ and $\Vol(\A(V))\ge\pi^{2k}/M_-(V)$.
\end{theoremC}

The points of $V$ with arguments congruent to $\theta$ modulo $\pi$ are the fixed points of the twisted conjugation $\tau_\theta(z)=e^{2i\theta}\bar z$ on the finite set $V\cap\tau_\theta(V)$, which by Bernstein's theorem has $M_+(V)$ points for almost every $\theta$. For curves the condition says that $\Delta$ has an odd number of boundary lattice points, which is a special case of Theorem B. In higher dimension, Theorem C covers for instance all $k$-planes in general position, for which $M_+=1$.

The second higher-dimensional result concerns fibrations. Write points of $(\C^*)^{2k}$ as $(z',z'')$ with $z'\in(\C^*)^{2k-2}$ and $z''\in(\C^*)^2$, and let $\rho:\R^{2k}\to\R^2$ be the projection to the last two coordinates.

\begin{definition}\label{L:def:tower}
A tower is a sequence of Laurent polynomials $f_1,\dots,f_k$ such that $f_j$ involves only the variables $z_1,\dots,z_{2j}$, so that $\Delta_j\subset\R^{2j}\times\{0\}$, satisfying the following conditions. Let $V_j\subset(\C^*)^{2j}$ be the zero set of $f_1,\dots,f_j$ and $\rho_j$ the projection of $\R^{2j}$ to the coordinates $2j-1,2j$. We require that $V_j$ has pure dimension $j$, that $\rho_j\Delta_j$ is two-dimensional, and that for all $p$ in $V_{j-1}$ outside an algebraic subset $B_j$ of dimension at most $j-2$, the Laurent polynomial $f_j(p,\cdot)$ in two variables is reduced, has Newton polygon $\rho_j\Delta_j$, and defines a smooth curve $C_p\subset(\C^*)^2$. We also require that $f_j(p,\cdot)$ does not vanish identically for any $p\in V_{j-1}$. For $j=1$ these conditions mean that $f_1$ is reduced with smooth zero set. The tower is called smooth if every $V_j$ is a smooth complete intersection of dimension $j$. The polygons $\rho_j\Delta_j$ are the fibre polygons of the tower, and $V=V_k$.
\end{definition}

\begin{theoremD}
Let $V$ be a tower whose fibre polygons $\rho_j\Delta_j$ are not centrally symmetric. Then $M_-(V)=\prod_j\MV(\rho_j\Delta_j,-\rho_j\Delta_j)$ and
$$\Vol(\A(V))\ \ge\ \frac{\pi^{2k}}{M_-(V)} .$$
If the tower is smooth, then moreover $\mathcal M(V)\ge\pi^{2k}$. More generally, the class of smooth complete intersections with $\mathcal M(V)\ge\pi^{2k}$ is closed under products, pullbacks by isogenies of the torus, and fibrations whose base and generic fibres belong to the class.
\end{theoremD}

\subsection{The equality case}\label{L:s1.6}

For $k=1$ the extremal curves are the trinomial curves, by Theorem B. Products of $k$ trinomial curves attain the bound in every dimension, and so do their images under monomial changes of coordinates, multiplicative translations and pullbacks by isogenies. One might expect that these are all the extremal varieties. This is not the case.

\begin{definition}\label{L:def:twisted}
Let $W\subset(\C^*)^{2k-2}$ be a variety, let $T\subset\R^2$ be a lattice triangle with vertices $\beta_0,\beta_1,\beta_2$, and let $c_0,c_1,c_2$ be Laurent polynomials in $z'$ which do not vanish identically on any component of $W$ and have no common zero on $W$. The variety
$$V=\{(p,q)\in W\times(\C^*)^2:\ c_0(p)q^{\beta_0}+c_1(p)q^{\beta_1}+c_2(p)q^{\beta_2}=0\}$$
is a twisted trinomial extension of $W$ if, for almost every $x'\in\A(W)$, the vector $(\log|c_1/c_0|,\log|c_2/c_0|)$ takes the same value at all points $p\in W$ over $x'$ at which it is defined. An iterated twisted product of trinomial curves is a variety obtained from a trinomial curve by $k-1$ successive twisted trinomial extensions.
\end{definition}

\begin{theoremE}
Let $V$ be a tower whose fibre polygons are not centrally symmetric and whose varieties $V_1,\dots,V_k$ are irreducible. Then $\Vol(\A(V))=\pi^{2k}/M_-(V)$ if and only if every fibre polygon is a triangle, $V_1$ is a trinomial curve, and each $V_j$ is a twisted trinomial extension of $V_{j-1}$ with triangle $\rho_j\Delta_j$; that is, if and only if $V$ is an iterated twisted product of trinomial curves, the successive extensions being those of the filtration $V_1\subset V_2\subset\dots\subset V_k$ of the tower.
\end{theoremE}

\begin{theoremF}
The smooth surface $V\subset(\C^*)^4$ defined by $1+z_1+z_2=0$ and $1+(2+z_1)z_3+z_4=0$ satisfies $\Vol(\A(V))=\pi^4/4=\pi^4/M_-(V)$, but it is not a toric product of two curves: it cannot be obtained from a product $C_1\times C_2$ of curves in $(\C^*)^2\times(\C^*)^2$ by any finite sequence of multiplicative translations, monomial changes of coordinates, isogenies and their inverses. An open dense subset of it is isomorphic, as an abstract variety, to an open subset of a product of two lines, through a nonmonomial isomorphism.
\end{theoremF}

The amoeba of this surface is a sheared product: over every point $x'\neq(\log2,0)$ of the amoeba of the line, its slice is a translate of the amoeba of the line by the vector $(-\frac12\log(2e^{2x'_2}-e^{2x'_1}+2),0)$. We conjecture that the iterated twisted products of trinomial curves are the only smooth complete intersections attaining the bound among those for which it holds (Conjecture \ref{L:conj:equality}).

\subsection{$k$-planes and the equality case}\label{sec:intro-planes}

For $k$-planes in general position in $(\C^*)^{2k}$ one has $M_+=1$, $M_-=\binom{2k}k$ and $\mathcal M(V)=\pi^{2k}$ [NP, Theorem 1.2], [Jo], and Mikhalkin's bound $\Vol(\A(V))\ge\pi^{2k}/\binom{2k}k$ [Mik, (3.1)] is the case of planes in Theorem A(b), and also in Theorem C. Generic real planes have amoeba volume $\pi^{2k}/2^k$ [NP, Theorem 1.2], [MN2, Corollary 7.1]. Equality would force $N=\binom{2k}k$ almost everywhere on the amoeba. In Sections~\ref{K:s1} to \ref{K2:s3} we show that this never happens. The tool is the notion of a nondegenerate corner (Definition~\ref{K:def:CN}): a point of the plane near which one block of coordinates is small, and where a certain real $(k+1)\times(k+1)$ matrix is invertible. Near such a corner the regular fibres of $\Log|_V$ have at most $2^k$ points, the number for real planes (Theorem~\ref{K:thm:main}). We first establish the existence of nondegenerate corners for $k=2$ (Theorem~\ref{K:thm:CN}). We then prove in general (Theorem~\ref{K2:thm:corner}) that every $k$-plane in general position has a nondegenerate corner, for every $k\ge2$, by an explicit expansion of the determinant near a point of $L_M$ where one coordinate vanishes. Together with the nondegeneracy of planes with nonzero coefficients (Theorem~\ref{K2:thm:nondeg}), this gives the strict inequality $\Vol(\A(V))>\pi^{2k}/\binom{2k}k$ for every $k$-plane in general position (Corollary~\ref{K2:cor:strict}), in accordance with Conjecture~\ref{L:conj:equality}.

\subsection{Varieties of higher codimension}\label{sec:intro-ngt}

When $n>2k$ the amoeba of a $k$-dimensional variety is a set of Hausdorff dimension $2k$ in $\R^n$, and its natural measure is $\Hm^{2k}$. By [MN1, Theorem 1.1] this measure is finite as soon as no component of $V$ lies in a subtorus of dimension less than $2k$, and [MN1, Theorem 4.1] gives an explicit upper bound for rational curves. Part~\ref{part:two} studies this measure. The idea is to project: an integer $2k\times n$ matrix $B$ of rank $2k$ defines a monomial map $\phi_B$ with $\Log\circ\phi_B=B\circ\Log$, so the image of $\A(V)$ under $B$ is the amoeba of the half-dimensional variety $\phi_B(V)$, to which Part~\ref{part:one} applies. In the other direction, the Cauchy--Binet formula writes the $2k$-dimensional Jacobian of $\Log|_V$ as the Euclidean norm of the vector of Jacobians of all coordinate projections. Theorem~\ref{N:thm:main} combines these facts into the bounds (\ref{N:e1.1}), (\ref{N:e1.2}) and (\ref{N:e1.3}). They involve the coamoeba multivolumes $M_I(V)$ of the coordinate projections, the essential supremum and infimum $N_{\max}$, $N_{\min}$ of the number of points of $V$ over a point of the amoeba, and mixed volumes $\mu_I^\pm(V)$ in $\R^{2n}$, which bound the fibre counts and which we compute by a graph formula (Proposition~\ref{N:prop:graph}). Corollaries~\ref{N:cor:planes} and \ref{N:cor:curves} make the bounds explicit for planes and for curves.

The fibre counts can be removed for generic varieties. Theorem~\ref{F:thm:injective} shows that if the Newton polytopes are full-dimensional and $n\ge2k+1$, then for almost every choice of coefficients almost every fibre of the amoeba map is a single orbit of the finite group $H$ dual to the lattice spanned by the differences of exponents. In particular the amoeba map is injective almost everywhere when this lattice is $\Z^n$, and Corollary~\ref{F:cor:noNmax} gives, for almost every $k$-plane in $(\C^*)^n$, $n\ge2k+1$, the bounds $\binom n{2k}^{1/2}\pi^{2k}\le\Hm^{2k}(\A(V))\le\binom n{2k}\pi^{2k}$. Theorem~\ref{F:thm:birational} shows that a generic complete intersection curve projects birationally onto plane curves with full Newton triangles, which makes Corollary~\ref{N:cor:curves} unconditional for generic coefficients.

Lines are the simplest case, and here we obtain much finer information. For a line $\ell_q=\{(\rho-q_1,\dots,\rho-q_n)\}$ the Jacobian of the amoeba map has the closed form (\ref{F:e3.1}). Nonreal lines have injective amoeba maps and real lines have two-to-one amoeba maps (Lemma~\ref{F:lem:lineN}). In $(\C^*)^3$, real lines have area at most $\frac32\pi^2$ and nonreal lines at least $\sqrt3\pi^2$ (Theorem~\ref{F:thm:lines3}). Theorem~\ref{G:thm:A} reduces the area of real lines to a continuous function $g$ of one parameter, proves its limit $\frac{1+\sqrt2}2$ at the degenerate end, and shows that the infimum of the area over all lines in $(\C^*)^3$ is an infimum over real lines, which lies between $\frac{\sqrt3}2\pi^2$ and $\frac{1+\sqrt2}2\pi^2$, and which, when it is attained, is attained only by real lines. Conjecture~\ref{G:conj:mono} predicts that it is attained exactly at the equispaced configurations. For random lines the area grows like $\frac\pi2n\log n$: we prove this, in mean and in probability, for points uniform in the disc (Theorem~\ref{G:thm:B}), for Gaussian points (Theorem~\ref{D:thm:E}) and for real lines through uniform points of a segment (Theorem~\ref{D:thm:F}). The bound $\binom n2^{1/2}\pi^2$ for nonreal lines is never attained, and it is off by a factor of order $\log n$ for almost all lines in these models (Theorem~\ref{G:thm:C}). Finally, Theorem~\ref{D:thm:D} computes the limit of the area when the points split into nested clusters organised by a rooted tree, as a sum over the internal vertices of weighted integrals, a pair of children of sizes $m,m'$ contributing $\sqrt{mm'}\,\pi^2$. For binary trees the limit of $\int J$, which is the area for nonreal families and twice the area for real ones, is at least $\frac{\pi^2}2n\log_2n$ (Corollary~\ref{D:cor:binary}), which supports the conjecture that every line satisfies a lower bound of order $n\log n$ (Conjecture~\ref{G:conj:nlogn}).

\subsection{Topology of the extremal varieties}\label{sec:intro-top}
Sections~\ref{T:s1} and \ref{T:s2} describe the varieties at which the bounds are attained. For $n=2k$, Theorem~\ref{T:thm:parity} shows, by a cobordism argument, that every regular fibre of $\Log|_V$ has an even number of points, and a number divisible by $4$ when $k\ge2$ and $V$ is real up to multiplicative translation. This explains the factor $\frac12$ in Mikhalkin's upper bound $\Vol(\A(V))\le\frac{\pi^{2k}}2M_+(V)$, improves it to $\frac{\pi^{2k}}4M_+(V)$ for real varieties with $k\ge2$, a bound attained by products of simple Harnack curves, and sharpens the lower bound of Theorem A(b) (Corollary~\ref{T:cor:bounds}). Theorem~\ref{T:thm:extremal} characterises the volume-maximal varieties as the multiHarnack varieties with two-to-one amoeba maps, which are never real when $k\ge2$, and the volume-minimal ones by the constancy of the fibre count. The group of translations by the compact torus preserving $V$ has order at most $M_+(V)$ (Proposition~\ref{T:prop:group}), and for plane curves with a Newton polygon which is not centrally symmetric the volume-minimal curves are exactly those for which this order is maximal (Theorem~\ref{T:thm:curvegroup}). Theorem~\ref{T:thm:curves} compares the two extremal classes of curves: volume-maximal curves are simple Harnack curves, whose amoebas have all the holes allowed by the Newton polygon, while the amoebas of volume-minimal curves have no holes at all. For $n>2k$, Theorems~\ref{T:thm:upper} and \ref{T:thm:lower} show that the upper and lower bounds of Theorem~\ref{N:thm:main} are attained only by varieties contained in a translate of a coordinate subtorus, respectively of an arbitrary subtorus, of dimension $2k$. So the bounds are strict for $k$-planes in general position and for all lines in general position (Corollary~\ref{T:cor:strict}).

\subsection{Organisation}\label{sec:organisation}

Part~\ref{part:one} treats the half-dimensional case. Section~\ref{L:sec:mv} collects the facts on mixed volumes, Bernstein's theorem and lattice polygons that we use. Section~\ref{L:sec:log} develops the logarithmic geometry of half-dimensional varieties and proves Theorem A. Section~\ref{L:sec:curves} proves the coamoeba bound for curves, and Section~\ref{L:sec:curve-equality} the equality case of Theorem B. Section~\ref{L:sec:parity} proves Theorem C, Section~\ref{L:sec:fibrations} proves Theorem D, and Section~\ref{L:sec:tower-equality} proves Theorems E and F. Section~\ref{L:sec:examples} gives examples and open questions. Sections~\ref{K:s1} to \ref{K:s6} introduce corners of $k$-planes and the corner criterion, and Sections~\ref{K2:s1} to \ref{K2:s3} prove that every $k$-plane in general position has a nondegenerate corner. Part~\ref{part:two} treats the case $n>2k$. Sections~\ref{N:s1} to \ref{N:s9} prove the main bounds through projections, the graph formula and the explicit corollaries. Sections~\ref{F:s1} to \ref{F:s4} prove generic injectivity, the bounds for lines in $(\C^*)^3$ and the birationality of generic space curves. Sections~\ref{G:s1} to \ref{G:s5} study lines: the minimum in $(\C^*)^3$, the growth for random lines, and lower bounds for nonreal lines. Sections~\ref{D:s1} to \ref{D:s5} prove the degeneration formula and extend the growth law to Gaussian and to real random lines. Sections~\ref{T:s1} and \ref{T:s2} describe the topology of the volume-extremal varieties, for $n=2k$ and for $n>2k$. Section~\ref{sec:open} collects the open problems. Apart from classical results (Bernstein's theorem, Kleiman's transversality theorem, the Brunn--Minkowski and Rogers--Shephard inequalities, Whitney's theorem on even functions, the area formula, Sard's theorem and basic projective geometry), all proofs are given in full.

\part{The half-dimensional case $n=2k$}\label{part:one}

\section{Polytopes, mixed volumes and Bernstein's theorem}\label{L:s2}\label{L:sec:mv}

\subsection{Mixed volumes}\label{L:s2.1}

The normalised mixed volume of convex bodies $P_1,\dots,P_n\subset\R^n$ is the unique function $\MV(P_1,\dots,P_n)$ which is symmetric, multilinear with respect to Minkowski addition and dilations by nonnegative numbers, and satisfies $\MV(P,\dots,P)=n!\Vol(P)$. Explicitly, for $\lambda_1,\dots,\lambda_n\ge0$ the volume $\Vol(\lambda_1P_1+\dots+\lambda_nP_n)$ is a homogeneous polynomial of degree $n$, and $\MV(P_1,\dots,P_n)$ is the coefficient of $\lambda_1\cdots\lambda_n$ in it. In particular
$$\Vol(\lambda P+\mu R)=\frac1{n!}\sum_{a=0}^n\binom na\MV(P^{(a)},R^{(n-a)})\lambda^a\mu^{n-a},$$
where $P^{(a)}$ denotes $a$ copies of $P$. The mixed volume is invariant under translations of each argument, nonnegative and monotone in each argument, and $\MV(gP_1,\dots,gP_n)=|\det g|\MV(P_1,\dots,P_n)$ for $g\in GL_n(\R)$ [Sch, Chapter 5]. For $n=2$, $\MV(P,R)=\Area(P+R)-\Area(P)-\Area(R)$. If $P$ is a two-dimensional convex body and $u\neq0$, then $\MV(P,[0,u])$ is $|u|$ times the width of $P$ in the direction orthogonal to $u$, which is positive.

\begin{lemma}\label{L:lem:block}
Let $L=\R^a\times\{0\}\subset\R^n$, $b=n-a$, and let $\varrho:\R^n\to\R^b$ be the projection to the last $b$ coordinates. If $P_1,\dots,P_a\subset L$ and $R_1,\dots,R_b\subset\R^n$ are convex bodies, then
$$\MV(P_1,\dots,P_a,R_1,\dots,R_b)=\MV_a(P_1,\dots,P_a)\,\MV_b(\varrho R_1,\dots,\varrho R_b),$$
with the normalised mixed volumes of $\R^a$ and $\R^b$.
\end{lemma}

\begin{proof}
First let $P_1=\dots=P_a=P$ and $R_1=\dots=R_b=R$. For $\lambda,\mu>0$ and $y\in\R^b$, the fibre of $\lambda P+\mu R$ over $y$ under $\varrho$ is $\lambda P+(\mu R)_y$, where $(\mu R)_y=\mu R\cap\varrho^{-1}(y)$ is identified with a convex subset of $\R^a$. It is empty unless $y\in\mu\varrho(R)$. By Fubini's theorem,
$\Vol_n(\lambda P+\mu R)=\int_{\mu\varrho(R)}\Vol_a(\lambda P+(\mu R)_y)\,dy$. For fixed $y$, $\Vol_a(\lambda P+(\mu R)_y)$ is a polynomial in $\lambda$ of degree at most $a$, with coefficient of $\lambda^a$ equal to $\Vol_a(P)$. Hence the coefficient of $\lambda^a$ in $\Vol_n(\lambda P+\mu R)$ is $\Vol_a(P)\Vol_b(\mu\varrho R)=\Vol_a(P)\Vol_b(\varrho R)\mu^b$. Comparing with the expansion above, $\frac1{n!}\binom na\MV(P^{(a)},R^{(b)})=\Vol_a(P)\Vol_b(\varrho R)$, so $\MV(P^{(a)},R^{(b)})=a!\Vol_a(P)\,b!\Vol_b(\varrho R)=\MV_a(P^{(a)})\MV_b((\varrho R)^{(b)})$. Now apply this identity with $P=\sum_i\lambda_iP_i$ and compare the coefficients of $\lambda_1\cdots\lambda_a$. Both sides are multilinear in the $P_i$, so we obtain the lemma with distinct $P_i$ and equal $R_j$. Then apply the result with $R=\sum_j\mu_jR_j$ and compare the coefficients of $\mu_1\cdots\mu_b$, using that $\varrho$ commutes with Minkowski sums.
\end{proof}

\begin{corollary}\label{L:cor:mvproduct}
If $P_1,\dots,P_a\subset\R^a\times\{0\}$ and $R_1,\dots,R_b\subset\{0\}\times\R^b$, then
$$\MV(P_1,\dots,P_a,R_1,\dots,R_b)=\MV_a(P_1,\dots,P_a)\,\MV_b(R_1,\dots,R_b).$$
\end{corollary}

\subsection{Bernstein's theorem}\label{L:s2.2}

Let $g_1,\dots,g_n$ be Laurent polynomials in $n$ variables with Newton polytopes $P_1,\dots,P_n$. For $w\in\R^n$, the initial form $\operatorname{in}_wg$ of $g=\sum c_\alpha z^\alpha$ is the sum of the terms $c_\alpha z^\alpha$ for which $\langle w,\alpha\rangle$ is maximal on the support of $g$. Bernstein's theorem has two parts, and we use both. First, the number of isolated common zeros of $g_1,\dots,g_n$ in $(\C^*)^n$, counted with multiplicity, is at most $\MV(P_1,\dots,P_n)$, even in the presence of positive-dimensional components of the common zero set [Be], [HS, Theorem 1]; see also [Fu, Section 5.5]. Second, if for every $w\in\R^n\setminus\{0\}$ the initial forms $\operatorname{in}_wg_1,\dots,\operatorname{in}_wg_n$ have no common zero in $(\C^*)^n$, then all common zeros are isolated and their number, counted with multiplicity, equals $\MV(P_1,\dots,P_n)$ [Be, Theorem B].

\begin{definition}\label{L:def:nondeg}
A smooth complete intersection $X\subset(\C^*)^n$ of dimension $a$ given by $g_1,\dots,g_{n-a}$ is nondegenerate if, for every $w\in\R^{n}\setminus\{0\}$, the initial forms $\operatorname{in}_wg_1,\dots,\operatorname{in}_wg_{n-a}$ define in $(\C^*)^{n}$ either the empty set or a smooth complete intersection of dimension $a$. We use it mainly for $n=2k$ and $a=k$, and in Lemma \ref{L:lem:translate} only through the consequence that each initial variety has dimension at most $a$.
\end{definition}

This is Khovanskii's nondegeneracy condition [Kh]. Only finitely many initial systems occur, and for generic coefficients with given supports each of them defines a smooth complete intersection or the empty set, by Bertini's theorem [Kh].

\subsection{Monomial maps and isogenies}\label{L:s2.3}

For an integer $n\times n$ matrix $P$ with $\det P\neq0$, let $\phi_P:(\C^*)^n\to(\C^*)^n$ be the homomorphism $\phi_P(z)_i=\prod_jz_j^{P_{ij}}$. Since $\log|\phi_P(z)_i|=\sum_jP_{ij}\log|z_j|$, one has $\Log\circ\phi_P=P\circ\Log$, and similarly $\Arg\circ\phi_P=P\circ\Arg$ modulo $2\pi$. By the Smith normal form, $\phi_P$ is surjective with finite kernel of order $|\det P|$, so it is an unramified covering of degree $|\det P|$. If $|\det P|=1$ it is an automorphism, a monomial change of coordinates. For a monomial, $(z^\alpha)\circ\phi_P=z^{P^{\mathsf T}\alpha}$, so $\Newt(f\circ\phi_P)=P^{\mathsf T}\Newt(f)$. Multiplicative translations $z\mapsto tz$, $t\in(\C^*)^n$, translate amoebas and coamoebas.

\begin{lemma}\label{L:lem:subgroup}
Let $H\subset(\C^*)^n$ be a finite subgroup. There is an integer matrix $P$ with $|\det P|=|H|$ such that $\phi_P$ is surjective with kernel $H$.
\end{lemma}

This is the classical description of finite diagonalizable groups through their character groups, see for instance [Hu, Section 16]; we include the short proof.

\begin{proof}
If $N$ is the exponent of $H$, then $H\subset\mu_N^n$, where $\mu_N$ is the group of $N$-th roots of unity. Every character of $H$ extends to $\mu_N^n$, because $\C^*$ is divisible, and every character of $\mu_N^n$ is the restriction of a monomial $z^\alpha$, $\alpha\in\Z^n$. So restriction gives a surjection $\Z^n\to\operatorname{Hom}(H,\C^*)$, whose kernel $\Lambda$ is a sublattice of index $|\operatorname{Hom}(H,\C^*)|=|H|$. Let the rows of $P$ be a basis of $\Lambda$, so that $|\det P|=|H|$. The kernel of $\phi_P$ contains $H$, since $z^\alpha=1$ on $H$ for $\alpha\in\Lambda$, and it has order $|\det P|$, so it equals $H$.
\end{proof}

\subsection{Lattice polygons}\label{L:s2.4}

For a lattice polygon $\Delta\subset\R^2$ and a primitive vector $n\in\Z^2$, let $\Delta^n$ be the face of $\Delta$ on which $\langle\cdot,n\rangle$ is maximal, and let $\ell_n(\Delta)$ be its lattice length: $\#(\Delta^n\cap\Z^2)-1$ if $\Delta^n$ is an edge, and $0$ if it is a vertex.

\begin{lemma}\label{L:lem:symmetric}
Let $\Delta$ be a two-dimensional lattice polygon. The following are equivalent: $\Delta$ is centrally symmetric; $\Delta$ is a translate of a Minkowski sum of lattice segments; $\ell_n(\Delta)=\ell_{-n}(\Delta)$ for every primitive $n$. If $\Delta$ is not centrally symmetric, there are two non-parallel primitive vectors $n,n'$ with $\ell_n(\Delta)\neq\ell_{-n}(\Delta)$ and $\ell_{n'}(\Delta)\neq\ell_{-n'}(\Delta)$.
\end{lemma}

The equivalence of the first two conditions is the classical fact that a centrally symmetric convex polygon is a zonogon, see for instance [Sch, Section 3.5]; we give a proof adapted to lattice lengths.

\begin{proof}
Let $J(n_1,n_2)=(-n_2,n_1)$. When the boundary is traversed counterclockwise, the edge with primitive outer normal $n$ has edge vector $\ell_n(\Delta)Jn$, and the sum of all edge vectors is zero. A centrally symmetric polygon has $\ell_n=\ell_{-n}$ for all $n$, since the symmetry maps the face $\Delta^n$ onto $\Delta^{-n}$. A Minkowski sum of segments is centrally symmetric. Suppose that $\ell_n=\ell_{-n}$ for all $n$. Then the edge normals come in pairs $\pm n_1,\dots,\pm n_r$, with common lattice lengths $\ell_i>0$. The zonotope $Z=\sum_i[0,\ell_iJn_i]$ has, for each $i$, exactly two edges parallel to $Jn_i$, with outer normals $\pm n_i$ and lattice length $\ell_i$, because the face of $[0,\ell_iJn_i]$ in a direction $\nu$ is the whole segment only when $\nu$ is parallel to $n_i$. So $Z$ and $\Delta$ have the same cyclically ordered edge vectors, and a convex polygon is determined up to translation by this list. Hence $\Delta=Z+c$, and $c$ is a lattice vector because vertices of both are lattice points. For the last assertion, group the edges in pairs of opposite normals, allowing one length to vanish. The vanishing of the sum of edge vectors reads $\sum_{\{n,-n\}}(\ell_n-\ell_{-n})Jn=0$. If exactly one pair had a nonzero coefficient, the sum would be a nonzero multiple of $Jn$. So at least two pairs have nonzero coefficients, and their normals are not parallel.
\end{proof}

\begin{lemma}\label{L:lem:normalize}
Let $n$ be primitive. There is $M\in SL_2(\Z)$ with $M^{\mathsf T}(0,1)=n$. For any Laurent polynomial $f$ in two variables, the polynomial $f'=f\circ\phi_{M^{\mathsf T}}$ has Newton polygon $M\Delta$, where $\Delta=\Newt(f)$. Its zero set is the image of $\{f=0\}$ under the automorphism $\phi_{M^{\mathsf T}}^{-1}$, and the top and bottom faces of $M\Delta$, with outer normals $(0,1)$ and $(0,-1)$, have lattice lengths $\ell_n(\Delta)$ and $\ell_{-n}(\Delta)$.
\end{lemma}

\begin{proof}
Since $\gcd(n_1,n_2)=1$ there are integers $u_1,u_2$ with $u_1n_2-u_2n_1=1$, and the matrix $M$ with rows $(u_1,u_2)$ and $(n_1,n_2)$ lies in $SL_2(\Z)$ and satisfies $M^{\mathsf T}(0,1)=n$. The Newton polygon of $f\circ\phi_{M^{\mathsf T}}$ is $M\Delta$. For primitive $n'$, $\langle n',M\alpha\rangle=\langle M^{\mathsf T}n',\alpha\rangle$, so the face of $M\Delta$ with outer normal $n'$ is $M(\Delta^{M^{\mathsf T}n'})$, with the same lattice length since $M$ preserves $\Z^2$. Taking $n'=(0,\pm1)$ gives the claim.
\end{proof}

We also use two classical inequalities for a convex body $K\subset\R^2$: $4\Area(K)\le\Area(K-K)\le6\Area(K)$. The left one follows from the Brunn--Minkowski inequality, with equality if and only if $K$ is centrally symmetric, and the right one is the Rogers--Shephard inequality [RS], with equality if and only if $K$ is a triangle. Hence
$$2\Area(\Delta)\le\MV(\Delta,-\Delta)=\Area(\Delta-\Delta)-2\Area(\Delta)\le4\Area(\Delta),$$
with equality on the left exactly for centrally symmetric polygons and on the right exactly for triangles. Finally, by Pick's formula $\Area(\Delta)=i+b/2-1$, where $i$ and $b$ are the numbers of interior and boundary lattice points, one has $2\Area(\Delta)\equiv b$ modulo $2$. A centrally symmetric lattice polygon $\Delta=2c-\Delta$ has $2c\in\Z^2$, because the symmetry maps vertices to vertices. So $p\mapsto2c-p$ is an involution of the boundary lattice points, without fixed points since its only fixed point $c$ is interior, and $b$ is even.

\section{Logarithmic geometry of half-dimensional varieties}\label{L:s3}\label{L:sec:log}

Let $V\subset(\C^*)^{2k}$ be a smooth complete intersection of dimension $k$. Near a point $z\in V$ we use logarithmic coordinates $\zeta=\log z$, defined locally up to $2\pi i\Z^{2k}$, and we let $L_z=\{(u_j/z_j)_j:u\in T_zV\}\subset\C^{2k}$, a complex subspace of dimension $k$. In these coordinates $\Log=\operatorname{Re}\zeta$ and $\Arg=\operatorname{Im}\zeta$, so the differentials of $\Log|_V$ and $\Arg|_V$ at $z$ are the real linear maps $\operatorname{Re}|_{L_z}$ and $\operatorname{Im}|_{L_z}$ from $L_z$ to $\R^{2k}$. The map $\Log:(\C^*)^{2k}\to\R^{2k}$ is proper, and since $V$ is closed, $\Log|_V$ is proper and $\A(V)$ is closed.

\begin{lemma}\label{L:lem:jac}
Fix a Riemannian metric on $V$, and let $J_{\Log}$ and $J_{\Arg}$ be the absolute values of the Jacobian determinants of $\Log|_V$ and $\Arg|_V$. Then $J_{\Log}=J_{\Arg}$ at every point of $V$. In particular $\Log|_V$ and $\Arg|_V$ have the same critical points.
\end{lemma}

This identity is due to Madani and the author [MN2, Proposition 3.1], who proved it for every coordinate projection of a $k$-dimensional complex submanifold of $(\C^*)^n$, and it is also used by Mikhalkin [Mik]; we include a proof for completeness.

\begin{proof}
Since $\operatorname{Im}u=\operatorname{Re}(-iu)$ and multiplication by $-i$ maps $L_z$ to itself, $\operatorname{Im}|_{L_z}=\operatorname{Re}|_{L_z}\circ m$, where $m$ is multiplication by $-i$ on $L_z$. As a real linear map, $m$ has determinant $|\det_\C m|^2=|(-i)^k|^2=1$. So the determinants of $\operatorname{Re}|_{L_z}$ and $\operatorname{Im}|_{L_z}$ have the same absolute value, for any volume form on $L_z$.
\end{proof}

For curves this lemma has a familiar form. If $t$ is a local coordinate on a curve $C$ and $d\zeta_j|_C=a_j\,dt$, then the holomorphic $2$-form $d\zeta_1\wedge d\zeta_2$ vanishes on $C$, and its real part gives $dx_1\wedge dx_2=d\theta_1\wedge d\theta_2$ on $C$. Writing $t=r+is$, direct expansion gives
$$dx_1\wedge dx_2|_C=-\operatorname{Im}(\overline{a_1}a_2)\,dr\wedge ds .\eqtag{L:e3.1}$$

\begin{lemma}[Area formula]\label{L:lem:area}
One has $\int_{\R^{2k}}N(x)\,dx=\int_VJ_{\Log}=\int_VJ_{\Arg}=\int_{\T^{2k}}N_{\Arg}(\theta)\,d\theta=\mathcal M(V)$. Let $S$ be the critical set of $\Log|_V$ and $U=\R^{2k}\setminus\Log(S)$. Then $\Log(S)$ is a closed null set, $N$ is finite and locally constant on $U$, and over every ball in $U$ the preimage in $V$ is a disjoint union of $N$ sheets mapped diffeomorphically onto the ball.
\end{lemma}

The identity $\int N=\mathcal M(V)$ of the multivolumes of the amoeba and of the coamoeba is [Mik, Corollary 5 and Lemma 6].

\begin{proof}
The set $S$ is closed and $\Log|_V$ is proper, so $\Log(S)$ is closed, and it is null by Sard's theorem. Over $U$, $\Log|_V$ is a proper local diffeomorphism, hence a finite covering, which gives the statements on $N$. On each sheet the measure $J_{\Log}\,d\mathrm{vol}$ is the pullback of the Lebesgue measure, and $\Log^{-1}(\Log(S))\setminus S$ is null, being covered by countably many open sets on which $\Log$ is a diffeomorphism. The set $S$ itself has $J_{\Log}=0$. This gives $\int N=\int_VJ_{\Log}$, the manifold form of the area formula [Fe, Section 3.2]. The map $\Arg|_V$ is not proper, so the covering argument does not apply to it. However, the general area formula for $C^1$ maps on $\sigma$-compact manifolds [Fe, Theorems 3.2.3 and 3.2.5] gives $\int_{\T^{2k}}N_{\Arg}=\int_VJ_{\Arg}$ directly, and Lemma \ref{L:lem:jac} gives the remaining equality.
\end{proof}

For $x\in\R^{2k}$ let $\sigma_x(z)=(e^{2x_j}/\bar z_j)_j$, an antiholomorphic involution of $(\C^*)^{2k}$ whose fixed point set is the torus $\Log^{-1}(x)$. For $\theta\in\T^{2k}$ let $\tau_\theta(z)=e^{2i\theta}\bar z=(e^{2i\theta_j}\bar z_j)_j$, an antiholomorphic involution whose fixed point set is $e^{i\theta}(\R^*)^{2k}$, the set of points whose arguments are congruent to $\theta$ modulo $\pi$. Indeed $z=e^{2i\theta}\bar z$ means that $e^{-i\theta}z$ is real. In logarithmic coordinates $\sigma_x(\zeta)=2x-\bar\zeta$ and $\tau_\theta(\zeta)=2i\theta+\bar\zeta$, so their differentials are $u\mapsto-\bar u$ and $u\mapsto\bar u$.

\begin{lemma}\label{L:lem:transverse}
Let $z\in V$ be a regular point of $\Log|_V$. If $\Log z=x$, then $V$ and $\sigma_x(V)$ meet transversally at $z$. If $\Arg z\equiv\theta$ modulo $\pi$, then $V$ and $\tau_\theta(V)$ meet transversally at $z$.
\end{lemma}

\begin{proof}
The logarithmic tangent spaces of $\sigma_x(V)$ and $\tau_\theta(V)$ at $z$ are $-\bar L_z$ and $\bar L_z$, and transversality means that they meet $L_z$ only in $0$, the dimensions being complementary. If $0\neq u\in L_z\cap(\pm\bar L_z)$, then $\pm\bar u\in L_z$, so $\operatorname{Re}u$ and $\operatorname{Im}u$ lie in $L_z$ and one of them is nonzero. Multiplying it by $i$ gives a nonzero purely imaginary vector of $L_z$, which is killed by $\operatorname{Re}|_{L_z}=d(\Log|_V)_z$. So $z$ is critical.
\end{proof}

\begin{proposition}[{[Mik, Theorem 1 and Proposition 2]}]\label{L:prop:fibre}
For every $x\in U$, $N(x)\le M_-(V)$.
\end{proposition}

\begin{proof}
The points of $V$ over $x$ are fixed by $\sigma_x$, so they lie in $V\cap\sigma_x(V)$. The variety $\sigma_x(V)$ is the zero set of $g_j(w)=\overline{f_j(\sigma_x(w))}=\sum_\alpha\bar c_{j,\alpha}e^{2\langle\alpha,x\rangle}w^{-\alpha}$, whose Newton polytopes are $-\Delta_j$. By Lemma \ref{L:lem:transverse}, the points of the fibre are transverse, hence isolated, common zeros of $f_1,\dots,f_k,g_1,\dots,g_k$. Bernstein's bound gives at most $M_-(V)$ of them.
\end{proof}

Let $Q_\pi=[0,\pi)^{2k}$ and, for $\theta\in Q_\pi$, $N_\pi(\theta)=\sum_{\epsilon\in\{0,\pi\}^{2k}}N_{\Arg}(\theta+\epsilon)$, the number of points of $V$ with arguments congruent to $\theta$ modulo $\pi$. Since the translates of $Q_\pi$ by the vectors $\epsilon$ tile $\T^{2k}$, $\mathcal M(V)=\int_{Q_\pi}N_\pi$. By the above, $N_\pi(\theta)$ is the number of fixed points of $\tau_\theta$ on $V\cap\tau_\theta(V)$.

\begin{proposition}[{[Mik, Theorem 1 and Proposition 2]}]\label{L:prop:argbound}
For almost every $\theta$, $N_\pi(\theta)\le M_+(V)$. Hence $\mathcal M(V)\le\pi^{2k}M_+(V)$.
\end{proposition}

\begin{proof}
The variety $\tau_\theta(V)$ is the zero set of $\overline{f_j(\tau_\theta(w))}=\sum_\alpha\bar c_{j,\alpha}e^{-2i\langle\alpha,\theta\rangle}w^\alpha$, with Newton polytopes $\Delta_j$. By Sard's theorem applied to $\Arg|_V$ and Lemma \ref{L:lem:jac}, for almost every $\theta$ all points of $V$ with arguments congruent to $\theta$ modulo $\pi$ are regular points of $\Log|_V$. By Lemma \ref{L:lem:transverse} they are isolated common zeros, and Bernstein's bound applies.
\end{proof}

Part (a) of Theorem A is therefore a restatement of results of Mikhalkin [Mik, Theorem 1, Proposition 2, Corollary 5 and Lemma 6]. Part (b) follows from it by integration. It refines, for complete intersections in $(\C^*)^{2k}$, the lower half of the comparison theorem [MN2, Theorem 3.1], which holds when $\Log|_V$ and $\Arg|_V$ are locally finite coverings: the bound $N\le M_-(V)$ makes the constant explicit in terms of the Newton polytopes.

\begin{proof}[Proof of Theorem A]
Part (a) is Proposition \ref{L:prop:fibre}, Lemma \ref{L:lem:area} and Proposition \ref{L:prop:argbound}, since $\Log(S)$ is null. For part (b), $N\le M_-(V)$ on $U$ by Proposition \ref{L:prop:fibre}, $N=0$ outside $\A(V)$, and $\R^{2k}\setminus U$ is null. Hence $\mathcal M(V)=\int N\le M_-(V)\Vol(\A(V))$ by Lemma \ref{L:lem:area}. If $M_-(V)=0$, then $N=0$ on $U$, so $\mathcal M(V)=0$, which justifies the convention in the statement.
\end{proof}

For curves, Proposition \ref{L:prop:argbound} gives $\Area(\A(C))\le\mathcal M(C)\le2\pi^2\Area(\Delta)$, weaker by a factor $2$ than the bound of Passare and Rullg\aa rd [PR]. Mikhalkin [Mik, Theorem 1] proves the sharper bound $\Vol(\A(V))\le\frac{\pi^{2k}}2M_+(V)$ in all dimensions, which for curves is the bound of [PR]; it follows from the fact that the amoeba map has degree zero, so that the two halves of $V$ on which it preserves and reverses orientation cover the amoeba equally. We use Proposition \ref{L:prop:argbound} only to know that all the quantities involved are finite.

\begin{lemma}\label{L:lem:isogeny}
Let $V_0$ be a smooth complete intersection of dimension $k$ and $V=\phi_P^{-1}(V_0)$ for an integer matrix $P$ with $\det P\neq0$. Then $V$ is a smooth complete intersection defined by the $f_j\circ\phi_P$, with Newton polytopes $P^{\mathsf T}\Delta_j$, and $N_V(x)=|\det P|N_{V_0}(Px)$, $\A(V)=P^{-1}\A(V_0)$, $\mathcal M(V)=\mathcal M(V_0)$, $M_\pm(V)=|\det P|M_\pm(V_0)$. If $V$ is the image of $V_0$ under a multiplicative translation, then $N_V$ and $\A(V)$ are translates of $N_{V_0}$ and $\A(V_0)$, and $\mathcal M$ and $M_\pm$ do not change.
\end{lemma}

\begin{proof}
Since $\phi_P$ is an unramified covering, $V$ is smooth, and $d(f_j\circ\phi_P)=\phi_P^*df_j$ have nonvanishing wedge product on $V$. A point $z'$ lies over $x'$ exactly when $\phi_P(z')$ lies over $Px'$, and each point of $V_0$ has $|\det P|$ preimages, all over the same point. This gives the formula for $N_V$ and $\A(V)=P^{-1}\A(V_0)$. The change of variables $x=Px'$ gives $\int N_V=\int N_{V_0}$, hence $\mathcal M(V)=\mathcal M(V_0)$ by Lemma \ref{L:lem:area}. The formula for $M_\pm$ follows from $\MV(gP_1,\dots)=|\det g|\MV(P_1,\dots)$ with $g=P^{\mathsf T}$. The last statement is clear.
\end{proof}

\begin{lemma}\label{L:lem:product}
Let $V_1\subset(\C^*)^{2k_1}$ and $V_2\subset(\C^*)^{2k_2}$ be smooth complete intersections of dimensions $k_1,k_2$. Then $V=V_1\times V_2$ satisfies $N_V(x_1,x_2)=N_{V_1}(x_1)N_{V_2}(x_2)$, $\A(V)=\A(V_1)\times\A(V_2)$, $\mathcal M(V)=\mathcal M(V_1)\mathcal M(V_2)$ and $M_\pm(V)=M_\pm(V_1)M_\pm(V_2)$.
\end{lemma}

\begin{proof}
The first two statements follow from $\Log(z_1,z_2)=(\Log z_1,\Log z_2)$, and the third follows by integration and Lemma \ref{L:lem:area}. The polytopes of $V_1$ and their negatives lie in $\R^{2k_1}\times\{0\}$, and those of $V_2$ in $\{0\}\times\R^{2k_2}$. After reordering the arguments, Corollary \ref{L:cor:mvproduct} gives $M_\pm(V)=M_\pm(V_1)M_\pm(V_2)$.
\end{proof}

\section{Curves: the coamoeba bound}\label{L:s4}\label{L:sec:curves}

In this section and the next, $f$ is a Laurent polynomial in two variables $(z,w)$, $\Delta=\Newt(f)$ is two-dimensional, and $C=\{f=0\}\subset(\C^*)^2$. We assume that $f$ is reduced, that is, without repeated irreducible factors in the factorial ring $\C[z^{\pm1},w^{\pm1}]$. Then $C$ has finitely many singular points, and we put $\mathcal M(C)=\int_{C_{\mathrm{reg}}}J_{\Arg}$, which is also the integral of the number of points of $C_{\mathrm{reg}}$ with given argument.

\subsection{The variation of the arguments along rays}\label{L:s4.1}

Let $n$ be primitive with $\ell_n(\Delta)\neq\ell_{-n}(\Delta)$. By Lemma \ref{L:lem:normalize}, after a monomial change of coordinates, which preserves $\mathcal M(C)$ by Lemma \ref{L:lem:isogeny}, we may assume that the bottom and top faces of $\Delta$, with outer normals $(0,-1)$ and $(0,1)$, have different lattice lengths $a$ and $b$. Let $y_0<y_1$ be the minimal and maximal second coordinates on $\Delta$, $m=y_1-y_0\ge1$, and write $f(z,w)=\sum_{q=y_0}^{y_1}P_q(z)w^q$. Let $c(z)$ be a greatest common divisor of the Laurent polynomials $P_q$ in $\C[z^{\pm1}]$ and $f=c\,h$. Then $\Newt(f)=\Newt(c)+\Newt(h)$, and $\Newt(c)$ is a horizontal segment. So the bottom and top faces of $\Newt(h)$ have lattice lengths $a-\deg c$ and $b-\deg c$, whose difference is still $a-b\neq0$, where $\deg c$ is the lattice length of $\Newt(c)$. The zero set of $c$ is a finite union of vertical lines $\{z=z_0\}$, on which $\Arg$ has one-dimensional image, so that they do not contribute to $\mathcal M(C)$. Outside these lines, $C$ coincides with $\{h=0\}$. Moreover $\Newt(h)$ is two-dimensional, because a segment has opposite faces of equal lattice length while the bottom and top faces of $\Newt(h)$ have different lengths. Hence we may and do assume that the $P_q$ have no common factor, keeping the notation $a\neq b$ for the lengths of the bottom and top faces.

The bottom and top faces are $[i_0,i_0+a]\times\{y_0\}$ and $[i_1,i_1+b]\times\{y_1\}$, and the coefficients of $f$ at their endpoints are nonzero, since these are vertices of $\Delta$.

\begin{lemma}\label{L:lem:roots}
There is a finite set $Z_0\subset\C^*$ such that, over $\C^*\setminus Z_0$, the projection $(z,w)\mapsto z$ restricted to $C$ is an $m$-sheeted covering by smooth points of $C$. On every simply connected open subset $\Sigma$ of $\C^*\setminus Z_0$ the sheets are graphs of holomorphic functions $w_1,\dots,w_m$ without zeros, and $\prod_jw_j(z)=R(z)=(-1)^mP_{y_0}(z)/P_{y_1}(z)$.
\end{lemma}

\begin{proof}
Since $f$ is reduced and its coefficients $P_q$ have no common factor, $w^{-y_0}f$ is squarefree as a polynomial in $w$ over the field $\C(z)$, by Gauss's lemma. So its discriminant is a nonzero Laurent polynomial in $z$. Let $Z_0$ be the set of zeros of $P_{y_0}P_{y_1}$ and of the discriminant, together with the projections of the finitely many singular points of $C$. For $z\notin Z_0$, the polynomial $w^{-y_0}f(z,w)$ has degree $m$, nonzero leading coefficient $P_{y_1}(z)$, nonzero constant term $P_{y_0}(z)$ and no multiple root. So it has $m$ distinct roots in $\C^*$, and the holomorphic implicit function theorem gives local root functions whose graphs exhaust the preimage. This is a covering, trivial over simply connected sets, and Vieta's formula gives the product.
\end{proof}

\begin{lemma}\label{L:lem:phase}
If $\Theta:(0,\infty)\to\R$ is continuous and $e^{i\Theta(r)}\to e^{i\gamma}$ as $r\to\infty$, then $\Theta(r)$ has a finite limit in $\gamma+2\pi\Z$. The same holds as $r\to0$.
\end{lemma}

\begin{proof}
For large $r$, the values $\Theta(r)$ lie in the union of the disjoint intervals $(\gamma+2\pi q-\pi/3,\gamma+2\pi q+\pi/3)$, $q\in\Z$. By connectedness, the image of a tail interval lies in one of them, and exponentiation is injective there.
\end{proof}

\begin{proposition}\label{L:prop:rays}
If $\Delta$ is not centrally symmetric, then $\mathcal M(C)\ge\pi^2$.
\end{proposition}

\begin{proof}
By Lemma \ref{L:lem:symmetric} there is $n$ with $\ell_n(\Delta)\neq\ell_{-n}(\Delta)$, and we use the normalisation above. Cut the $z$-plane along the finitely many rays from the origin through the points of $Z_0$, and along the positive real ray. The complementary open sectors $\Sigma=\{re^{it}:r>0,\ \alpha<t<\beta\}$ are simply connected and avoid $Z_0$, and the preimage of the omitted rays is a null set for $J_{\Arg}$. Fix a sector and the root functions of Lemma \ref{L:lem:roots}, and let $\vartheta_j(r,t)$ be a smooth real lift of $\arg w_j(re^{it})$. On the graph of $w_j$, parametrised by $(r,t)$, one has $\Arg=(t,\vartheta_j)$, so $d\theta_1\wedge d\theta_2=dt\wedge d\vartheta_j=-\partial_r\vartheta_j\,dr\wedge dt$ and $J_{\Arg}\,d\mathrm{vol}=|\partial_r\vartheta_j|\,dr\,dt$. Tonelli's theorem gives
$$\int_{C\cap(\Sigma\times\C^*)}J_{\Arg}=\int_\alpha^\beta\mathcal V(t)\,dt,\qquad\mathcal V(t)=\sum_{j=1}^m\int_0^\infty|\partial_r\vartheta_j(r,t)|\,dr .$$
Put $\Theta_t(r)=\sum_j\vartheta_j(r,t)$, a continuous lift of the argument of $R(re^{it})$. As $r\to\infty$, the dominant terms of $P_{y_0}$ and $P_{y_1}$ are those of degrees $i_0+a$ and $i_1+b$, and as $r\to0$ those of degrees $i_0$ and $i_1$, with nonzero coefficients. Hence $R(re^{it})/|R(re^{it})|$ tends to $e^{i(\gamma_\infty+(i_0+a-i_1-b)t)}$ as $r\to\infty$ and to $e^{i(\gamma_0+(i_0-i_1)t)}$ as $r\to0$, for constants $\gamma_\infty,\gamma_0$. By Lemma \ref{L:lem:phase}, $\Theta_t$ has limits at $0$ and $\infty$, and $J(t)=\Theta_t(\infty)-\Theta_t(0)\in\kappa+(a-b)t+2\pi\Z$ with $\kappa=\gamma_\infty-\gamma_0$. For $0<\varepsilon<\rho$, $\mathcal V(t)\ge\sum_j|\vartheta_j(\rho,t)-\vartheta_j(\varepsilon,t)|\ge|\Theta_t(\rho)-\Theta_t(\varepsilon)|$, and in the limit $\mathcal V(t)\ge|J(t)|\ge\dist(\kappa+(a-b)t,2\pi\Z)$. Summing over the sectors,
$$\mathcal M(C)\ge\int_0^{2\pi}\dist(\kappa+dt,2\pi\Z)\,dt,\qquad d=a-b\neq0 .$$
The function $g(s)=\dist(s,2\pi\Z)$ is $2\pi$-periodic with $\int_0^{2\pi}g=2\int_0^\pi s\,ds=\pi^2$. As $t$ runs over $[0,2\pi]$, $s=\kappa+dt$ runs over $|d|$ periods with $dt=ds/|d|$, so the integral equals $\pi^2$.
\end{proof}

\begin{theorem}\label{L:thm:curves}
Let $f$ be reduced with two-dimensional Newton polygon $\Delta$ not centrally symmetric, and suppose that $C=\{f=0\}$ is smooth. Then $\Area(\A(C))\ge\pi^2/\MV(\Delta,-\Delta)\ge\pi^2/(4\Area(\Delta))$, the last inequality being strict unless $\Delta$ is a triangle.
\end{theorem}

\begin{proof}
Since $f$ is reduced and $C$ is smooth, $df\neq0$ on $C$, by the Jacobian criterion, so $C$ is a smooth complete intersection of dimension one. Theorem A and Proposition \ref{L:prop:rays} give the first inequality, and Section \ref{L:sec:mv} gives $\MV(\Delta,-\Delta)\le4\Area(\Delta)$, with equality exactly for triangles.
\end{proof}

\subsection{Lines and trinomials}\label{L:s4.2}

\begin{lemma}\label{L:lem:line}
The line $L=\{1+z+w=0\}$ satisfies $\mathcal M(L)=\pi^2$, $N=2$ almost everywhere on $\A(L)$, and $\Area(\A(L))=\pi^2/2$. Its Newton polygon $T_0=\conv(0,e_1,e_2)$ has $\MV(T_0,-T_0)=2$ and $M_+(L)=1$.
\end{lemma}

\begin{proof}
A point of $L$ with argument $\theta$ is $(r_1e^{i\theta_1},r_2e^{i\theta_2})$ with $r_1,r_2>0$ and $1+r_1e^{i\theta_1}+r_2e^{i\theta_2}=0$. If $\theta_1\not\equiv\theta_2$ modulo $\pi$, this real linear system has a unique solution $(r_1,r_2)\in\R^2$, so $N_{\Arg}(\theta)\le1$, with equality exactly when both entries are positive. That happens exactly on the two open triangles described in Section \ref{L:sec:intro}, of total area $\pi^2$, so $\mathcal M(L)=\pi^2$. The line is defined over $\R$, so complex conjugation preserves $L$ and $\Log$, and its fixed points form the real curve $L(\R)$, whose image under $\Log$ is null. Over a point of $U\setminus\Log(L(\R))$, the fibre is a union of pairs $\{z,\bar z\}$ with $z\neq\bar z$, so $N$ is even and positive on the amoeba. The hexagon $T_0-T_0$ has area $3$, so $\MV(T_0,-T_0)=3-\frac12-\frac12=2$, and Proposition \ref{L:prop:fibre} gives $N\le2$. Hence $N=2$ almost everywhere on $\A(L)$, and $\Area(\A(L))=\frac12\int N=\frac12\mathcal M(L)=\pi^2/2$. Finally $M_+(L)=2\Area(T_0)=1$.
\end{proof}

\begin{proposition}\label{L:prop:trinomial}
Let $T$ be a lattice triangle with vertices $\alpha_0,\alpha_1,\alpha_2$ and $f=c_0z^{\alpha_0}+c_1z^{\alpha_1}+c_2z^{\alpha_2}$ with $c_0c_1c_2\neq0$. Then $C=\{f=0\}$ is smooth and irreducible, $\MV(T,-T)=4\Area(T)$, every regular fibre of $\Log|_C$ over $\A(C)$ has exactly $4\Area(T)$ points, $\mathcal M(C)=\pi^2$, and $\Area(\A(C))=\pi^2/(4\Area(T))=\pi^2/\MV(T,-T)$.
\end{proposition}

As the proof shows, the proposition follows from the case of the line (Lemma \ref{L:lem:line}) by a monomial map. We state it because the fibre count and the irreducibility are used in the equality case.

\begin{proof}
Write $T=gT_0+c$ with $g\in GL_2(\R)$. Then $T-T=g(T_0-T_0)$ has area $3|\det g|=6\Area(T)$, so $\MV(T,-T)=4\Area(T)$. Let $u=\alpha_1-\alpha_0$, $v=\alpha_2-\alpha_0$, and let $P$ be the matrix with rows $u,v$, so that $|\det P|=2\Area(T)$ and $\phi_P(z)=(z^u,z^v)$. Then $f=c_0z^{\alpha_0}\cdot(1+\frac{c_1}{c_0}y_1+\frac{c_2}{c_0}y_2)\circ\phi_P$, so $C=\phi_P^{-1}(L')$, where $L'$ is the line $1+\frac{c_1}{c_0}y_1+\frac{c_2}{c_0}y_2=0$, a multiplicative translate of $L$. By Lemmas \ref{L:lem:isogeny} and \ref{L:lem:line}, $C$ is smooth, $N_C=2\Area(T)\cdot N_{L'}\circ P=4\Area(T)$ almost everywhere on $\A(C)=P^{-1}\A(L')$, $\mathcal M(C)=\pi^2$, and $\Area(\A(C))=(\pi^2/2)/|\det P|$. For irreducibility, the line $L'$ is a projective line minus three points. Small loops around the points where $y_1=0$ and $y_2=0$ have classes $(1,0)$ and $(0,1)$ in $\pi_1((\C^*)^2)=\Z^2$. So $\pi_1(L')\to\pi_1((\C^*)^2)$ is surjective, and the covering $\phi_P^{-1}(L')\to L'$ is connected. A smooth connected curve is irreducible.
\end{proof}

\subsection{Centrally symmetric polygons}\label{L:s4.3}

\begin{proposition}\label{L:prop:zonotope}
Let $\Delta$ be a two-dimensional centrally symmetric lattice polygon. For every $\eta>0$ there is a reduced $f$ with $\Newt(f)=\Delta$ and $C=\{f=0\}$ smooth, such that $\mathcal M(C)<\eta$ and $\Area(\A(C))<\eta$. For the unit square, $f$ can be taken irreducible, of the form $(1+z)(1+w)-\varepsilon$.
\end{proposition}

\begin{proof}
By Lemma \ref{L:lem:symmetric}, after a translation, $\Delta=\sum_{i=1}^r[0,u_i]$ with $u_i\in\Z^2\setminus\{0\}$. For nonzero constants $c_i$ and $\varepsilon$ put $g=\prod_i(1-c_iz^{u_i})$ and $f_\varepsilon=g-\varepsilon$. A vertex of $\Delta$ is $\sum_{i\in I}u_i$, where $I$ is the set of $i$ with $\langle\nu,u_i\rangle>0$ for $\nu$ in the interior of its normal cone, and its coefficient in $g$ is $\prod_{i\in I}(-c_i)\neq0$, since only one choice of terms produces it. Subtracting $\varepsilon$ changes only the coefficient of the vertex $0$, which becomes $1-\varepsilon$. Here we have chosen the signs of the $u_i$, which does not change the zonotope up to translation, so that all $u_i$ lie in an open half-plane; then $0$ is a vertex of $\Delta$. So $\Newt(f_\varepsilon)=\Delta$ for $\varepsilon\neq1$. The regular function $g$ on $(\C^*)^2$ has finitely many critical values, so for all but finitely many $\varepsilon$ the curve $g^{-1}(\varepsilon)$ is smooth and $f_\varepsilon$ is reduced.

Let $p\in C_\varepsilon=\{f_\varepsilon=0\}$ and $\theta=\Arg p$. Then $z^{u_i}(p)=t_ie^{i\langle u_i,\theta\rangle}$ with $t_i>0$, and since the distance from $1$ to the line $\R e^{i\psi}$ is $|\sin\psi|$, $|1-c_iz^{u_i}(p)|\ge|\sin(\langle u_i,\theta\rangle+\arg c_i)|$. Hence $\gamma(\theta)=\prod_i|\sin(\langle u_i,\theta\rangle+\arg c_i)|\le|g(p)|=|\varepsilon|$, and the coamoeba of $C_\varepsilon$ lies in $\{\gamma\le|\varepsilon|\}$. These sets decrease to the zero set of $\gamma$, a finite union of closed curves, so their measure tends to zero. By Proposition \ref{L:prop:argbound}, $N_{\Arg}\le M_+=2\Area(\Delta)$ almost everywhere. So $\mathcal M(C_\varepsilon)\le2\Area(\Delta)\Area\{\gamma\le|\varepsilon|\}\to0$, and $\Area(\A(C_\varepsilon))\le\int N=\mathcal M(C_\varepsilon)$. For the unit square, the curve $(1+z)(1+w)=\varepsilon$ is the graph of $w=\varepsilon/(1+z)-1$ over $\C^*\setminus\{-1,\varepsilon-1\}$, the points where $z$ or $w$ would vanish being removed, hence irreducible, and it is smooth for $\varepsilon\neq0$.
\end{proof}

Proposition \ref{L:prop:rays} and Proposition \ref{L:prop:zonotope} together give the first and last assertions of Theorem B. Theorem \ref{L:thm:curves} gives the inequality, and the equality case is proved in the next section.

\section{Curves: the equality case}\label{L:s5}\label{L:sec:curve-equality}

In this section $f$ is irreducible, $C=\{f=0\}$ is smooth, $\Delta=\Newt(f)$ is two-dimensional and not centrally symmetric, $\mu=\MV(\Delta,-\Delta)$, and we assume $\Area(\A(C))=\pi^2/\mu$. We prove that $\Delta$ is a triangle and $f$ a trinomial. Let $S$ be the critical set of $\Log|_C$, $B=\Log(S)$, $U=\R^2\setminus B$, and $N(x)$ the number of points of $C$ over $x\in U$. Put $A^\circ=\A(C)\cap U$, an open set.

\begin{lemma}\label{L:lem:eq-first}
One has $\mathcal M(C)=\pi^2$, and $N(x)=\mu$ for every $x\in A^\circ$.
\end{lemma}

\begin{proof}
By Lemma \ref{L:lem:area}, Proposition \ref{L:prop:rays} and Proposition \ref{L:prop:fibre}, $\pi^2\le\mathcal M(C)=\int_{A^\circ}N\le\mu\Area(A^\circ)=\mu\Area(\A(C))=\pi^2$. So $N=\mu$ almost everywhere on $A^\circ$, and since $N$ is locally constant on $A^\circ$, everywhere.
\end{proof}

Fix a primitive $n$ with $\ell_n(\Delta)\neq\ell_{-n}(\Delta)$, let $\Phi=\phi_{M^{\mathsf T}}$ be the automorphism of Lemma \ref{L:lem:normalize}, and $C'=\Phi^{-1}(C)$, with unequal top and bottom lengths. By Lemma \ref{L:lem:isogeny}, $\mathcal M(C')=\pi^2$, and $\Phi$ maps the critical set $S'$ of $\Log|_{C'}$ onto $S$. On $C'$, with coordinates $(z,w)$, the logarithmic Gauss map [Ka], [Mi00], [Ni] is the meromorphic function $\gamma=d\zeta_2/d\zeta_1|_{C'}=a_2/a_1$ in the notation of (\ref{L:e3.1}). It is not constant. Otherwise $d\zeta_2=c\,d\zeta_1$ on $C'$. Near a puncture of the smooth compactification of $C'$, $z$ and $w$ behave like $s^{\alpha}$ and $s^{\beta}$ in a local parameter $s$, with integers $(\alpha,\beta)\neq(0,0)$, so $c=\beta/\alpha$ is rational. Then $w^{\alpha}z^{-\beta}$ is locally constant on $C'$, and $C'$ is a coset of a one-dimensional subtorus, with a segment as Newton polygon. By (\ref{L:e3.1}), a point $p$ with $a_1(p)\neq0$ is critical if and only if $\operatorname{Im}\gamma(p)=0$, since $\operatorname{Im}(\overline{a_1}a_2)=|a_1|^2\operatorname{Im}\gamma$, and points with $a_1=0$ are critical. So $S'=\gamma^{-1}(\R\cup\{\infty\})$. For $t\in\R/2\pi\Z$ let $C'_t=\{p\in C':\arg z(p)=t\}$.

\begin{lemma}\label{L:lem:monotone}
For almost every $t$, the function $\operatorname{Im}\gamma$ does not take both a positive and a negative value on $C'_t$.
\end{lemma}

\begin{proof}
With the notation of the proof of Proposition \ref{L:prop:rays} for $C'$, along the graph of $w_j$ over a ray, $\partial_r\vartheta_j(r,t)=\operatorname{Im}\partial_r\log w_j(re^{it})=\frac1r\operatorname{Im}\gamma(re^{it},w_j(re^{it}))$, since $d\log w_j/d\log z=\gamma$ and $d\log z/dr=1/r$. Since $\mathcal M(C')=\pi^2$, all the inequalities in that proof are equalities, so $\mathcal V(t)=|J(t)|$ for almost every $t$. But $\mathcal V(t)=\int_0^\infty\sum_j|\partial_r\vartheta_j|\,dr\ge|\int_0^\infty\sum_j\partial_r\vartheta_j\,dr|=|J(t)|$, with equality only if the continuous functions $\partial_r\vartheta_j(\cdot,t)$ are all nonnegative or all nonpositive. Over the sectors, $C'_t$ is the union of the graphs of the $w_j$ over the ray of angle $t$.
\end{proof}

\begin{lemma}\label{L:lem:critical-structure}
There is a finite set $F'\subset S'$ such that $S'\setminus F'$ is a one-dimensional real analytic submanifold, and $S'$ has no isolated points. Every $p\in S'\setminus F'$ has a neighbourhood $W$ in which $S'\cap W$ is an arc separating $W$ into two components, on which $\operatorname{Im}\gamma$ has opposite signs.
\end{lemma}

\begin{proof}
Let $F'$ be the set of poles of $\gamma$ and of points of $\gamma^{-1}(\R)$ where $d\gamma=0$, a finite set since $\gamma$ is a nonconstant rational function on the smooth affine curve $C'$. At $p\in S'\setminus F'$, $\gamma$ is a local biholomorphism onto a neighbourhood of a real number, so $S'$ is locally the preimage of a real interval, and its two sides map to the two half-planes. At a point of $F'$, $\gamma$ has a local normal form $u\mapsto u^m$ as a map to the Riemann sphere, and the preimage of the real projective line consists of $2m$ arcs issuing from the point.
\end{proof}

\begin{lemma}\label{L:lem:constant-arg}
Let $n^\perp=(n_2,-n_1)$. On every connected component of $S\setminus\Phi(F')$, the function $\langle n^\perp,\Arg\rangle$ is constant modulo $2\pi$.
\end{lemma}

\begin{proof}
Since $\Arg\circ\Phi=M^{\mathsf T}\Arg$ and the first column of $M^{-1}$ is $(n_2,-n_1)$, the function $\arg z$ on $C'$ corresponds under $\Phi$ to $\langle n^\perp,\Arg\rangle$ on $C$. So it suffices that $\arg z$ is locally constant on $S'\setminus F'$. Otherwise there is an arc $\sigma\subset S'\setminus F'$ on which $\arg z$ has nonzero derivative at a point $p_1$. Then the level curves $C'_t$ near $p_1$ form a smooth foliation transverse to $\sigma$. By Lemma \ref{L:lem:critical-structure}, for every $t$ in an interval, $C'_t$ crosses $\sigma$ near $p_1$ and contains points with $\operatorname{Im}\gamma>0$ and points with $\operatorname{Im}\gamma<0$. This contradicts Lemma \ref{L:lem:monotone}.
\end{proof}

Let $\G$ be the group of the maps $z\mapsto\vartheta z$ and $z\mapsto\vartheta\bar z$ of $(\C^*)^2$, $\vartheta\in(S^1)^2$, with coordinatewise operations. Every $g\in\G$ satisfies $\Log\circ g=\Log$. Let $\Gamma_C=\{g\in\G:g(C)=C\}$, $\Gamma_C^+$ its subgroup of translations and $\Gamma_C^-=\Gamma_C\setminus\Gamma_C^+$.

\begin{lemma}\label{L:lem:extend}
If $g\in\G$ and $g(C)\cap C$ is infinite, then $g\in\Gamma_C$. The group $\Gamma_C$ is finite.
\end{lemma}

\begin{proof}
If $g(z)=\vartheta z$, then $g(C)$ is the zero set of $f(\vartheta^{-1}z)$, and if $g(z)=\vartheta\bar z$, then $g(C)$ is the zero set of the Laurent polynomial $\overline{f(\vartheta\bar z)}$. Both are irreducible, being images of $f$ under ring automorphisms, possibly composed with conjugation of coefficients. Two irreducible curves with infinitely many common points coincide. The set $\{t\in(\C^*)^2:tC=C\}$ is a Zariski closed subgroup. If it were infinite, its identity component would be a subtorus of positive dimension, and $C$ would be a coset of a one-dimensional subtorus, with a segment as Newton polygon. So $\Gamma_C^+$ is finite, and $\Gamma_C^-$ is empty or a coset of $\Gamma_C^+$.
\end{proof}

\begin{lemma}\label{L:lem:fixed}
Every antiholomorphic $g\in\G$ is an involution. If $g\in\Gamma_C^-$, its fixed point set $\operatorname{Fix}_C(g)$ in $C$ is a closed real analytic subset of dimension at most one contained in $S$. A point of $C$ is fixed by at most one element of $\Gamma_C\setminus\{\mathrm{id}\}$, which is antiholomorphic.
\end{lemma}

\begin{proof}
For $g(z)=\vartheta\bar z$, $g(g(z))=\vartheta\bar\vartheta z=z$. The fixed point set is real analytic, and it is not all of $C$, since $g|_C$ is antiholomorphic. At a fixed point $p$, $dg_p$ reverses the orientation of $T_pC$, and $d\Log_p\circ dg_p=d\Log_p$, so $\det d\Log_p=-\det d\Log_p$ and $p\in S$. A translation with $\vartheta\neq(1,1)$ has no fixed point, and if two antiholomorphic elements $g,g'$ fix $p$, then $g'g$ is a translation fixing $p$, so $g'=g$.
\end{proof}

\begin{proposition}\label{L:prop:S-fixed}
The set $\Gamma_C^-$ is nonempty, and $S=\bigcup_{g\in\Gamma_C^-}\operatorname{Fix}_C(g)$.
\end{proposition}

\begin{proof}
Let $n,n'$ be non-parallel primitive vectors as in Lemma \ref{L:lem:symmetric}, and let $F$ be the union of the two finite sets of Lemma \ref{L:lem:constant-arg} for the two normalisations. On a connected component $\Sigma$ of $S\setminus F$, the functions $\langle n^\perp,\Arg\rangle$ and $\langle n'^\perp,\Arg\rangle$ are constant modulo $2\pi$, and since $n^\perp,n'^\perp$ are independent, $\Arg$ takes finitely many values on the connected set $\Sigma$, hence is constant, equal to some $\theta_0$. Then $\Sigma\subset e^{i\theta_0}(0,\infty)^2$, which is fixed pointwise by $g_0(z)=e^{2i\theta_0}\bar z$. So $g_0(C)\cap C$ is infinite and $g_0\in\Gamma_C^-$ by Lemma \ref{L:lem:extend}. Hence $S\setminus F$ is contained in the union of the fixed point sets. The set $S$ is not empty. Indeed, for a vertex $v$ of $\Delta$ and $\nu$ in the interior of its normal cone, the monomial $c_vz^v$ dominates the others on $\Log^{-1}(s\nu)$ for $s$ large, so the amoeba is a closed proper subset of $\R^2$. Its boundary points are critical values, because $\Log|_C$ is open near regular points. Since $S$ has no isolated points, $S\setminus F$ is dense in $S$, so $\Gamma_C^-\neq\emptyset$. The union of the fixed point sets is closed, since $\Gamma_C$ is finite, so it contains $S$. The reverse inclusion is Lemma \ref{L:lem:fixed}.
\end{proof}

Let $Q=C/\Gamma_C$, with quotient map $q$. By Lemma \ref{L:lem:fixed}, $\Gamma_C$ acts freely on $C\setminus S$, and a point $p\in S$ has stabiliser $\{\mathrm{id},g\}$ for a unique $g\in\Gamma_C^-$. Near $p$ there is a holomorphic coordinate $u$ in which $g$ is $u\mapsto\bar u$, by the local linearisation of an antiholomorphic involution with a fixed point. Hence $Q$ is a topological surface with boundary $\partial Q=q(S)$, whose interior $q(C\setminus S)$ is covered by $C\setminus S$ with deck group $\Gamma_C$. Since $\Log\circ g=\Log$ on $\Gamma_C$, there is a continuous proper map $\Lambda:Q\to\R^2$ with $\Lambda\circ q=\Log|_C$, and $\Lambda(\partial Q)=B$.

\begin{lemma}\label{L:lem:L-local}
The map $\Lambda$ is a local diffeomorphism on the interior of $Q$. There is a finite set $F_Q\subset\partial Q$ such that every $b\in\partial Q\setminus F_Q$ has a neighbourhood which $\Lambda$ maps homeomorphically onto a closed half-disc, with the boundary mapped onto the diameter, a real analytic arc.
\end{lemma}

\begin{proof}
The first claim holds because $q$ and $\Log|_C$ are local diffeomorphisms on $C\setminus S$. Let $p\in S$ with $q(p)=b$, and $g$, $u=u_1+iu_2$ as above, so that $\operatorname{Fix}_C(g)=\{u_2=0\}$ near $p$. Since $\Log(u)=\Log(\bar u)$, each component of $\Log$ is even in $u_2$, and since $\Log$ is real analytic its expansion in $u_2$ contains only even powers, so $\Log(u)=\Psi(u_1,u_2^2)$ for a real-analytic map $\Psi$; in the smooth category this is Whitney's theorem on even functions [Wh]. Here $\gamma$, $a_1$ and the finite set $F'$ of the previous lemmas are transported from $C'$ to $C$ by the monomial automorphism $\Phi$, which preserves all the properties used. In the coordinates $(u_1,v)=(u_1,u_2^2)$, $v\ge0$, on the half-disc $\{u_2\ge0\}$, $\Lambda$ is the restriction of $\Psi$, so it suffices that $\det D\Psi(u_1,0)\neq0$. The Jacobian of $\Log$ in the coordinate $u$ is $2u_2\det D\Psi(u_1,u_2^2)$, and by (\ref{L:e3.1}) it is $-|a_1|^2\operatorname{Im}\gamma$. At points of $S$ where $a_1\neq0$ and $d\gamma\neq0$, $\gamma$ is real on $\{u_2=0\}$ with real nonzero derivative there, so $\operatorname{Im}\gamma$ vanishes to order exactly one in $u_2$, and $\det D\Psi(u_1,0)\neq0$. The excluded points form a finite set, whose image is $F_Q$.
\end{proof}

\begin{lemma}\label{L:lem:one-branch}
For all $y\in\Lambda(\partial Q)$ outside a countable set, the fibre $\Lambda^{-1}(y)$ contains exactly one point of $\partial Q$, and it is not in $F_Q$.
\end{lemma}

\begin{proof}
Let $b_1\neq b_2$ in $\partial Q\setminus F_Q$ with $\Lambda(b_1)=\Lambda(b_2)=y$, and let $\sigma_1,\sigma_2$ be boundary arcs through them as in Lemma \ref{L:lem:L-local}. Suppose that $\Lambda(\sigma_1)$ and $\Lambda(\sigma_2)$ have infinitely many common points near $y$. Lift the arcs to arcs $\tilde\sigma_i\subset S$ through points $p_i$ over $b_i$. By the proof of Proposition \ref{L:prop:S-fixed}, after shrinking, $\tilde\sigma_i\subset e^{i\theta_i}(0,\infty)^2$. For a common image point $x$, the points of $\tilde\sigma_1$ and $\tilde\sigma_2$ over $x$ are $e^{i\theta_1}e^x$ and $e^{i\theta_2}e^x$, which differ by the translation $\tau=e^{i(\theta_2-\theta_1)}$, independent of $x$. So $\tau C\cap C$ is infinite, $\tau\in\Gamma_C^+$, and $q(p_2)=q(\tau p_1)=q(p_1)$, a contradiction. So the images of distinct boundary arcs meet in discrete sets. Covering $\partial Q\setminus F_Q$ by countably many arcs, and adding the finite set $\Lambda(F_Q\cup q(F))$, gives the countable exceptional set.
\end{proof}

\begin{proposition}\label{L:prop:single-orbit}
One has $|\Gamma_C|=\mu$, and every fibre of $\Log|_C$ over $A^\circ$ is a single free orbit of $\Gamma_C$.
\end{proposition}

\begin{proof}
For $x\in U$, the fibre lies in $C\setminus S$, where $\Gamma_C$ acts freely, so $N(x)=|\Gamma_C|\cdot\#\Lambda^{-1}(x)$. By Lemma \ref{L:lem:eq-first}, $\#\Lambda^{-1}$ takes only the values $0$ and $\mu/|\Gamma_C|$ on $U$. The set $B=\Lambda(\partial Q)$ is uncountable, since $\partial Q\setminus F_Q$ is a nonempty union of arcs on which $\Lambda$ is locally injective. Choose $y\in B$ outside the exceptional sets of Lemma~\ref{L:lem:one-branch}, and let $b$ be the unique boundary point over $y$. The other points of $\Lambda^{-1}(y)$ are interior, finitely many, and $\Lambda$ is a local diffeomorphism near each of them. Near $b$, $\Lambda$ maps a neighbourhood onto a closed half-disc bounded by an arc of $B$ through $y$. By properness, for $y'\in U$ near $y$ the number $\#\Lambda^{-1}(y')$ equals the number of interior preimages of $y$, plus one on one side of the arc and plus zero on the other. The two values differ by one and lie in $\{0,\mu/|\Gamma_C|\}$, so $\mu/|\Gamma_C|=1$.
\end{proof}

\begin{theorem}\label{L:thm:curve-equality}
Let $f$ be irreducible with smooth zero set $C$ and two-dimensional Newton polygon $\Delta$, not centrally symmetric. Then $\Area(\A(C))=\pi^2/\MV(\Delta,-\Delta)$ if and only if $\Delta$ is a triangle and $f$ is a trinomial.
\end{theorem}

\begin{proof}
The converse is Proposition \ref{L:prop:trinomial}. Assume equality and let $H=\Gamma_C^+$. By Propositions \ref{L:prop:S-fixed} and \ref{L:prop:single-orbit}, $|H|=\mu/2$. Let $P$ be as in Lemma \ref{L:lem:subgroup}, with kernel $H$, and $D=\phi_P(C)$. Since $\phi_P$ is finite and $C$ is $H$-invariant, $D$ is a closed irreducible curve with $C=\phi_P^{-1}(D)$, and $D$ is smooth because $\phi_P$ is a local biholomorphism. Let $h$ be an irreducible Laurent polynomial generating the ideal of $D$, so that $dh\neq0$ on $D$. Every irreducible factor of $h\circ\phi_P$ vanishes on a curve contained in $C$, hence is associate to $f$, so $h\circ\phi_P=uf^e$ with $u$ a unit, and $e=1$ because $d(h\circ\phi_P)=\phi_P^*dh\neq0$ on $C$. So $\Delta$ is a translate of $P^{\mathsf T}\Newt(h)$, and $\Area(\Newt(h))=\Area(\Delta)/|\det P|=2\Area(\Delta)/\mu$. By Section \ref{L:sec:mv}, $2\Area(\Delta)<\mu\le4\Area(\Delta)$, the first inequality being strict because $\Delta$ is not centrally symmetric. So $\frac12\le\Area(\Newt(h))<1$, and since twice the area of a lattice polygon is an integer, $\Area(\Newt(h))=\frac12$ and $\mu=4\Area(\Delta)$, so that $\Delta$ is a triangle. By Pick's formula a lattice polygon of area $\frac12$ has only its three vertices as lattice points, so $h$ is a trinomial, and so is $f=u^{-1}h\circ\phi_P$, with support the vertices of $\Delta$.
\end{proof}

\begin{corollary}\label{L:cor:reduced-equality}
Let $f$ be reduced with smooth zero set and Newton polygon $\Delta$, two-dimensional and not centrally symmetric. Then $\Area(\A(\{f=0\}))=\pi^2/\MV(\Delta,-\Delta)$ if and only if $f$ is a monomial times a trinomial whose Newton polygon is the triangle $\Delta$.
\end{corollary}

\begin{proof}
Write $f=z^\beta\prod_if_i$ with $f_i$ irreducible, pairwise non-associated, not monomials, and $\Delta_i=\Newt(f_i)$, so that $\Delta=\beta+\sum_i\Delta_i$. A sum of centrally symmetric polygons is centrally symmetric, so some $\Delta_{i_0}$ is not, and it is two-dimensional. The curves $\{f_i=0\}$ are smooth and pairwise disjoint, since a common point would be singular. With $Q=\sum_{i\neq i_0}\Delta_i$, multilinearity and nonnegativity of mixed areas give $\MV(\Delta,-\Delta)\ge\MV(\Delta_{i_0},-\Delta_{i_0})+\MV(\Delta_{i_0},-Q)$. Hence, by Theorem \ref{L:thm:curves} for $f_{i_0}$,
$$\Area(\A(\{f=0\}))\ge\Area(\A(\{f_{i_0}=0\}))\ge\frac{\pi^2}{\MV(\Delta_{i_0},-\Delta_{i_0})}\ge\frac{\pi^2}{\MV(\Delta,-\Delta)} .$$
If $Q$ is not a point, it contains a segment, and $\MV(\Delta_{i_0},-Q)>0$, so the last inequality is strict. Under equality, $Q$ is a point, $f=z^{\beta}f_{i_0}$ up to a constant, and Theorem \ref{L:thm:curve-equality} applies.
\end{proof}

This completes the proof of Theorem B.

\section{The parity argument}\label{L:s6}\label{L:sec:parity}

\begin{lemma}\label{L:lem:translate}
Let $X$ and $Y$ be nondegenerate smooth complete intersections in $(\C^*)^n$, of dimensions $a$ and $b$ with $a+b=n$, defined by polynomials with Newton polytopes $P_1,\dots,P_{n-a}$ and $R_1,\dots,R_{n-b}$. There is a proper Zariski closed subset $Z\subset(\C^*)^n$ such that, for $t\notin Z$, the intersection $X\cap tY$ is transverse and has exactly $\MV(P_1,\dots,P_{n-a},R_1,\dots,R_{n-b})$ points.
\end{lemma}

\begin{proof}
The variety $tY$ is defined by the polynomials $h(t^{-1}w)$, $h$ running over the equations of $Y$. They have the same Newton polytopes, and $\operatorname{in}_w(h(t^{-1}\cdot))=(\operatorname{in}_wh)(t^{-1}\cdot)$, because the substitution multiplies each coefficient by a nonzero constant. So the initial system in the direction $w$ defines $X_w\cap tY_w$, where $X_w$ and $Y_w$ are defined by the initial forms of the equations of $X$ and $Y$. Only finitely many initial systems occur, one for each relatively open cone of the common refinement of the normal fans of the polytopes involved, and each cone contains a primitive integer vector. For primitive $w\in\Z^n$, every initial form $g$ satisfies $g(s^wz)=s^cg(z)$ for $s\in\C^*$ and some integer $c$, where $s^w=(s^{w_1},\dots,s^{w_n})$. So $X_w$ and $Y_w$ are invariant under the one-parameter subgroup $G_w=\{s^w\}$. By nondegeneracy each of them is empty or of dimension $a$, respectively $b$. The set $B_w$ of $t$ with $X_w\cap tY_w\neq\emptyset$ is the image of $X_w\times Y_w$ under $(x,y)\mapsto xy^{-1}$. That morphism is constant on the one-dimensional orbits $\{(gx,gy):g\in G_w\}$, so $B_w$ is constructible of dimension at most $n-1$. For $t$ outside the closures of the $B_w$, the second part of Bernstein's theorem [Be, Theorem B], see also [Kh], gives exactly the mixed volume as number of common zeros counted with multiplicity. By Kleiman's transversality theorem for the transitive action of $(\C^*)^n$ on itself [Kl, Theorem 2], $X\cap tY$ is transverse for $t$ in a dense Zariski open set. Let $Z$ be the union of the closures of the $B_w$ and the complement of that open set.
\end{proof}

\begin{proposition}\label{L:prop:parity}
Let $V$ be nondegenerate. For almost every $\theta\in Q_\pi$, $V\cap\tau_\theta(V)$ has exactly $M_+(V)$ points and $N_\pi(\theta)\equiv M_+(V)$ modulo $2$.
\end{proposition}

\begin{proof}
We have $\tau_\theta(V)=e^{2i\theta}\cdot\bar V$, where $\bar V$ is defined by the polynomials with conjugated coefficients. It has the same Newton polytopes, and it is nondegenerate because its initial systems are the conjugates of those of $V$. Apply Lemma \ref{L:lem:translate} to $X=V$, $Y=\bar V$. The unit torus $\{|t_j|=1\}$ is Zariski dense in $(\C^*)^{2k}$, since a Laurent polynomial $\sum a_\beta t^\beta$ vanishing on it has vanishing Fourier coefficients $a_\beta$. So the intersection of $Z$ with the unit torus lies in the zero set of a Laurent polynomial not vanishing identically on it, a proper real analytic subset of a connected real analytic manifold, hence a null set. Its preimage under $\theta\mapsto e^{2i\theta}$ is null. For the other $\theta$, $V\cap\tau_\theta(V)$ has $M_+(V)$ points, and the involution $\tau_\theta$ maps this finite set to itself. Its orbits have one or two elements, so the number $N_\pi(\theta)$ of its fixed points has the parity of $M_+(V)$.
\end{proof}

\begin{proof}[Proof of Theorem C]
By Proposition \ref{L:prop:parity}, $N_\pi$ is odd, hence at least one, almost everywhere on $Q_\pi$, so $\mathcal M(V)=\int_{Q_\pi}N_\pi\ge\pi^{2k}$. Theorem A gives the volume bound.
\end{proof}

\begin{corollary}\label{L:cor:Mone}
If $V$ is nondegenerate and $M_+(V)=1$, then $N_\pi=1$ almost everywhere and $\mathcal M(V)=\pi^{2k}$.
\end{corollary}

\begin{proof}
By Propositions \ref{L:prop:parity} and \ref{L:prop:argbound}, $N_\pi$ is odd and at most one almost everywhere.
\end{proof}

For $k$-planes in general position, the equality $\mathcal M(V)=\pi^{2k}$ is due to Nisse and Passare [NP, Theorem 1.2 and Corollary 5.5] and, independently, to Johansson [Jo]; see also [Mik, Proposition 8], and [MN2, Theorem 7.1] for the injectivity of the argument map on generic linear spaces. Corollary \ref{L:cor:Mone} gives a different proof, through the parity argument, which applies to every nondegenerate complete intersection with $M_+(V)=1$.

\begin{example}\label{L:ex:planes}
For a $k$-plane $V$ in general position, given by $k$ affine linear equations, all $\Delta_j$ equal the standard simplex $S$ of $\R^{2k}$. Then $M_+=(2k)!\Vol(S)=1$, and Corollary \ref{L:cor:Mone} gives $\mathcal M(V)=\pi^{2k}$. Moreover $M_-=\MV(S^{(k)},(-S)^{(k)})=\binom{2k}k$, by the formula $\Vol(\lambda S-\mu S)=\frac1{n!}\sum_a\binom na^2\lambda^a\mu^{n-a}$ for the standard simplex of $\R^n$ (see [RS, Section 4]; the formula can also be checked directly by decomposing $\lambda S-\mu S$ into products of simplices). Theorem C gives $\Vol(\A(V))\ge\pi^{2k}/\binom{2k}k$, a bound first obtained by Mikhalkin [Mik, (3.1)]. Every plane in general position, defined by equations with all coefficients nonzero, is nondegenerate by Theorem~\ref{K2:thm:nondeg}.
\end{example}

For curves, $M_+$ is odd exactly when $\Delta$ has an odd number of boundary lattice points, and such polygons are not centrally symmetric (Section \ref{L:sec:mv}). So for curves Theorem C is contained in Theorem B. The triangle $2T_0=\conv(0,2e_1,2e_2)$, which has six boundary lattice points, shows that Theorem B is strictly stronger.

\section{Fibrations and towers}\label{L:s7}\label{L:sec:fibrations}

Write points of $(\C^*)^{2a}\times(\C^*)^{2b}$ as $(z',z'')$, and let $\Pi(z',z'')=z'$.

\begin{theorem}\label{L:thm:fibration}
Let $f_1,\dots,f_a$ be Laurent polynomials in $z'$ and $f_{a+1},\dots,f_{a+b}$ Laurent polynomials in $(z',z'')$. Let $W$ be the zero set of $f_1,\dots,f_a$ and $V$ that of all of them, and assume that $W$ and $V$ are smooth complete intersections of dimensions $a$ and $a+b$. Let $B\subset W$ be an algebraic subset containing no component of $W$, such that $\Pi^{-1}(B)\cap V$ contains no component of $V$, and such that for $p\in W\setminus B$ the fibre $V_p=\{q:f_{a+i}(p,q)=0,\ 1\le i\le b\}$ is a smooth complete intersection of dimension $b$ in $(\C^*)^{2b}$. Then
$$\mathcal M(V)=\int_{W\setminus B}\mathcal M(V_p)\,d\mu_W(p),$$
where $\mu_W=J_{\Arg}\,d\mathrm{vol}$ is the coamoeba measure of $W$, of total mass $\mathcal M(W)$. In particular, if $\mathcal M(V_p)\ge\pi^{2b}$ for $\mu_W$-almost every $p$, then $\mathcal M(V)\ge\pi^{2b}\mathcal M(W)$.
\end{theorem}

\begin{proof}
Let $R\subset V$ be the set of points $(p,q)$ with $p\in W\setminus B$. Its complement is the proper analytic subset $\Pi^{-1}(B)\cap V$ of each component, which has measure zero for the smooth measure $J_{\Arg}d\mathrm{vol}$ on $V$. At $(p,q)\in R$, the differentials $d_qf_{a+1},\dots,d_qf_{a+b}$ are independent, because $V_p$ is a smooth complete intersection. By the implicit function theorem on the manifold $W\times(\C^*)^{2b}$, the set $V$ is near $(p,q)$ a submanifold of dimension $a+b$, and $\Pi|_V$ is a submersion onto a neighbourhood of $p$ in $W$ whose fibres are open subsets of the $V_{p'}$. Choose real coordinates $(s,t)$ on $V$ near $(p,q)$, with $s$ coordinates on $W$ and $\Pi(s,t)=s$. The argument form of $V$ is $\omega_V=\Arg^*(d\theta_1\wedge\dots\wedge d\theta_{2a+2b})$. Since $\theta_1,\dots,\theta_{2a}$ are functions of $z'$, we have $\omega_V=\Pi^*\omega_W\wedge\omega''$, where $\omega_W$ is the argument form of $W$ and $\omega''=d\theta_{2a+1}\wedge\dots\wedge d\theta_{2a+2b}$. Write $\Pi^*\omega_W=A(s)\,ds$ and $\omega''=\mathcal B(s,t)\,dt+\eta$, where every term of $\eta$ contains some $ds_i$. Then $\omega_V=A(s)\mathcal B(s,t)\,ds\wedge dt$, and $\mathcal B(s,\cdot)\,dt$ is the argument form of the fibre $V_p$ in the coordinates $t$. Hence $|\omega_V|=|A(s)|\,|\mathcal B(s,t)|\,ds\,dt$. Cover $R$ by countably many such charts $U_1,U_2,\dots$ and let $R_1=U_1$, $R_n=U_n\setminus(U_1\cup\dots\cup U_{n-1})$, a measurable partition of $R$. On each chart, Tonelli's theorem in the coordinates $(s,t)$ gives $\int_{R_n}|\omega_V|=\int\big(\int\mathbf 1_{R_n}|\mathcal B(s,t)|\,dt\big)|A(s)|\,ds$. Summing over $n$, and using that for fixed $p$ the sets $R_n\cap V_p$ partition $V_p$, gives $\int_R|\omega_V|=\int_{W\setminus B}\mathcal M(V_p)\,d\mu_W(p)$. Finally, $B$ contains no component of $W$, so it is a proper analytic subset of each component of $W$, hence $\mu_W$-null. Without this hypothesis the formula would still hold, but $B$ could carry positive $\mu_W$-mass, for instance if $W$ has a component over which $V$ is empty, and the last assertion would fail.
\end{proof}

\begin{corollary}\label{L:cor:class}
Let $\mathcal C$ be the class of smooth complete intersections $V$ of dimension $k$ in $(\C^*)^{2k}$, $k\ge1$, with $\mathcal M(V)\ge\pi^{2k}$. It contains the nondegenerate $V$ with $M_+(V)$ odd, and the smooth reduced curves whose Newton polygons are not centrally symmetric. It is closed under products, under pullbacks by isogenies and multiplicative translations, and under fibrations as in Theorem \ref{L:thm:fibration} with $W\in\mathcal C$ and $V_p\in\mathcal C$ for $\mu_W$-almost every $p$. Every $V\in\mathcal C$ satisfies $\Vol(\A(V))\ge\pi^{2k}/M_-(V)$.
\end{corollary}

\begin{proof}
The first two assertions are Theorem C and Proposition \ref{L:prop:rays}. Closure under products follows from Lemma \ref{L:lem:product}, under isogenies and translations from Lemma \ref{L:lem:isogeny}, and under fibrations from Theorem \ref{L:thm:fibration}. The last assertion is Theorem A.
\end{proof}

We now turn to towers (Definition \ref{L:def:tower}). The requirement that $f_j(p,\cdot)$ never vanishes identically guarantees that $\Pi^{-1}(B_j)\cap V_j$ has dimension at most $j-1$. We write $\mu_j=\MV(\rho_j\Delta_j,-\rho_j\Delta_j)$ and, for $p\in V_{j-1}$, $C_p=\{q\in(\C^*)^2:f_j(p,q)=0\}$, so that $V_j=\{(p,q):p\in V_{j-1},\ q\in C_p\}$.

\begin{lemma}\label{L:lem:tower-mv}
For a tower, $M_-(V)=\prod_{j=1}^k\mu_j$ and $M_+(V)=\prod_{j=1}^k2\Area(\rho_j\Delta_j)$.
\end{lemma}

\begin{proof}
Reorder the arguments of $M_-(V)$ as $(\Delta_1,-\Delta_1,\Delta_2,-\Delta_2,\dots)$. Since $\pm\Delta_1\subset\R^2\times\{0\}$, Lemma \ref{L:lem:block} with $a=2$ gives $M_-(V)=\MV_2(\Delta_1,-\Delta_1)\cdot\MV_{2k-2}(\varrho\Delta_2,-\varrho\Delta_2,\dots)$, where $\varrho$ forgets the first two coordinates. Now $\pm\varrho\Delta_2\subset\R^2\times\{0\}\subset\R^{2k-2}$, with $\varrho\Delta_2=\rho_2\Delta_2$, and the claim follows by induction. The same computation applies to $M_+$.
\end{proof}

\begin{theorem}\label{L:thm:towers}
Let $V$ be a tower whose fibre polygons $\rho_j\Delta_j$ are not centrally symmetric. Then for every $j\ge2$
$$\Vol(\A(V_j))\ \ge\ \frac{\pi^2}{\mu_j}\Vol(\A(V_{j-1})),$$
and $\Vol(\A(V_1))\ge\pi^2/\mu_1$. Hence $\Vol(\A(V))\ge\pi^{2k}/M_-(V)$. If the tower is smooth, then $\mathcal M(V_j)\ge\pi^2\mathcal M(V_{j-1})$ and $\mathcal M(V)\ge\pi^{2k}$.
\end{theorem}

\begin{proof}
The set $\A(V_j)$ is closed, since $\Log|_{V_j}$ is proper. For $x'\in\R^{2j-2}$ its slice $\A_{x'}=\{x''\in\R^2:(x',x'')\in\A(V_j)\}$ is the union of the amoebas $\A(C_p)$ over the points $p\in V_{j-1}$ with $\Log p=x'$, and it is empty if $x'\notin\A(V_{j-1})$. By Fubini's theorem, $\Vol(\A(V_j))=\int_{\A(V_{j-1})}\Area(\A_{x'})\,dx'$. The set $B_j$ is an algebraic set of complex dimension at most $j-2$, a finite union of smooth real manifolds of dimension at most $2j-4$, so $\Log(B_j)$ is null in $\R^{2j-2}$. For $x'\in\A(V_{j-1})\setminus\Log(B_j)$, choose $p$ over $x'$, necessarily outside $B_j$. Then $f_j(p,\cdot)$ is reduced with Newton polygon $\rho_j\Delta_j$, not centrally symmetric, and $C_p$ is smooth. By Theorem \ref{L:thm:curves}, $\Area(\A_{x'})\ge\Area(\A(C_p))\ge\pi^2/\mu_j$. Integrating gives the first inequality. The case $j=1$ is Theorem \ref{L:thm:curves}. By induction and Lemma \ref{L:lem:tower-mv}, $\Vol(\A(V))\ge\pi^{2k}/\prod_j\mu_j=\pi^{2k}/M_-(V)$. For smooth towers, apply Theorem \ref{L:thm:fibration} with $a=j-1$, $b=1$, $W=V_{j-1}$ and $B=B_j$. For $p\notin B_j$, $C_p$ is a smooth curve defined by a reduced polynomial, so $d_qf_j\neq0$ on $C_p$ by the Jacobian criterion, and $\mathcal M(C_p)\ge\pi^2$ by Proposition \ref{L:prop:rays}.
\end{proof}

Theorem D follows from Lemma \ref{L:lem:tower-mv}, Theorem \ref{L:thm:towers} and Corollary \ref{L:cor:class}. The slice argument gives the volume bound for towers without using the coamoeba. It also gives a more general statement: for a fibration as in Theorem \ref{L:thm:fibration} in which $B$ contains no component of $W$, $\Vol(\A(V))\ge\Vol(\A(W))\cdot\inf_{p\in W\setminus B}\Vol(\A(V_p))$, because $\Log(B)$ is then a null subset of $\R^{2a}$.

\section{Equality for towers, and a twisted product}\label{L:s8}\label{L:sec:tower-equality}

\begin{proof}[Proof of Theorem E]
Suppose first that $\Vol(\A(V))=\pi^{2k}/M_-(V)$. In the chain of inequalities of Theorem \ref{L:thm:towers} all inequalities become equalities, so $\Vol(\A(V_1))=\pi^2/\mu_1$, and for each $j\ge2$, $\Vol(\A(V_j))=(\pi^2/\mu_j)\Vol(\A(V_{j-1}))$. By Corollary \ref{L:cor:reduced-equality}, $\rho_1\Delta_1=\Delta_1$ is a triangle and $f_1$ is a monomial times a trinomial, so $V_1$ is a trinomial curve. Fix $j\ge2$. The equality $\int_{\A(V_{j-1})}\big(\Area(\A_{x'})-\pi^2/\mu_j\big)\,dx'=0$, with a nonnegative integrand outside a null set, shows that for almost every $x'\in\A(V_{j-1})$ and for every $p$ over $x'$ outside $B_j$, $\Area(\A(C_p))=\Area(\A_{x'})=\pi^2/\mu_j$. By Corollary \ref{L:cor:reduced-equality} applied to $f_j(p,\cdot)$, whose Newton polygon is $\rho_j\Delta_j$, the polygon $\rho_j\Delta_j$ is a triangle $T_j$ with vertices $\beta_0,\beta_1,\beta_2$, and $f_j(p,\cdot)$ has exactly three monomials, at the vertices. Let $c_i(z')$ be the sum of the terms of $f_j$ whose exponents project to $\beta_i$ under $\rho_j$, divided by $z''^{\beta_i}$, and let $c_\beta$ be defined similarly for the other lattice points $\beta$ of $T_j$. Then $c_\beta(p)=0$ for all $p\notin B_j$ lying over the full-measure set of good $x'$. To see that these points are dense in $V_{j-1}$, work on the regular part $V_{j-1}^{\mathrm{reg}}$, which is a connected manifold because $V_{j-1}$ is irreducible. The Jacobian $J_{\Log}$ of $\Log$ on it is real analytic, and it does not vanish identically, since otherwise, as the image under $\Log$ of the singular locus of $V_{j-1}$ is null by a count of dimensions, Sard's theorem would give $\Vol(\A(V_{j-1}))=0$, contradicting the lower bound of Theorem \ref{L:thm:towers}. So its zero set is nowhere dense, and outside it $\Log$ is a local diffeomorphism, under which preimages of null sets are null. Hence the points lying over the full-measure set of good $x'$ are dense. Hence $c_\beta$ vanishes on $V_{j-1}$, and $V_j$ is the variety of Definition \ref{L:def:twisted} built from $W=V_{j-1}$, $T_j$ and the $c_i$. The $c_i$ do not vanish identically on any component, since $c_i(p)\neq0$ for $p\notin B_j$.

The $c_i$ have no common zero on $V_{j-1}$, because $f_j(p,\cdot)=\sum_ic_i(p)z''^{\beta_i}$ on $V_{j-1}$ and the tower axioms forbid $f_j(p,\cdot)$ to vanish identically. It remains to check the condition on absolute values. For $p\notin B_j$, put $a(p)=(c_1(p)/c_0(p),c_2(p)/c_0(p))$ and let $P_j$ be the matrix with rows $\beta_1-\beta_0$ and $\beta_2-\beta_0$. As in the proof of Proposition \ref{L:prop:trinomial}, $C_p=\phi_{P_j}^{-1}(a(p)^{-1}L)$, where $a^{-1}L=\{y:(a_1y_1,a_2y_2)\in L\}$, so $\A(C_p)=P_j^{-1}(\A(L)-\Log a(p))$. Let $p,p'$ be two points outside $B_j$ over the same $x'$, with $x'$ as above. Both $\A(C_p)$ and $\A(C_{p'})$ are contained in $\A_{x'}$ and have the same area as $\A_{x'}$, which is finite, so their symmetric difference is null. They are translates of each other by $v=P_j^{-1}(\Log a(p')-\Log a(p))$. Let $K=\A(C_p)$. Then $\mathbf 1_K(x)=\mathbf 1_K(x-v)$ for almost every $x$. If $v\neq0$, decompose $\R^2$ into the translates $S+mv$, $m\in\Z$, of a strip $S$ of width $|v|$ transverse to $v$. Then $\Area(K\cap(S+mv))$ does not depend on $m$, and since $\Area(K)<\infty$ it vanishes, so $\Area(K)=0$, a contradiction. Hence $\Log a(p)=\Log a(p')$, which is the required condition.

Conversely, suppose that the tower $V$ is an iterated twisted product of trinomial curves along its filtration. Then $f_j=\sum_ic_i(z')z''^{\beta_i}$, fibre polygons the triangles $T_j$, $\mu_j=4\Area(T_j)$, and $B_j$ the set where some $c_i$ vanishes. For $x'\in\A(V_{j-1})$ outside the null set $\Log(B_j)$ and outside the null set where the condition of Definition \ref{L:def:twisted} fails, all the points $p$ over $x'$ give the same set $\A(C_p)=P_j^{-1}(\A(L)-\Log a(p))$. So $\A_{x'}=\A(C_p)$ has area $\pi^2/\mu_j$ by Proposition \ref{L:prop:trinomial}. By Fubini, $\Vol(\A(V_j))=(\pi^2/\mu_j)\Vol(\A(V_{j-1}))$, and by induction and Lemma \ref{L:lem:tower-mv}, $\Vol(\A(V))=\pi^{2k}/M_-(V)$.
\end{proof}

\begin{remark}
Products of trinomial curves are, up to isogenies and monomial changes of coordinates, twisted products in which the $c_i$ are monomials times constants. Conversely, if $c_i=\kappa_iz'^{\gamma_i}$, a substitution $z''\mapsto z'^{A}z''$, where $(z'^Az'')_l=z''_l\prod_mz_m'^{\,A_{ml}}$, transforms $c_iz''^{\beta_i}$ into $\kappa_iz'^{\gamma_i+A\beta_i}z''^{\beta_i}$. With $A$ a rational matrix solving $A(\beta_i-\beta_0)=-(\gamma_i-\gamma_0)$ for $i=1,2$, this substitution, realised by a monomial change of coordinates composed with an isogeny, makes $V_j$ a product of $V_{j-1}$ with a trinomial curve. The condition of Definition \ref{L:def:twisted} is satisfied by monomials, because $|z'^\gamma|=e^{\langle\gamma,x'\rangle}$ depends only on $x'$. It is also satisfied when $V_{j-1}$ is real up to a multiplicative translation, its amoeba fibres are orbits of the conjugation, and the $c_i$ have real coefficients in the corresponding coordinates. The next example is of this kind.
\end{remark}

\begin{proof}[Proof of Theorem F]
Let $f_1=1+z_1+z_2$ and $f_2=1+(2+z_1)z_3+z_4$. The variety $V$ is the graph of $(z_1,z_3)\mapsto(z_2,z_4)=(-1-z_1,-1-(2+z_1)z_3)$ over the open set of $(\C^*)^2$ where these values are nonzero, so it is a smooth irreducible surface. It is a smooth complete intersection, since $df_1=dz_1+dz_2$ and $df_2=z_3dz_1+(2+z_1)dz_3+dz_4$ are independent. It is a tower with $V_1=L$, $\rho_2\Delta_2=\conv(0,e_1,e_2)=T_0$, and $B_2=\{(-2,1)\}$, where $f_2(p,\cdot)=1+z_4$ is not identically zero. For $p\notin B_2$, $f_2(p,\cdot)$ is the trinomial $1+c\,z_3+z_4$ with $c=2+z_1(p)\neq0$. By Lemma \ref{L:lem:line}, for almost every $x'\in\A(L)$ the fibre of $\Log|_L$ over $x'$ consists of two conjugate points $p,\bar p$, and $|2+\bar z_1|=|2+z_1|$. So $V$ is a twisted trinomial extension of $L$, and by Theorem E and Lemma \ref{L:lem:tower-mv}, $\Vol(\A(V))=\pi^4/(\mu_1\mu_2)=\pi^4/4$ and $M_-(V)=4$. Explicitly, on $L$ one has $\operatorname{Re}z_1=\frac12(|1+z_1|^2-|z_1|^2-1)$, so $|2+z_1|^2=2e^{2x'_2}-e^{2x'_1}+2$. The slice of $\A(V)$ over $x'$ is the translate of $\A(L)$ by $(-\frac12\log(2e^{2x'_2}-e^{2x'_1}+2),0)$.

Say that a variety $X\subset(\C^*)^4$ has property $(\Pi)$ if there are two complementary rational planes $\Lambda_1,\Lambda_2\subset\Q^4$ such that, for $i=1,2$ and all $a,b\in\Lambda_i\cap\Z^4$, the monomials $z^a$ and $z^b$ are functionally dependent on $X$, that is, $(z^a,z^b)$ maps each component of $X$ into a curve. A product $C_1\times C_2$ has property $(\Pi)$ with the two coordinate planes. Property $(\Pi)$ is preserved by multiplicative translations, which multiply monomials by constants. It is preserved by pullbacks $\phi_P^{-1}$, with the planes $P^{\mathsf T}\Lambda_i$, because $z^a\circ\phi_P=z^{P^{\mathsf T}a}$. It is preserved by images $\phi_P(X)$, with the planes $(P^{\mathsf T})^{-1}\Lambda_i$: for $b\in\Z^4$, $w^b$ on $\phi_P(X)$ pulls back to $z^{P^{\mathsf T}b}$ on $X$, and $\phi_P:X\to\phi_P(X)$ is finite and surjective, so dependence on $X$ implies dependence on the image. The integer points of the new planes need not be images of integer points, but they have multiples which are, and if $z^{Nc}$ and $z^{Nd}$ are functionally dependent for some integer $N\ge1$, so are $z^c$ and $z^d$, because $(u,v)\mapsto(u^N,v^N)$ is finite. It suffices therefore to show that $V$ does not have property $(\Pi)$. If two holomorphic functions $g,h$ on a connected complex surface map it into a curve, then $dg\wedge dh=0$: otherwise $(g,h)$ would be a local biholomorphism at some point, and its image would contain an open set. So property $(\Pi)$ gives two independent vectors in $\Lambda_1$, and two in $\Lambda_2$, satisfying $d\log z^a\wedge d\log z^b=0$ on $V$. We show that for linearly independent $a,b\in\Q^4$, this condition forces $a,b\in\Pi_0=\{v\in\Q^4:v_3=v_4=0\}$. Then $\Lambda_1=\Lambda_2=\Pi_0$, which contradicts complementarity.

Use $(z_1,z_3)$ as coordinates on $V$. Then
$$d\log z_1=\frac{dz_1}{z_1},\quad d\log z_2=\frac{dz_1}{1+z_1},\quad d\log z_3=\frac{dz_3}{z_3},\quad d\log z_4=\frac{z_3\,dz_1+(2+z_1)\,dz_3}{D},$$
with $D=1+(2+z_1)z_3$. For $v\in\Q^4$ put $\lambda_v=\sum_jv_j\,d\log z_j$ and $E_v=v_1/z_1+v_2/(1+z_1)$. Suppose that $a,b$ are independent with $\lambda_a\wedge\lambda_b=0$, and that their span is not contained in $\Pi_0$. If the span is contained in $\{v_4=0\}$, it contains a nonzero $e$ with $e_3=e_4=0$ and a vector $g$ with $g_4=0$, $g_3\neq0$. Then $\lambda_e\wedge\lambda_g=E_e\,g_3\,dz_1\wedge dz_3/z_3=0$, so $E_e\equiv0$, which forces $e_1=e_2=0$ and $e=0$, a contradiction. Otherwise the span contains a vector $g$ with $g_4\neq0$ and a nonzero vector $e$ with $e_4=0$. The coefficient of $dz_1\wedge dz_3$ in $\lambda_e\wedge\lambda_g$ is
$$E_e\Big(\frac{g_3}{z_3}+\frac{g_4(2+z_1)}{D}\Big)-\frac{e_3}{z_3}\Big(E_g+\frac{g_4z_3}{D}\Big)=\frac{E_eg_3-e_3E_g}{z_3}+\frac{g_4\big((2+z_1)E_e-e_3\big)}{D}.$$
For fixed generic $z_1$, the functions $1/z_3$ and $1/D$ of $z_3$ are linearly independent, having poles at $0$ and at $-1/(2+z_1)$. So both coefficients vanish, and since $g_4\neq0$, $(2+z_1)(e_1/z_1+e_2/(1+z_1))=e_3$ identically in $z_1$. The residues at $z_1=0$ and $z_1=-1$ give $2e_1=0$ and $e_2=0$, and then $e_3=0$, so $e=0$, a contradiction. Hence the span of $a,b$ is contained in $\Pi_0$.
\end{proof}

\begin{remark}
The map $(z_1,z_2,z_3,z_4)\mapsto((z_1,z_2),((2+z_1)z_3,z_4))$ is an algebraic isomorphism from the open dense subset $\{z_1\neq-2\}$ of the surface of Theorem F onto $(L\setminus\{(-2,1)\})\times L$; the complement is the curve $\{z_1=-2,\ z_2=1,\ z_4=-1\}$. So the surface is, up to this curve, abstractly a product of two curves, and Theorem F concerns products in the toric sense, which is the sense relevant for amoebas: the isomorphism is not monomial. Over the point $x'=(\log2,0)$ of $\A(L)$, whose only preimage in $L$ is $(-2,1)$, the slice of $\A(V)$ is the line $\{x_4=0\}$, so the description of the amoeba as a sheared product holds outside this null set. Note also that $M_+(V)=1$, so that, if $V$ is nondegenerate, Corollary \ref{L:cor:Mone} gives $\mathcal M(V)=\pi^4$ independently of the tower structure. Its amoeba is the image of $\A(L)\times\A(L)$ under the volume preserving shear $(x',x'')\mapsto(x',x''+(h(x'),0))$, where $h(x')=-\frac12\log(2e^{2x'_2}-e^{2x'_1}+2)$, a nonlinear function.
\end{remark}

\section{Examples and open questions}\label{L:s9}\label{L:sec:examples}

The class of towers is not contained in the class covered by Theorem C, and the planes covered by Theorem C are not towers in their given coordinates. A $k$-plane with generic coefficients satisfies Theorem C, but for $k\ge2$ it is not a tower in the given coordinates, because its projection to $(z_1,z_2)$ is dominant, so that no nonzero Laurent polynomial in $z_1,z_2$ vanishes on it. We do not know whether it becomes a tower after a monomial change of coordinates. On the other hand, the tower $\{f_1(z_1,z_2)=0,\ f_2(z_1,\dots,z_4)=0\}$ with both fibre polygons equal to $2T_0$ satisfies Theorem D, while $M_+=4\cdot4=16$ is even.

The class $\mathcal C$ of Corollary \ref{L:cor:class} contains all these examples, their products, isogeny pullbacks and fibrations. We do not know a description of the tuples $(\Delta_1,\dots,\Delta_k)$ for which every nondegenerate $V$ belongs to $\mathcal C$. For curves, Theorem B gives the answer: the polygon must not be centrally symmetric. In higher dimension, the products $C_\varepsilon\times L^{k-1}$, with $C_\varepsilon$ as in Proposition \ref{L:prop:zonotope}, have $\mathcal M\to0$ and $\Vol(\A)\to0$ while all polytopes but one are triangles. So a condition on the tuple cannot only ask that some polytope be not centrally symmetric: a single centrally symmetric polytope, here the square of the first factor, already allows the bound to fail.

\begin{conjecture}\label{L:conj:equality}
Let $V$ be a smooth complete intersection of dimension $k$ in $(\C^*)^{2k}$ with $\mathcal M(V)\ge\pi^{2k}$. If $\Vol(\A(V))=\pi^{2k}/M_-(V)$, then, after a monomial change of coordinates, $V$ is the pullback by an isogeny of an iterated twisted product of trinomial curves.
\end{conjecture}

Theorem E proves this for towers, and Theorem F shows that twisted products cannot be avoided. By Theorem~\ref{T:thm:extremal}(b), equality forces $N=M_-(V)$ on the regular values of the amoeba and $\mathcal M(V)=\pi^{2k}$. When $V$ is nondegenerate and $M_+(V)$ is odd, it also forces $N_\pi=1$ almost everywhere on $Q_\pi$. Without the parity hypothesis this fails: the extremal trinomial curve $1+z^2+w^2=0$ is the pullback of a line by $(z,w)\mapsto(z^2,w^2)$, and its rolled count $N_\pi$ takes only the values $0$ and $4$. For curves, the classification used the condition $N=\MV(\Delta,-\Delta)$ together with the monotonicity of the arguments along rays. For $k$-planes in general position, $\mathcal M(V)=\pi^{2k}$ by Corollary \ref{L:cor:Mone}, and equality would require $N=\binom{2k}k$ almost everywhere. Conjecture \ref{L:conj:equality} predicts that it cannot, and this is proved in Corollary~\ref{K2:cor:strict}: for every $k\ge2$, no $k$-plane in general position attains the bound, because near a nondegenerate corner the regular fibres of $\Log|_V$ have at most $2^k<\binom{2k}k$ points.

\section{$k$-planes and corners: the equality case}\label{K:s1}

We return to the setting of Sections~\ref{L:sec:log} to \ref{L:sec:examples}: $V\subset(\C^*)^{2k}$ is a smooth complete intersection of dimension $k$, defined by $f_1,\dots,f_k$ with Newton polytopes $\Delta_1,\dots,\Delta_k$, $M_\pm(V)=\MV(\Delta_1,\dots,\Delta_k,\pm\Delta_1,\dots,\pm\Delta_k)$, $N(x)$ is the number of points of $V$ over $x\in\R^{2k}$, and $\mathcal M(V)=\int N$ is the volume of the coamoeba counted with multiplicity. Theorem A says that $N\le M_-(V)$ almost everywhere and $\Vol(\A(V))\ge\mathcal M(V)/M_-(V)$. Theorem C says that $\mathcal M(V)\ge\pi^{2k}$ when $V$ is nondegenerate and $M_+(V)$ is odd. Equality $\Vol(\A(V))=\pi^{2k}/M_-(V)$ in the resulting inequality forces two conditions: $N=M_-(V)$ almost everywhere on $\A(V)$, and $\mathcal M(V)=\pi^{2k}$. Conjecture~\ref{L:conj:equality} predicts that only isogeny pullbacks of iterated twisted products of trinomial curves attain the bound. Theorem E proves it for towers. The first class beyond towers is that of $k$-planes, which we treat in this section and in Sections~\ref{K:s2} to \ref{K2:s3}: the present section introduces corners and states the corner criterion, Sections~\ref{K:s2} to \ref{K:s4} prove it and establish nondegenerate corners for $k=2$, Section~\ref{K:s6} discusses what the method says in general, and Sections~\ref{K2:s1} to \ref{K2:s3} remove the restriction on $k$.

A $k$-plane is an affine subspace $V\subset\C^{2k}$ of dimension $k$, and we also write $V$ for $V\cap(\C^*)^{2k}$. Its closure $\bar V\simeq\PP^k$ in $\PP^{2k}$ meets the $2k$ coordinate hyperplanes and the hyperplane at infinity in $2k+1$ hyperplanes $\lambda_1,\dots,\lambda_{2k},\lambda_0$ of $\PP^k$. The plane is in general position if any $k+1$ of them have empty intersection. We may define $V$ by $k$ affine linear equations all of whose coefficients are nonzero. Indeed, generic linear combinations of a basis of equations have this property, because no coordinate is constant on $V$ and $V$ does not contain the origin, as follows from general position. Then every $\Delta_j$ is the standard simplex $S$, so $M_+(V)=1$ and $M_-(V)=\binom{2k}k$ (Example~\ref{L:ex:planes}). Such a plane is nondegenerate by Theorem~\ref{K2:thm:nondeg}, and Corollary~\ref{L:cor:Mone} gives $\mathcal M(V)=\pi^{2k}$, and the bound reads $\Vol(\A(V))\ge\pi^{2k}/\binom{2k}k$. Equality would require $N=\binom{2k}k$ almost everywhere on the amoeba. For $k=2$ this means that almost every fibre has six points, the maximal number allowed by Bernstein's theorem.

For a set $M\subset\{1,\dots,2k\}$ of $k-1$ indices, let $O$ be the complementary set of $k+1$ indices, and write $z=(z_M,z_O)$. Let $L_M=\{z\in\bar V\cap\C^{2k}:z_M=0\}$, which by general position is an affine line in the coordinate subspace $\{z_M=0\}\simeq\C^{k+1}$ (Lemma \ref{K:lem:geometry}). Let $d$ be a direction vector of $L_M$, and for $j\in M$ let $u_j$ be the vector of the direction space $D$ of $V$ whose $M$-coordinates are $(u_j)_i=\delta_{ij}$, $i\in M$. It exists, and it is unique up to adding multiples of $d$ (Lemma \ref{K:lem:geometry}). For $h\in L_M\cap(\C^*)^O$ write $w/h=(w_l/h_l)_{l\in O}$ for $w\in\C^{2k}$.

\begin{definition}\label{K:def:CN}
The corner $(M,h)$ of $V$ is nondegenerate if there is $\theta\in(\R/2\pi\Z)^M$ such that the real $(k+1)\times(k+1)$ matrix with columns
$$\operatorname{Re}\frac dh,\quad\operatorname{Im}\frac dh,\quad\operatorname{Im}\frac{e^{i\theta_j}u_j}h\ \ (j\in M)$$
is invertible. This does not depend on the choice of the $u_j$ modulo $d$, nor on the choice of $d$. The plane $V$ has a nondegenerate corner if some $(M,h)$ is nondegenerate.
\end{definition}

\begin{theorem}\label{K:thm:main}
Let $V$ be a $k$-plane in general position in $(\C^*)^{2k}$, $k\ge2$, with a nondegenerate corner. Then there is a set $O_V\subset\A(V)$ of positive measure on which every regular fibre of $\Log|_V$ has at least one and at most $2^k$ points. Consequently, if $V$ is moreover nondegenerate, which is always the case by Theorem~\ref{K2:thm:nondeg},
$$\Vol(\A(V))\ \ge\ \frac{\pi^{2k}+\big(\binom{2k}k-2^k\big)|O_V|}{\binom{2k}k}\ >\ \frac{\pi^{2k}}{\binom{2k}k}.$$
\end{theorem}

\begin{theorem}\label{K:thm:CN}
Every $2$-plane in general position in $(\C^*)^4$ has a nondegenerate corner.
\end{theorem}

\begin{corollary}\label{K:cor:strict}
For every $2$-plane in general position in $(\C^*)^4$, the inequality $\Vol(\A(V))\ge\pi^4/6$ is strict. In particular no such plane attains the bound $\pi^4/6=\pi^{2k}/\binom{2k}k$ of Theorem C for $k=2$, in accordance with Conjecture~\ref{L:conj:equality}.
\end{corollary}

\begin{remark}
The bound $2^k$ is sharp. For a generic real $k$-plane, every regular fibre of $\Log|_V$ over the amoeba has exactly $2^k$ points, by Nisse and Passare [NP, Proposition 5.6]. So the corner count is the fibre count of real planes, and Theorem \ref{K:thm:main} says that near a nondegenerate corner every plane behaves at most like a real plane. Together with the identity $\int N=\pi^{2k}$, this gives the amoeba volume $\pi^{2k}/2^k$ of generic real planes [NP, Theorem 1.2], [MN2, Corollary 7.1], and it is consistent with the conjecture that $\pi^{2k}/2^k$ is the minimal volume of amoebas of $k$-planes, attained by real planes.
\end{remark}

\section{The geometry of a corner}\label{K:s2}\label{K:sec:geometry}

\begin{lemma}\label{K:lem:geometry}
Let $V$ be a $k$-plane in general position and $M$ a set of $k-1$ coordinate indices, with complement $O$. Then the following hold.
The coordinate functions $z_j$, $j\in M$, are affine functions on $V$ with linearly independent linear parts, and $L_M$ is an affine line. For $j\in M$ there is a vector $u_j\in D$ with $(u_j)_M=e_j$, unique modulo $d$, and every point of $V$ can be written uniquely as $h+\sum_{j\in M}t_ju_j$ with $h\in L_M$ and $t_j=z_j$.
The line $L_M$ meets each hyperplane $\{z_l=0\}$, $l\in O$, in exactly one point, and these points are distinct and finite. In a parameter $\rho\in\C$, $L_M=\{(\gamma_l(\rho-q_l))_{l\in O}\}$ with $\gamma_l\neq0$ and distinct $q_l$.
The projection $\pi_O:V\to\C^O$ is injective, and its image $V_O$ is an affine hyperplane $\{\sum_{l\in O}b_lz_l=e\}$ with $e\neq0$.
\end{lemma}

\begin{proof}
The closures of the hyperplanes $\{z_j=0\}\cap V$, $j\in M$, are $k-1$ of the hyperplanes $\lambda_j$ of $\bar V\simeq\PP^k$. Adding any two further ones gives $k+1$ hyperplanes with empty intersection, so the $k-1$ hyperplanes $\lambda_j$, $j\in M$, are in general position and meet in a projective line $\bar L_M$, not contained in $\lambda_0$, and the linear parts of the $z_j$ are independent. The affine part of $\bar L_M$ is $L_M$, and the linear parts of the $z_j$, $j\in M$, map $D$ onto $\C^M$ with kernel $\C d$. This gives the vectors $u_j$, unique modulo $d$, and the decomposition. Each $\lambda_l$, $l\in O$, meets $\bar L_M$ in one point, because $\bar L_M\not\subset\lambda_l$, otherwise $k$ hyperplanes $\lambda_j$, $j\in M$, and $\lambda_l$ together with any other would be concurrent. For the same reason the points for different $l$ are distinct, and different from the point $\bar L_M\cap\lambda_0$, which is the point at infinity. So $z_l|_{L_M}=\gamma_l(\rho-q_l)$ with distinct $q_l$. A nonzero vector $w\in D$ with $w_O=0$ would define a point at infinity of $\bar V$ lying on all $\lambda_l$, $l\in O$, which are $k+1$ hyperplanes. So $\pi_O$ is injective on $V$, and its image is an affine subspace of dimension $k$ in $\C^{k+1}$, a hyperplane $\{\sum b_lz_l=e\}$. If $e=0$, the origin of $\C^O$ would be in $V_O$, giving a point of $V$ on the $k+1$ hyperplanes $\lambda_l$, $l\in O$.
\end{proof}

We denote by $\Log_O$ the map $z\mapsto(\log|z_l|)_{l\in O}$.

\begin{lemma}\label{K:lem:line}
The map $\Log_O|_{L_M}$ is at most two-to-one, and it is an immersion at every point $h(\rho)$ with $\rho$ outside the line through the $q_l$, when these are collinear, and at every point when they are not.
\end{lemma}

This is the statement on amoebas of lines of [NP, Section 4.1] and [MN2, Remark 7.1], see also Lemma \ref{F:lem:lineN}; we include the short proof.

\begin{proof}
Two parameters $\rho\neq\rho'$ have the same image if and only if $|\rho-q_l|=|\rho'-q_l|$ for all $l\in O$, that is, if all $q_l$ lie on the perpendicular bisector of $[\rho,\rho']$. Since $k+1\ge3$, this forces the $q_l$ to be collinear, and then $\rho'$ is the reflection of $\rho$ in their line, so fibres have at most two points. The derivative of $\Log_O\circ h$ with respect to $\rho$ is given by the vector $(1/(\rho-q_l))_l$. The map is an immersion at $\rho$ unless $e^{i\phi}/(\rho-q_l)$ is real for some $\phi$ and all $l$, that is, unless $\rho$ and all $q_l$ lie on one line.
\end{proof}

\begin{lemma}\label{K:lem:regular}
The set of $h\in L_M$ at which $\Log_O|_{V_O}$ has rank $k+1$ at $\pi_O(h)$ is the complement of a proper real-analytic subset of $L_M$.
\end{lemma}

\begin{proof}
At $z\in V_O\cap(\C^*)^O$, the logarithmic tangent space is $L_z=\{w:\sum_lb_lz_lw_l=0\}$. The rank of $\operatorname{Re}|_{L_z}$ is less than $k+1$ if and only if a nonzero real vector $\eta$ annihilates $\operatorname{Re}(L_z)$. The real part of the complex linear form $w\mapsto\sum\eta_lw_l$ vanishes on the complex subspace $L_z$ if and only if the form itself vanishes there, that is, if and only if $\eta=\lambda(b_lz_l)_l$ for some $\lambda\in\C^*$. Then all $b_lz_l$ lie on the real line $\lambda^{-1}\R$, and since they add up to $e\neq0$, that line is $\R e$. So $z$ is critical exactly when $b_lz_l/e\in\R$ for all $l$. If this held on an open subset of $L_M$, the holomorphic functions $b_lz_l/e$ on $L_M$ would be real, hence constant, for all $l$. Some $b_l$ is nonzero, and the corresponding $z_l$ would be constant on $L_M$, while it vanishes at $\rho=q_l$ and not elsewhere.
\end{proof}

\section{The fibre count near a corner}\label{K:s3}\label{K:sec:corner}

\begin{lemma}\label{K:lem:trig}
Let $n\ge1$, $A_1,\dots,A_n,B_1,\dots,B_n\in\R^n$, and $F(\theta)=\sum_{j=1}^n(\cos\theta_j\,A_j-\sin\theta_j\,B_j)$ for $\theta\in(\R/2\pi\Z)^n$. If $c\in\R^n$ is a regular value of $F$, then $\#F^{-1}(c)\le2^n$.
\end{lemma}

\begin{proof}
Consider the polynomial system in the $2n$ complex unknowns $(C_j,S_j)$ given by the $n$ linear equations $\sum_j(C_jA_j-S_jB_j)=c$ and the $n$ quadratic equations $C_j^2+S_j^2=1$. The points $(\cos\theta_j,\sin\theta_j)$ for $\theta\in F^{-1}(c)$ are distinct real solutions. At such a point, the Jacobian matrix of the system is invertible. Indeed, a vector $(\dot C,\dot S)$ in its kernel satisfies $C_j\dot C_j+S_j\dot S_j=0$, so $(\dot C_j,\dot S_j)=\dot\theta_j(-S_j,C_j)$ for real or complex $\dot\theta_j$, and then the linear equations give $dF_\theta(\dot\theta)=0$, hence $\dot\theta=0$, the complexification of an injective real matrix being injective. So these solutions are isolated, and by B\'ezout's theorem the number of isolated solutions in $\C^{2n}$ of a system of $n$ equations of degree one and $n$ of degree two is at most $2^n$ [Fu2, Example 8.4.6].
\end{proof}

\begin{proposition}\label{K:prop:corner}
Let $(M,h_0)$ be a nondegenerate corner of $V$ at which $\Log_O|_{L_M}$ is an immersion. Then there is a set $O_V\subset\R^{2k}$ of positive measure such that every $x\in O_V$ is the image of a point of $V$, and every regular value $x\in O_V$ of $\Log|_V$ has at most $2^k$ preimages.
\end{proposition}

\begin{proof}
Let $y_0=\Log_Oh_0$. By Lemma \ref{K:lem:line}, the preimages of $y_0$ in $L_M$ are $h_0$ and possibly one more point $h_0'$, with the same image. Since a reflection of the parameter is an antiholomorphic map of $L_M$ preserving $\Log_O$, $\Log_O|_{L_M}$ is an immersion at $h_0'$ as well. Fix $j_0\in M$ and write $M=\{m_1,\dots,m_{k-1}\}$. For $r>0$, $s\in[1,2]^{M}$ and $c\in\R^O$, consider the point
$$x(r,s,c)=\big((\log(rs_j))_{j\in M},\ y_0+rc\big)\in\R^M\times\R^O .$$
By Lemma \ref{K:lem:geometry}, a point of $V$ is $z=h+\sum_jt_ju_j$ with $h\in L_M$, $t_j=z_j$. It lies over $x(r,s,c)$ if and only if $|t_j|=rs_j$ for $j\in M$ and $\Log_O(h+\sum_jt_ju_j)=y_0+rc$.

We first show that, for $c$ and $s$ in compact sets and $r$ small, all such $z$ have $h$ close to $h_0$ or to $h_0'$, at distance $O(r)$. As $r\to0$, $|t_j|\to0$ and $\Log_Oz$ stays bounded, so $z_O$ stays in a compact subset of $(\C^*)^O$. Any limit point of such points $z$ lies in $\bar V\cap\{z_M=0\}=L_M$ and over $y_0$, hence equals $h_0$ or $h_0'$. Since $|z_O-h_O|=O(r)$ and $\Log_O$ is Lipschitz near these points, $|\Log_Oh-y_0|=O(r)$. Since $\Log_O|_{L_M}$ is an immersion at $h_0$ and $h_0'$, with these two points as the only preimages of $y_0$, we get $\operatorname{dist}(h,\{h_0,h_0'\})=O(r)$.

Near $h_0$, write $h=h_0+r\omega d$ with $\omega\in\C$ bounded, and $t_j=rs_je^{i\theta_j}$. Then
$$\Psi_r(\omega,\theta)=\frac1r\Big(\Log_O\big(h_0+r\omega d+r\textstyle\sum_js_je^{i\theta_j}u_j\big)-y_0\Big)=\operatorname{Re}\Big(\omega\frac d{h_0}+\sum_js_je^{i\theta_j}\frac{u_j}{h_0}\Big)+O(r),$$
with the error and its derivatives in $(\omega,\theta)$ uniformly $O(r)$ on compact sets, and uniformly in $s\in[1,2]^M$. The points of $V$ near $h_0$ over $x(r,s,c)$ correspond bijectively to the solutions of $\Psi_r(\omega,\theta)=c$ with $\omega$ bounded. Let $\Psi_0^s(\omega,\theta)=\operatorname{Re}(\omega d/h_0+\sum_js_je^{i\theta_j}u_j/h_0)$, a map from $\C\times(\R/2\pi\Z)^M$ to $\R^O$, both of real dimension $k+1$. Its differential at $(\omega,\theta)$ maps $(\dot\omega,\dot\theta)$ to $\operatorname{Re}(\dot\omega d/h_0)+\sum_js_j\dot\theta_j\operatorname{Re}(ie^{i\theta_j}u_j/h_0)$. Its matrix, in the basis given by $\dot\omega=1$, $\dot\omega=i$ and the $\dot\theta_j$, has columns $\operatorname{Re}(d/h_0)$, $-\operatorname{Im}(d/h_0)$ and $-s_j\operatorname{Im}(e^{i\theta_j}u_j/h_0)$. Up to signs and the positive factors $s_j$, these are the columns of Definition \ref{K:def:CN}. So $\Psi_0^s$ is a local diffeomorphism at $(\omega,\theta^*)$ for every $\omega$, where $\theta^*$ is given by the nondegeneracy of the corner.

The equation $\Psi_0^s(\omega,\theta)=c$ splits. Let $T=\operatorname{Re}(\C d/h_0)$, the tangent plane of $\Log_O(L_M)$ at $y_0$, of dimension two because of the immersion, and let $P$ be the projection of $\R^O$ onto $\R^O/T\simeq\R^{k-1}$. Then $\Psi_0^s(\omega,\theta)=c$ if and only if $F_s(\theta)=Pc$, where $F_s(\theta)=\sum_js_jP\operatorname{Re}(e^{i\theta_j}u_j/h_0)$, and $\omega$ is then uniquely determined by $\theta$. The map $F_s$ has the form of Lemma \ref{K:lem:trig} with $n=k-1$, $A_j=s_jP\operatorname{Re}(u_j/h_0)$ and $B_j=s_jP\operatorname{Im}(u_j/h_0)$, and $\theta^*$ is a regular point of $F_{s}$ for all $s$. The same holds near $h_0'$, with a map $F'_s$, which may have no regular point. Choose $s^*$ and a value $c^*$ such that $Pc^*$ is a regular value of $F_{s^*}$ and of $F'_{s^*}$ and lies in the open set $F_{s^*}(W)$, where $W$ is a neighbourhood of $\theta^*$ consisting of regular points. Such values form a nonempty open set, by Sard's theorem and since the critical values of a map from a compact manifold form a closed set. Then there are compact neighbourhoods $K_s$ of $s^*$ and $K_c$ of $c^*$ such that, for $(s,c)\in K_s\times K_c$, $Pc$ is a regular value of $F_s$ and $F'_s$ and has a preimage under $F_s$ in $W$.

For such $(s,c)$, the solutions of $\Psi_0^s=c$ are nondegenerate and finite in number, at most $2^{k-1}$ by Lemma \ref{K:lem:trig}, and they lie in a compact set of $(\omega,\theta)$. By the uniform $C^1$ convergence $\Psi_r\to\Psi_0^s$, for $r$ smaller than some $r_0$ independent of $(s,c)\in K_s\times K_c$, the solutions of $\Psi_r=c$ with $\omega$ in a fixed large disc are in bijection with those of $\Psi_0^s=c$. Near each nondegenerate solution there is exactly one solution, by the implicit function theorem. Away from the solutions $|\Psi_0^s-c|$ is bounded below, so there is none. The same holds near $h_0'$. Hence every point $x(r,s,c)$ with $0<r<r_0$ and $(s,c)\in K_s\times K_c$ has at least one preimage, and at most $2^{k-1}+2^{k-1}=2^k$. Fix $r\in(0,r_0)$. The points $x(r,s,c)$ with $(s,c)\in K_s\times K_c$ belong to $O_V$, and the map $(s,c)\mapsto x(r,s,c)$, between spaces of dimension $(k-1)+(k+1)=2k$, is a diffeomorphism onto its image, with Jacobian $r^{k+1}\prod_js_j^{-1}$. So $O_V$ has positive measure.
\end{proof}

\begin{proof}[Proof of Theorem \ref{K:thm:main}]
Let $(M,h_0)$ be a nondegenerate corner. The nondegeneracy of $(M,h)$ is an open condition in $h$, and by Lemma \ref{K:lem:line} the immersion condition fails only on a real line of parameters. So we may move $h_0$ slightly to make $\Log_O|_{L_M}$ an immersion at $h_0$. Proposition \ref{K:prop:corner} gives $O_V\subset\A(V)$ with $1\le N\le2^k$ at regular values. If $V$ is nondegenerate, $\int N=\mathcal M(V)=\pi^{2k}$ by Corollary~\ref{L:cor:Mone}, and $N\le\binom{2k}k$ almost everywhere by Theorem A. Hence $\pi^{2k}\le\binom{2k}k\Vol(\A(V))-\big(\binom{2k}k-2^k\big)|O_V|$. Finally $\binom{2k}k>2^k$ for $k\ge2$, since $\binom{2k}k=\prod_{i=1}^k\frac{k+i}i$ and every factor is at least $2$, the factor $i=1$ being $k+1>2$.
\end{proof}

\section{Nondegenerate corners}\label{K:s4}\label{K:sec:CN}

\begin{proof}[Proof of Theorem \ref{K:thm:CN} for $k=2$]
Let $k=2$, $M=\{j\}$ and $O$ the other three indices. By Lemma \ref{K:lem:regular}, we may choose $h\in L_M$ at which $\Log_O|_{V_O}$ has rank three and $\Log_O|_{L_M}$ is an immersion. The logarithmic tangent space $L_h$ of $V_O$ at $\pi_O(h)$ is spanned over $\C$ by $d/h$ and $u_j/h$, and the kernel $K$ of $\operatorname{Re}|_{L_h}$ is a real line, since $\dim_\R L_h=4$ and the rank is three. Since $\Log_O|_{L_M}$ is an immersion, $K\cap\C(d/h)=0$. The matrix of Definition \ref{K:def:CN} represents $\operatorname{Re}$ on the real subspace $R(\theta)=\C(d/h)\oplus\R\,ie^{i\theta}u_j/h$ of $L_h$, up to signs of columns, because $\operatorname{Re}(iw)=-\operatorname{Im}w$. It is invertible if and only if $R(\theta)\cap K=0$. Let $\bar K$ be the image of $K$ in the complex line $L_h/\C(d/h)$, a real line through $0$, and identify $L_h/\C(d/h)$ with $\C$ through the basis $u_j/h$. Then $R(\theta)\cap K=0$ if and only if $\bar K\neq\R ie^{i\theta}$, which holds for all but one value of $\theta$ modulo $\pi$.
\end{proof}

\begin{proof}[Proof of Corollary \ref{K:cor:strict}]
As explained at the beginning of Section~\ref{K:s1}, $V$ can be defined by two affine linear equations all of whose coefficients are nonzero. By Theorem~\ref{K2:thm:nondeg}, $V$ is then nondegenerate, with $M_+(V)=1$ and $M_-(V)=6$. By Theorem~\ref{K:thm:CN}, $V$ has a nondegenerate corner, and Theorem~\ref{K:thm:main} gives $\Vol(\A(V))>\pi^4/6$.
\end{proof}

For $k\ge3$ the same reasoning shows that $(M,h)$ is nondegenerate if and only if $\Log_O|_{V_O}$ has rank $k+1$ at $\pi_O(h)$, $K\cap\C(d/h)=0$, and the real $(k-1)$-dimensional image $\bar K$ of $K$ in $L_h/\C(d/h)\simeq\C^M$ meets the totally real subspace $\bigoplus_j\R\,ie^{i\theta_j}e_j$ only in $0$ for some $\theta$. The first two conditions hold for generic $h$ by Lemmas \ref{K:lem:regular} and \ref{K:lem:line}. The third is a transversality condition between a real subspace and a family of coordinate totally real subspaces. Write $\bar K=A(\R^{k-1})$ with $A$ a complex $(k-1)\times(k-1)$ matrix. The condition fails for all $\theta$ if and only if the real determinant $\det\big(\cos\theta_j\operatorname{Re}A_j+\sin\theta_j\operatorname{Im}A_j\big)_j$, with rows indexed by $j$, vanishes identically in $\theta$. By the theorem of Rado on independent transversals [Ra], this happens if and only if there is a set $J$ of rows with $\dim_\R\operatorname{span}\{\operatorname{Re}A_j,\operatorname{Im}A_j:j\in J\}<|J|$. This can happen for special subspaces, for instance when $\bar K$ contains a complex coordinate line. So for $k\ge3$ this reasoning does not show directly that every plane in general position has a nondegenerate corner. A different argument, given in Section~\ref{K2:s1}, proves it for every $k$.

\subsection{What the corner method says about the general equality case}\label{K:s6}

The definition of a corner makes sense for any smooth complete intersection $V$ of dimension $k$ in $(\C^*)^{2k}$. Let $M$ be a set of $k-1$ coordinates such that $V$ meets $\{z_M=0\}$ in a curve $C_M$, transversally at a smooth point $h_0\in(\C^*)^O$. Replace $d/h_0$ by the logarithmic tangent vector of $C_M$ and the vectors $u_j/h_0$ by the logarithmic derivatives of $V$ in the directions of $z_j$, $j\in M$. We expect that the proof of Proposition~\ref{K:prop:corner} then extends: near a nondegenerate corner every regular fibre should have at most $2^{k-1}N_{C_M}$ points, where $N_{C_M}$ is the number of points of $C_M$ over $\Log_Oh_0$. We have not carried this out, and it is not used in the paper. If it holds, equality in the bound $\Vol(\A(V))\ge\pi^{2k}/M_-(V)$ requires $2^{k-1}N_{C_M}\ge M_-(V)$ at every nondegenerate corner. For twisted products of trinomial curves the corners are degenerate in the sense of Definition~\ref{K:def:CN}, because the coordinates of different factors separate, so no obstruction arises there.

\section{Every $k$-plane in general position has a nondegenerate corner}\label{K2:s1}

We now complete the argument of Sections~\ref{K:s1} to \ref{K:s6}, whose notation we keep, and recall it briefly. Let $V\subset\C^{2k}$ be an affine $k$-plane, $k\ge2$, with direction space $D$, and also write $V$ for $V\cap(\C^*)^{2k}$. Its closure $\bar V\simeq\PP^k$ in $\PP^{2k}$ meets the coordinate hyperplanes and the hyperplane at infinity in $2k+1$ hyperplanes $\lambda_1,\dots,\lambda_{2k},\lambda_0$. The plane is in general position if any $k+1$ of them have empty intersection. A point at infinity of $\bar V$ is the class $[y]$ of a nonzero $y\in D$, and it lies on $\lambda_l$, $l\ge1$, exactly when $y_l=0$.

Fix a set $M$ of $k-1$ coordinate indices and let $O$ be the complementary set of $k+1$ indices. By Lemma~\ref{K:lem:geometry}, $L_M=\{z\in V:z_M=0\}$, taken in $\C^{2k}$, is an affine line, $L_M=\{p_0+\tau d\}$ with $d\in D$ and $d_M=0$. For each $l\in O$ the coordinate $z_l$ restricted to $L_M$ is a nonconstant affine function of $\tau$, vanishing at a single point, and these points are distinct. For $j\in M$ there is $u_j\in D$ whose $M$-coordinates are $(u_j)_i=\delta_{ij}$, unique modulo $\C d$. For $h\in L_M$ with $h_l\neq0$ for $l\in O$, and $w\in\C^{2k}$, write $w/h=(w_l/h_l)_{l\in O}\in\C^O$. The corner $(M,h)$ is nondegenerate (Definition~\ref{K:def:CN}) if, for some $\theta\in(\R/2\pi\Z)^M$, the real $(k+1)\times(k+1)$ matrix
$$\mathcal N(h,\theta)=\Big[\operatorname{Re}\frac dh,\ \operatorname{Im}\frac dh,\ \operatorname{Im}\frac{e^{i\theta_j}u_j}h\ (j\in M)\Big]$$
is invertible. Its invertibility does not depend on the choice of $d$ in $\C^*d$, nor on the choice of the $u_j$ modulo $\C d$. Indeed, replacing $d$ by $\alpha d$ replaces the first two columns by real combinations of them, with determinant $|\alpha|^2$. And $\operatorname{Im}(e^{i\theta_j}(u_j+c\,d)/h)$ differs from $\operatorname{Im}(e^{i\theta_j}u_j/h)$ by $\operatorname{Im}(\beta d/h)=\operatorname{Re}\beta\operatorname{Im}(d/h)+\operatorname{Im}\beta\operatorname{Re}(d/h)$, with $\beta=ce^{i\theta_j}$, a combination of the first two columns.

\begin{theorem}\label{K2:thm:corner}
Let $V$ be a $k$-plane in general position in $(\C^*)^{2k}$, $k\ge2$. For every set $M$ of $k-1$ coordinate indices and every $s\in O$, there are phases $\psi\in\R/2\pi\Z$ and $\theta\in(\R/2\pi\Z)^M$ such that, if $p$ denotes the point of $L_M$ with $p_s=0$, the corner $(M,p+re^{i\psi}d)$ is nondegenerate for all sufficiently small $r>0$. Moreover $\det\mathcal N(p+re^{i\psi}d,\theta)=\pm\,\mathcal D/r+O(1)$ as $r\to0$, where $\mathcal D\neq0$ is the $k\times k$ determinant defined in (\ref{K2:eq:D}) below.
\end{theorem}

\begin{theorem}\label{K2:thm:nondeg}
Let $V$ be a $k$-plane in general position, defined by $k$ affine linear equations each of which has all its $2k+1$ coefficients nonzero. Then $V$ is a nondegenerate complete intersection in the sense of Khovanskii, every Newton polytope is the standard simplex $S$ of $\R^{2k}$, $M_+(V)=1$ and $M_-(V)=\binom{2k}k$.
\end{theorem}

\begin{corollary}\label{K2:cor:strict}
For every $k\ge2$ and every $k$-plane $V$ in general position in $(\C^*)^{2k}$, there is a set $O_V\subset\A(V)$ of positive measure on which every regular fibre of $\Log|_V$ has at most $2^k$ points, and
$$\Vol(\A(V))\ \ge\ \frac{\pi^{2k}+\big(\binom{2k}k-2^k\big)|O_V|}{\binom{2k}k}\ >\ \frac{\pi^{2k}}{\binom{2k}k}.$$
\end{corollary}

\section{Proof of Theorem \ref{K2:thm:corner}}\label{K2:s2}

\begin{lemma}\label{K2:lem:independent}
Let $s\in O$ and $O'=O\setminus\{s\}$, a set of $k$ indices. Let $p\in L_M$ be the point with $p_s=0$, so that $p_l\neq0$ for $l\in O'$. Normalise $d$ by $d_s=1$, and put $\tilde u_j=u_j-(u_j)_s\,d$ for $j\in M$. Then the $k$ vectors
$$w_0=\Big(\frac{d_l}{p_l}\Big)_{l\in O'},\qquad w_j=\Big(\frac{(\tilde u_j)_l}{p_l}\Big)_{l\in O'}\quad(j\in M)$$
are linearly independent over $\C$ in $\C^{O'}\simeq\C^k$.
\end{lemma}

\begin{proof}
The vectors $d,\tilde u_j$, $j\in M$, lie in $D$ and are linearly independent. Indeed, if $a\,d+\sum_jc_j\tilde u_j=0$, the $M$-coordinates give $\sum_jc_je_j=0$, since $d_M=0$ and $(\tilde u_j)_M=e_j$, so all $c_j$ vanish, and then $a=0$. Next, the projection $D\to\C^{O'}$, $y\mapsto y_{O'}$, is injective. If $y\in D\setminus\{0\}$ had $y_{O'}=0$, the point at infinity $[y]$ would lie on the $k$ hyperplanes $\lambda_l$, $l\in O'$, and on $\lambda_0$, which are $k+1$ hyperplanes with a common point, contradicting general position. So the projections of $d,\tilde u_j$ to $\C^{O'}$ are independent, and dividing the coordinates by the nonzero numbers $p_l$ preserves independence. That $p_l\neq0$ for $l\in O'$ follows from Lemma~\ref{K:lem:geometry}: the zeros of the $z_l$ on $L_M$ are distinct.
\end{proof}

\begin{lemma}\label{K2:lem:trigdet}
Let $w_1,\dots,w_n\in\C^n$ be linearly independent over $\C$. Then the real function $g(\varphi)=\det\big(\operatorname{Im}(e^{i\varphi_1}w_1),\dots,\operatorname{Im}(e^{i\varphi_n}w_n)\big)$ of $\varphi\in(\R/2\pi\Z)^n$ is not identically zero.
\end{lemma}

\begin{proof}
Since $\operatorname{Im}(e^{i\varphi}w)=\frac1{2i}(e^{i\varphi}w-e^{-i\varphi}\bar w)$, multilinearity of the determinant gives
$$g(\varphi)=\frac1{(2i)^n}\sum_{\epsilon\in\{\pm1\}^n}\Big(\prod_r\epsilon_r\Big)\,e^{i\langle\epsilon,\varphi\rangle}\det\big(w_1^{\epsilon_1},\dots,w_n^{\epsilon_n}\big),$$
where $w^{+1}=w$ and $w^{-1}=\bar w$. The exponents $\langle\epsilon,\varphi\rangle$ for distinct $\epsilon$ are distinct linear forms, so the functions $e^{i\langle\epsilon,\varphi\rangle}$ are linearly independent. The coefficient of $e^{i(\varphi_1+\dots+\varphi_n)}$ is $(2i)^{-n}\det(w_1,\dots,w_n)\neq0$. Hence $g\not\equiv0$.
\end{proof}

By the two lemmas, applied with $n=k$ and the vectors $w_0$ and $w_j$, $j\in M$, there are phases $\psi$ and $\theta_j$ such that
\begin{equation}\label{K2:eq:D}
\mathcal D=\det\Big(\operatorname{Im}(e^{i\psi}w_0),\ \operatorname{Im}(e^{i\theta_j}w_j)\ (j\in M)\Big)\neq0 .
\end{equation}

\begin{proof}[Proof of Theorem \ref{K2:thm:corner}]
Let $\tau=re^{i\psi}$ and $h=p+\tau d$, with $d_s=1$. Then $h_s=\tau$, and for $l\in O'$, $h_l=p_l+\tau d_l=p_l(1+O(r))$. By the invariance noted in Section~\ref{K2:s1} we may use $\tilde u_j$ instead of $u_j$. We have $(\tilde u_j)_s=0$. The coordinates of the relevant vectors are therefore the following. For $d/h$: the $s$-coordinate is $1/\tau$, and the $l$-coordinates are $(w_0)_l+O(r)$ for $l\in O'$. For $\tilde u_j/h$: the $s$-coordinate is $0$, and the $l$-coordinates are $(w_j)_l+O(r)$ for $l\in O'$.

Replace the first two columns of $\mathcal N(h,\theta)$, which span $\{\operatorname{Re}(xd/h):x\in\C\}$, by $\operatorname{Re}(\tau d/h)$ and $\operatorname{Re}(ie^{i\psi}d/h)$. Since $\operatorname{Re}(xw)=\operatorname{Re}x\operatorname{Re}w-\operatorname{Im}x\operatorname{Im}w$, the passage from $(\operatorname{Re}(d/h),\operatorname{Im}(d/h))$ to these two columns is given by the matrix $\begin{pmatrix}\operatorname{Re}\tau&\operatorname{Re}(ie^{i\psi})\\-\operatorname{Im}\tau&-\operatorname{Im}(ie^{i\psi})\end{pmatrix}$, whose determinant is $-\operatorname{Re}\tau\operatorname{Im}(ie^{i\psi})+\operatorname{Im}\tau\operatorname{Re}(ie^{i\psi})=-r\cos^2\psi-r\sin^2\psi=-r$. In the new matrix, ordering the rows with $s$ first, the first column is $\operatorname{Re}(\tau d/h)$. It has $s$-entry $\operatorname{Re}(1)=1$ and $O'$-entries $\operatorname{Re}(\tau(w_0+O(r)))=O(r)$. The second column $\operatorname{Re}(ie^{i\psi}d/h)$ has $s$-entry $\operatorname{Re}(ie^{i\psi}/\tau)=\operatorname{Re}(i/r)=0$ and $O'$-entries $\operatorname{Re}(ie^{i\psi}w_0)+O(r)=-\operatorname{Im}(e^{i\psi}w_0)+O(r)$. The remaining columns $\operatorname{Im}(e^{i\theta_j}\tilde u_j/h)$ have $s$-entry $0$ and $O'$-entries $\operatorname{Im}(e^{i\theta_j}w_j)+O(r)$. Expanding along the first row, the determinant of the new matrix is $-\mathcal D+O(r)$. Hence $\det\mathcal N(h,\theta)=(-\mathcal D+O(r))/(-r)=\mathcal D/r+O(1)$, which is nonzero for small $r>0$. For small $r$, $h_s=\tau\neq0$ and $h_l\neq0$ for $l\in O'$, so $h$ is an admissible point.
\end{proof}

\section{Proof of Theorem \ref{K2:thm:nondeg} and of the corollary}\label{K2:s3}

\begin{proof}[Proof of Theorem \ref{K2:thm:nondeg}]
Write the equations as $\sum_{m=1}^{2k}a_{jm}z_m+c_j=0$, $j=1,\dots,k$, with all $a_{jm}$ and $c_j$ nonzero, and let $a_m\in\C^k$ be the columns of $A=(a_{jm})$. We first show that any set $T$ of at most $k$ columns of $A$ is linearly independent. If $\sum_{m\in T}x_ma_m=0$ with $x\neq0$, then the vector $x$, extended by zero, lies in $D=\ker A$ and vanishes outside $T$. So the point at infinity $[x]$ lies on $\lambda_0$ and on the $2k-|T|\ge k$ hyperplanes $\lambda_l$, $l\notin T$, at least $k+1$ hyperplanes, a contradiction. Next, for $w\in\R^{2k}\setminus\{0\}$, let $T\subset\{0,1,\dots,2k\}$ be the set where the weight is maximal, the constant term having weight $0$ and $z_m$ weight $w_m$. All $k$ equations have full support, so every initial form is $\sum_{m\in T'}a_{jm}z_m+\epsilon c_j$, with $T'=T\setminus\{0\}$ and $\epsilon=1$ if $0\in T$, $\epsilon=0$ otherwise.

If $|T'|\ge k$, the $k\times|T'|$ matrix $(a_{jm})_{m\in T'}$ has rank $k$, since it contains $k$ independent columns. So the initial system defines a smooth complete intersection of dimension $k$, or the empty set, in $(\C^*)^{2k}$; the variables outside $T'$ are free. If $|T'|<k$ and $\epsilon=0$, the system $\sum_{m\in T'}a_{jm}z_m=0$ has only the solution $z_{T'}=0$, since the columns are independent, so there is no solution in the torus. If $|T'|<k$ and $\epsilon=1$, a solution $z_{T'}$ of $\sum_{m\in T'}a_{jm}z_m=-c_j$ gives the point of $V$ with coordinates $z_{T'}$ and $0$ elsewhere. That point lies on $2k-|T'|\ge k+1$ coordinate hyperplanes, a contradiction, so again there is no solution. Hence $V$ is nondegenerate. The Newton polytopes are the standard simplex $S=\operatorname{conv}(0,e_1,\dots,e_{2k})$. So $M_+(V)=(2k)!\Vol(S)=1$, and $M_-(V)=\MV(S^{(k)},(-S)^{(k)})=\binom{2k}k$ by the formula $\MV(S^{(j)},(-S)^{(n-j)})=\binom nj$ of Rogers and Shephard [RS, Section 4].
\end{proof}

Note that equations with full supports always exist. Generic linear combinations of a basis of equations of $V$ have all coefficients nonzero. A coefficient $a_{jm}$ vanishing for all combinations would mean that $z_m$ is constant on $V$, and $c_j$ vanishing for all combinations would mean that $V$ contains the origin. Both contradict general position.

\begin{proof}[Proof of Corollary \ref{K2:cor:strict}]
By Theorem \ref{K2:thm:corner}, $V$ has a nondegenerate corner $(M,h_0)$. The nondegeneracy is an open condition in $h_0$, and by Lemma~\ref{K:lem:line} the map $\Log_O|_{L_M}$ fails to be an immersion only on a real line of parameters. So we may move $h_0$ slightly so that it is also an immersion point. Proposition~\ref{K:prop:corner} then gives a set $O_V\subset\A(V)$ of positive measure on which every regular fibre has at least one and at most $2^k$ points. By Theorem \ref{K2:thm:nondeg} and Corollary~\ref{L:cor:Mone}, $\int N=\mathcal M(V)=\pi^{2k}$, and by Theorem A, $N\le\binom{2k}k$ almost everywhere. Hence $\pi^{2k}\le\binom{2k}k\Vol(\A(V))-\big(\binom{2k}k-2^k\big)|O_V|$. Finally $\binom{2k}k>2^k$ for $k\ge2$.
\end{proof}

This settles the equality case of the bound of Theorem C for $k$-planes in general position: none of them attains it.

\part{The case $n>2k$}\label{part:two}

\section{The case $n>2k$: projections and the main bounds}\label{N:s1}

We now turn to varieties of higher codimension. Let $V\subset(\C^*)^n$ be a smooth complete intersection of dimension $k$, defined by Laurent polynomials $f_1,\dots,f_{n-k}$ with Newton polytopes $\Delta_1,\dots,\Delta_{n-k}\subset\R^n$, so that $df_1\wedge\dots\wedge df_{n-k}\neq0$ on $V$. Its amoeba $\A(V)=\Log(V)$ is the image of a real $2k$-dimensional manifold under a smooth proper map. When $n=2k$ it is a closed subset of $\R^{2k}$ of finite volume, and Theorem A, combined with Theorems B, C and D, gives a lower bound $\pi^{2k}/\MV(\Delta_1,\dots,\Delta_k,-\Delta_1,\dots,-\Delta_k)$ for that volume under suitable hypotheses. When $n>2k$, the amoeba is a closed set of Hausdorff dimension at most $2k$ in a larger space, and the natural measure is the $2k$-dimensional Hausdorff measure $\Hm^{2k}$. In Part~\ref{part:two} we bound $\Hm^{2k}(\A(V))$ from below and from above. The present section states the main bounds and their corollaries, Sections~\ref{N:s2} to \ref{N:s7} prove them, and Section~\ref{N:s9} discusses them.

The idea is to project. An integer $2k\times n$ matrix $B$ of rank $2k$, with rows $b_1,\dots,b_{2k}$, defines the monomial map $\phi_B(z)=(z^{b_1},\dots,z^{b_{2k}})$, which satisfies $\Log\circ\phi_B=B\circ\Log$ and $\Arg\circ\phi_B=B\circ\Arg$. So the image of the amoeba of $V$ under the linear map $B$ is the amoeba of $\phi_B(V)$, a variety of dimension at most $k$ in $(\C^*)^{2k}$, which is half-dimensional when $\phi_B|_V$ is generically finite. A linear map does not increase $2k$-dimensional measures by more than the factor $\|B\|_{2k}=\sqrt{\det BB^{\mathsf T}}$, the largest $2k$-dimensional Jacobian of $B$. So lower bounds for amoebas of half-dimensional varieties give lower bounds for $\Hm^{2k}(\A(V))$. In the other direction, the Cauchy--Binet formula expresses the $2k$-dimensional Jacobian of $\Log|_V$ through the Jacobians of all coordinate projections. This gives upper bounds, and a second, better lower bound.

We use the following notation. For $B$ as above let $N_B(y)$, $y\in\R^{2k}$, be the number of points $z\in V$ with $B\Log z=y$, and let
$$M_B(V)=\int_{\R^{2k}}N_B(y)\,dy .$$
We show (Proposition \ref{N:prop:MB}) that $M_B(V)$ equals $\int_{\T^{2k}}\#\{z\in V:B\Arg z=\theta\}\,d\theta$. It is the multivolume of the coamoeba of $\phi_B(V)$ counted with the multiplicity of $\phi_B|_V$. For a set $I\subset\{1,\dots,n\}$ of $2k$ coordinates we write $M_I$, $N_I$ for the coordinate projection $B=\pi_I$. In $\R^{2n}$, with coordinates $(\alpha,\beta)\in\R^n\times\R^n$, consider the polytopes $\Delta_j\times\{0\}$ and $\{0\}\times\Delta_j$ for $j=1,\dots,n-k$, and the segments $[0,(b_i,b_i)]$ and $[0,(b_i,-b_i)]$ for $i=1,\dots,2k$. Define
$$\mu_B^-(V)=\MV_{2n}\big((\Delta_j\times0)_j,(0\times\Delta_j)_j,([0,(b_i,b_i)])_i\big),$$
$$\mu_B^+(V)=\MV_{2n}\big((\Delta_j\times0)_j,(0\times\Delta_j)_j,([0,(b_i,-b_i)])_i\big).$$
Finally let $N(x)$ be the number of points of $V$ over $x\in\R^n$, which is finite for $\Hm^{2k}$-almost every $x\in\A(V)$. Let $N_{\max}$ and $N_{\min}$ be its essential supremum and infimum on $\A(V)$ with respect to $\Hm^{2k}$.

\begin{theorem}\label{N:thm:main}
Let $V$ be an irreducible smooth complete intersection of dimension $k$ in $(\C^*)^n$, $n\ge2k$. Then $N_B\le\mu_B^-(V)$ almost everywhere for every $B$, $N_{\max}\le\mu_B^-(V)$ for every $B$ with $M_B(V)>0$, in particular $N_{\max}\le\min\{\mu_I^-(V):M_I(V)>0\}$, and
$$\Hm^{2k}(\A(V))\ \ge\ \sup_B\frac{M_B(V)}{\mu_B^-(V)\,\|B\|_{2k}},\eqtag{N:e1.1}$$
$$\Hm^{2k}(\A(V))\ \ge\ \frac1{N_{\max}}\Big(\sum_{|I|=2k}M_I(V)^2\Big)^{1/2}\ \ge\ \frac{\big(\sum_IM_I(V)^2\big)^{1/2}}{\min\{\mu_I^-(V):M_I(V)>0\}},\eqtag{N:e1.2}$$
$$\Hm^{2k}(\A(V))\ \le\ \frac1{N_{\min}}\sum_{|I|=2k}M_I(V)\ \le\ \frac{\pi^{2k}}{N_{\min}}\sum_{|I|=2k}\mu_I^+(V).\eqtag{N:e1.3}$$
In (\ref{N:e1.1}) the supremum runs over integer matrices of rank $2k$ with $M_B(V)>0$. If all $M_I(V)$ vanish, then $\Hm^{2k}(\A(V))=0$.
\end{theorem}

The upper bound (\ref{N:e1.3}) is obtained by the argument of [MN1, Section 3 and Theorem 4.1]: the Cauchy--Binet identity, the inequality $J\le\sum_IJ_I$, and division by the minimal fibre count. That argument gives, for rational curves, the estimate $p\,\Area(\A_f)\le\pi^2\sum_{j_1<j_2}n_{j_1}n_{j_2}$. The new ingredient in (\ref{N:e1.3}) is the bound $M_I(V)\le\pi^{2k}\mu_I^+(V)$ through the graph mixed volumes. For a line in general position, $n_j=1$ and $M_I(V)=\pi^2$ for every pair $I$, and both estimates give $\binom n2\pi^2/N_{\min}$. For the lower bounds, the first inequality of (\ref{N:e1.2}), $\Hm^{2k}(\A(V))\ge N_{\max}^{-1}\int_VJ$, combined with the identity between $\int_VJ$ and the multivolume of the coamoeba [MN2, Proposition 3.1 and its proof], is the lower half of the comparison theorem [MN2, Theorem 3.1]. The new ingredients are the Minkowski inequality, which gives $(\sum_IM_I^2)^{1/2}$, the bound of $N_{\max}$ by mixed volumes, and the single-projection bound (\ref{N:e1.1}).

For $n=2k$, the only coordinate set is $I=\{1,\dots,2k\}$, and by Proposition \ref{N:prop:graph}
$$\mu_I^\pm(V)=\MV(\Delta_1,\dots,\Delta_k,\pm\Delta_1,\dots,\pm\Delta_k).$$ So (\ref{N:e1.1}) with $B=I$, and (\ref{N:e1.3}), become
$$\frac{\mathcal M(V)}{M_-(V)}\le\Vol(\A(V))\le\mathcal M(V)\le\pi^{2k}M_+(V).$$
The left inequality is Theorem A(b), and the right ones follow from Theorem A(a), since $N\ge1$ on the amoeba.

The quantities $M_B(V)$ are controlled by the coamoeba bounds of Theorems B, C and D, applied to the projected varieties. Let $W_B$ be the Zariski closure of $\phi_B(V)$ in $(\C^*)^{2k}$. Let $d_B^+(V)$ be the number of points of a generic fibre of the map $V\times\bar V\to(\C^*)^{2k}$, $(z,w)\mapsto\phi_B(z)/\phi_B(w)$, where $\bar V$ is the complex conjugate variety. We set $d_B^+(V)=0$ if this map is not dominant.

\begin{theorem}\label{N:thm:coamoeba}
Let $V$ be irreducible. Then $M_B(V)=\deg(\phi_B|_V)\,\mathcal M(W_B)$ if $\phi_B|_V$ is generically finite, and $M_B(V)=0$ otherwise. Moreover $d_B^+(V)\le\mu_B^+(V)$, and if $d_B^+(V)$ is odd, then $M_B(V)\ge\pi^{2k}$. For $k=1$, if the Newton polygon of a reduced equation of the plane curve $W_B$ is not centrally symmetric, then $M_B(V)\ge\pi^2$.
\end{theorem}

\begin{corollary}\label{N:cor:explicit}
Let $V$ be irreducible, and let $G$ be the set of coordinate sets $I$ with $|I|=2k$ for which $M_I(V)\ge\pi^{2k}$, for instance those with $d_I^+(V)$ odd. Then
$$\frac{\pi^{2k}\sqrt{\#G}}{\min\{\mu_I^-(V):M_I(V)>0\}}\ \le\ \Hm^{2k}(\A(V))\ \le\ \frac{\pi^{2k}}{N_{\min}}\sum_{|I|=2k}\mu_I^+(V).$$
\end{corollary}

For coordinate projections the mixed volumes $\mu_I^\pm$ reduce to mixed volumes in $\R^{2n-2k}$ of graph polytopes. Let $J$ be the complement of $I$, and for $\Delta\subset\R^n$ let $\Gamma_J(\Delta)=\{(a_J,a):a\in\Delta\}\subset\R^J\times\R^n$, where $a_J$ is the vector of coordinates of $a$ indexed by $J$.

\begin{proposition}\label{N:prop:graph}
For every coordinate set $I$ of size $2k$, with complement $J$,
$$\mu_I^\pm(V)=\MV_{2n-2k}\big(\Gamma_J(\Delta_1),\dots,\Gamma_J(\Delta_{n-k}),\ 0\times(\pm\Delta_1),\dots,0\times(\pm\Delta_{n-k})\big).$$
\end{proposition}

\begin{corollary}\label{N:cor:planes}
Let $V$ be a $k$-plane in general position in $(\C^*)^n$, $n\ge2k$, that is, such that any $k+1$ of the $n+1$ hyperplanes $\lambda_0,\lambda_1,\dots,\lambda_n$ cut on $\bar V\simeq\PP^k$ by the hyperplane at infinity and the coordinate hyperplanes have empty intersection, and defined by $n-k$ affine linear equations with all coefficients nonzero. Then, for every $I$, $\pi_I|_V$ is an isomorphism onto a $k$-plane in general position in $(\C^*)^{2k}$, $M_I(V)=\pi^{2k}$, $N_I\le\binom{2k}k$ almost everywhere, and
$$\frac{\binom n{2k}^{1/2}}{\binom{2k}k}\,\pi^{2k}\ \le\ \Hm^{2k}(\A(V))\ \le\ \binom n{2k}\,\pi^{2k}.$$
If $V$ is real up to multiplicative translation, the upper bound can be divided by $2$.
\end{corollary}

\begin{corollary}\label{N:cor:curves}
Let $V\subset(\C^*)^n$ be a smooth irreducible complete intersection curve defined by $n-1$ polynomials with Newton polytopes $d_jS_n$, where $S_n$ is the standard simplex, and let $D=\prod_jd_j$. Suppose that for every pair $I$ of coordinates the projection $\pi_I|_V$ is birational onto a plane curve with Newton polygon $D\cdot\conv(0,e_1,e_2)$, which holds for generic coefficients by Theorem~\ref{F:thm:birational}. Then $\pi^2\le M_I(V)\le D^2\pi^2$ and $N_I\le2D^2$ almost everywhere, for all $I$, and
$$\frac{\binom n2^{1/2}}{2D^2}\,\pi^2\ \le\ \Hm^2(\A(V))\ \le\ \binom n2D^2\,\pi^2 .$$
\end{corollary}

The second inequality of (\ref{N:e1.2}) uses only $N_{\max}\le\min_I\mu_I^-$. For generic $V$ with $n\ge2k+1$ one has $N=1$ almost everywhere when the differences of exponents span $\Z^n$, since pairs of points with the same image form a set of real dimension $4k-n<2k$; this is Theorem~\ref{F:thm:injective} below. Then (\ref{N:e1.2}) and (\ref{N:e1.3}) become $(\sum_IM_I^2)^{1/2}\le\Hm^{2k}(\A(V))\le\sum_IM_I$. For instance, for a $k$-plane with $N=1$ almost everywhere, $\binom n{2k}^{1/2}\pi^{2k}\le\Hm^{2k}(\A(V))\le\binom n{2k}\pi^{2k}$. We keep $N_{\max}$ in the statements, since generic injectivity requires genericity hypotheses.

\section{Preliminaries}\label{N:s2}

\subsection{Rectifiability and the area formula}\label{N:s2.1}
The map $\Log|_V:V\to\R^n$ is smooth and proper, so $\A(V)$ is closed. Fix a Riemannian metric on $V$. The $2k$-dimensional Jacobian of $\Log|_V$ at $z$ is $J(z)=\sqrt{\det(G_z)}$, where $G_z$ is the Gram matrix of the images under $d\Log$ of an orthonormal basis of $T_zV$. The area formula for smooth maps from a $2k$-dimensional manifold into $\R^n$ [Fe, Theorems 3.2.3 and 3.2.5] gives, for every nonnegative Borel function $g$ on $\R^n$,
$$\int_Vg(\Log z)\,J(z)\,d\mathrm{vol}(z)=\int_{\R^n}g(x)\,N(x)\,d\Hm^{2k}(x).\eqtag{N:e2.1}$$
In particular $N$ is finite $\Hm^{2k}$-almost everywhere on $\A(V)$ whenever $\int_VJ<\infty$. The same formula holds for any smooth map $F:V\to\R^m$, $m\ge2k$, with the corresponding Jacobian. The finiteness of $\Hm^{2k}(\A(V))$ for $n\ge2k$ is due to Madani and the author [MN1, Theorem 1.1], for every algebraic variety none of whose components lies in a subtorus of dimension less than $2k$; the proof in [MN1] is written for generic varieties. For complete intersections it also follows from the Cauchy--Binet bound and the fibre bound below, which give the explicit estimate (\ref{N:e1.3}). The identity (\ref{N:e2.2}) and the inequality $J\le\sum_IJ_I$ of Section~\ref{N:s2.2} appear, for every $k$, in [MN1, Section 3, (8) and (9)], see also [MN2, Section 3]; in [MN1] they are applied to curves.

\subsection{Cauchy--Binet}\label{N:s2.2}
Let $A$ be a real $n\times2k$ matrix, the matrix of $d\Log_z$ in an orthonormal basis of $T_zV$. Then $\det(A^{\mathsf T}A)=\sum_{|I|=2k}\det(A_I)^2$, where $A_I$ is the $2k\times2k$ submatrix of the rows in $I$. So
$$J(z)=\Big(\sum_{|I|=2k}J_I(z)^2\Big)^{1/2},\eqtag{N:e2.2}$$
where $J_I(z)=|\det A_I|$ is the Jacobian of $\pi_I\circ\Log|_V$ at $z$. In particular $\frac1{\sqrt{\binom n{2k}}}\sum_IJ_I\le J$ and $\max_IJ_I\le J\le\sum_IJ_I$. More generally, for a real $2k\times n$ matrix $B$, Cauchy--Binet gives $\det(BA)=\sum_I\det(B^I)\det(A_I)$, where $B^I$ is the $2k\times2k$ submatrix of the columns in $I$, and the Cauchy--Schwarz inequality gives
$$J_B(z)=|\det(BA)|\le\sum_I|\det B^I|\,J_I(z)\le\|B\|_{2k}\,J(z),\eqtag{N:e2.2'}$$
where $\|B\|_{2k}=(\sum_I(\det B^I)^2)^{1/2}=\sqrt{\det BB^{\mathsf T}}$.

\subsection{Linear maps}\label{N:s2.3}
For a real $2k\times n$ matrix $B$ of rank $2k$, and a $2k$-dimensional subspace $T\subset\R^n$, the Jacobian of $B|_T$ is at most $\|B\|_{2k}=\sqrt{\det BB^{\mathsf T}}$. Indeed, by Cauchy--Binet applied to $BU$, where $U$ is an orthonormal basis matrix of $T$, $|\det(BU)|$ is at most the product of the singular values of $B$, which is $\sqrt{\det BB^{\mathsf T}}$, with equality when $T$ is the row space of $B$. Since $\A(V)$ is the image of a $\sigma$-compact $2k$-manifold under a locally Lipschitz map, it is countably $2k$-rectifiable. When $\Hm^{2k}(\A(V))<\infty$, which is proved in Theorem \ref{N:thm:main} from (\ref{N:e2.2}) and Proposition \ref{N:prop:fibreB}, the area formula for Lipschitz maps on rectifiable sets [Fe, Theorem 3.2.20] gives
$$\Vol_{2k}(B(\A(V)))\le\int_{\A(V)}J^{\A}_B\,d\Hm^{2k}\le\|B\|_{2k}\,\Hm^{2k}(\A(V)),\eqtag{N:e2.3}$$
where $J^{\A}_B\le\|B\|_{2k}$ is the Jacobian of $B$ on the approximate tangent planes of $\A(V)$. We shall in fact prove (\ref{N:e1.1}) directly from (\ref{N:e2.2'}), without using (\ref{N:e2.3}).

\subsection{Mixed volumes}\label{N:s2.4}
We use the normalised mixed volume of Section~\ref{L:s2}, with $\MV(P,\dots,P)=m!\Vol(P)$ in $\R^m$, its invariance under unimodular linear maps, and the block factorisation (Lemma~\ref{L:lem:block}). If $P_1,\dots,P_a$ lie in the coordinate subspace $\R^a\times0$, then $\MV(P_1,\dots,P_a,R_1,\dots,R_b)=\MV_a(P_1,\dots,P_a)\,\MV_b(\varrho R_1,\dots,\varrho R_b)$, where $\varrho$ forgets the first $a$ coordinates. We also use the first part of Bernstein's theorem: the number of isolated common zeros in $(\C^*)^m$ of $m$ Laurent polynomials, counted with multiplicity, is at most the mixed volume of their Newton polytopes [Be], [HS, Theorem 1], [Fu, Section 5.5]; the statement for the torus, with isolated zeros only, is Theorem 1 of [HS].

\section{Projections and the coamoeba}\label{N:s3}

Fix $B$ and put $F_B=B\circ\Log|_V$ and $\Theta_B=B\circ\Arg|_V$, maps from $V$ to $\R^{2k}$ and to $\T^{2k}$. Let $J_B$ and $J'_B$ be their Jacobians. At $z\in V$, let $L_z=\{(u_j/z_j)_j:u\in T_zV\}\subset\C^n$ be the logarithmic tangent space, a complex subspace of dimension $k$. Then $dF_B=B\operatorname{Re}|_{L_z}$ and $d\Theta_B=B\operatorname{Im}|_{L_z}$.

\begin{lemma}\label{N:lem:jacB}
$J_B=J'_B$ at every point of $V$.
\end{lemma}

This is a direct consequence of [MN2, Proposition 3.1] and its proof.

\begin{proof}
Since $B$ is real, $B\operatorname{Im}u=B\operatorname{Re}(-iu)$, so $B\operatorname{Im}|_{L_z}=B\operatorname{Re}|_{L_z}\circ m$, where $m$ is multiplication by $-i$ on $L_z$. As a real linear map, $m$ has determinant $|(-i)^k|^2=1$.
\end{proof}

\begin{proposition}\label{N:prop:MB}
One has $M_B(V)=\int_VJ_B=\int_{\T^{2k}}N'_B(\theta)\,d\theta$, where $N'_B(\theta)=\#\{z\in V:B\Arg z=\theta\}$. If $V$ is irreducible and $\phi_B|_V$ is generically finite, then $M_B(V)=\deg(\phi_B|_V)\,\mathcal M(W_B)$, where $\mathcal M(W_B)$ is the multivolume of the coamoeba of the $k$-dimensional variety $W_B\subset(\C^*)^{2k}$. If $\phi_B|_V$ is not generically finite, then $M_B(V)=0$.
\end{proposition}

\begin{proof}
The first two equalities are the area formulas for $F_B$ and $\Theta_B$, which are maps between manifolds of the same dimension, together with Lemma \ref{N:lem:jacB}. Suppose that $V$ is irreducible and $\phi_B|_V$ generically finite of degree $\delta$. Then $W_B$ is irreducible of dimension $k$. There is a proper subvariety $Z\subset W_B$, containing the singular locus of $W_B$ and $W_B\setminus\phi_B(V)$, outside which every fibre of $\phi_B|_V$ has exactly $\delta$ points. Since $\Arg\circ\phi_B=B\circ\Arg$, for $\theta$ outside the null set $\Arg(Z)$, of real dimension at most $2k-2$, one has $N'_B(\theta)=\delta\cdot\#\{w\in W_B:\Arg w=\theta\}$. Integrating gives $M_B(V)=\delta\,\mathcal M(W_B)$. If $\phi_B|_V$ is not generically finite, then $\phi_B(V)$ has dimension less than $k$. In characteristic zero, the maximal rank of $d\phi_B|_V$ equals the dimension of the image, so $\operatorname{rank}_\C d\phi_B|_{T_zV}<k$ at every point. Since $d\Theta_B=\operatorname{Im}\circ d\log\phi_B$, its real rank is less than $2k$ everywhere, so $J'_B\equiv0$ and $M_B(V)=0$. Here, for a possibly singular $W_B$ which need not be a complete intersection, $\mathcal M(W_B)$ is defined as $\int_{\T^{2k}}\#\{w\in W_B^{\mathrm{reg}}:\Arg w=\theta\}\,d\theta$, as for reduced curves in Part~\ref{part:one}.
\end{proof}

\section{The fibre bound for projections}\label{N:s4}

\begin{lemma}\label{N:lem:transB}
Let $\psi_\pm:V\times\bar V\to(\C^*)^{2k}$ be the maps $\psi_-(z,w)=\phi_B(z)\phi_B(w)$ and $\psi_+(z,w)=\phi_B(z)/\phi_B(w)$. A point $z\in V$ satisfies $B\Log z=y$ if and only if $(z,\bar z)\in\psi_-^{-1}(e^{2y})$. It satisfies $B\Arg z\equiv\theta$ modulo $\pi$ if and only if $(z,\bar z)\in\psi_+^{-1}(e^{2i\theta})$. If $z$ is a regular point of $F_B$, the corresponding point $(z,\bar z)$ is an isolated and nondegenerate point of the fibre of $\psi_-$, and also of the fibre of $\psi_+$.
\end{lemma}

\begin{proof}
Since $\phi_B(\bar z)=\overline{\phi_B(z)}$, we have $\psi_-(z,\bar z)=|\phi_B(z)|^2=e^{2B\Log z}$ and $\psi_+(z,\bar z)=\phi_B(z)/\overline{\phi_B(z)}=e^{2iB\Arg z}$, which gives the first two assertions. The logarithmic tangent space of $V\times\bar V$ at $(z,\bar z)$ is $L_z\times\bar L_z$. The logarithmic differentials of $\psi_\pm$ there are $(u,v)\mapsto B(u\pm'v)$, where $\pm'$ is $+$ for $\psi_-$ and $-$ for $\psi_+$. Let $(u,v)\in L_z\times\bar L_z$ be in the kernel. Put $a=u+\bar v$ and $b=u-\bar v$, both in $L_z$. Then $u=(a+b)/2$ and $v=(\bar a-\bar b)/2$, so $u+v=\operatorname{Re}a+i\operatorname{Im}b$ and $u-v=i\operatorname{Im}a+\operatorname{Re}b$. For $\psi_-$, the real and imaginary parts of $B(u+v)=0$ give $B\operatorname{Re}a=0$ and $B\operatorname{Im}b=0$. Since $z$ is regular, $B\operatorname{Re}|_{L_z}$ is injective, and so is $B\operatorname{Im}|_{L_z}=B\operatorname{Re}|_{L_z}\circ m$. Hence $a=b=0$ and $u=v=0$. The case of $\psi_+$ is the same with the roles of $a$ and $b$ exchanged. A point where the differential of a map between complex manifolds of equal dimension $2k$ is injective is an isolated nondegenerate point of its fibre.
\end{proof}

\begin{proposition}\label{N:prop:fibreB}
For every regular value $y$ of $F_B$, $N_B(y)\le\mu_B^-(V)$. For almost every $\theta$, the number $N^\pi_B(\theta)$ of points $z\in V$ with $B\Arg z\equiv\theta$ modulo $\pi$ satisfies $N^\pi_B(\theta)\le\mu_B^+(V)$. Hence $M_B(V)\le\pi^{2k}\mu_B^+(V)$.
\end{proposition}

\begin{proof}
The fibre $\psi_-^{-1}(e^{2y})$ is the common zero set in $(\C^*)^{2n}$, with coordinates $(z,w)$, of the $2n$ Laurent polynomials $f_j(z)$, $\bar f_j(w)$ for $j=1,\dots,n-k$, and $z^{b_i}w^{b_i}-e^{2y_i}$ for $i=1,\dots,2k$. Here $\bar f_j$ has the conjugate coefficients, and $\bar V$ is its zero set. Their Newton polytopes are $\Delta_j\times0$, $0\times\Delta_j$ and $[0,(b_i,b_i)]$. By Lemma \ref{N:lem:transB}, the points $(z,\bar z)$ with $z$ in the fibre of $F_B$ over a regular value are distinct isolated zeros, so Bernstein's bound gives $N_B(y)\le\mu_B^-(V)$. The fibre of $\psi_+$ over $e^{2i\theta}$ is defined in the same way with the binomials $z^{b_i}w^{-b_i}-e^{2i\theta_i}$, whose Newton polytopes are the segments $[0,(b_i,-b_i)]$. By Sard's theorem and Lemma \ref{N:lem:jacB}, for almost every $\theta$ all points $z$ with $B\Arg z\equiv\theta$ modulo $\pi$ are regular points of $F_B$, and the same argument gives $N^\pi_B(\theta)\le\mu_B^+(V)$. Finally $M_B(V)=\int_{\T^{2k}}N'_B=\int_{[0,\pi)^{2k}}N^\pi_B$, since the $2^{2k}$ translates of $[0,\pi)^{2k}$ by vectors of $\{0,\pi\}^{2k}$ tile the torus.
\end{proof}

\begin{proof}[Proof of Theorem \ref{N:thm:main}]
Let $B$ satisfy $M_B(V)>0$. The zero set $Z_B$ of $J_B$ is the zero set of the function $\det(B\operatorname{Re}|_{L_z})$, computed in local holomorphic coordinates of $V$ with a real basis of $L_z$ depending real-analytically on $z$. So $Z_B$ is a real-analytic subset of the connected manifold $V$. It is proper, because $\int_VJ_B=M_B(V)>0$. Hence it has dimension at most $2k-1$: it is a locally finite union of real-analytic submanifolds of dimension at most $2k-1$. Its image under the locally Lipschitz map $\Log$ is therefore $\Hm^{2k}$-null. For $x\in\A(V)\setminus\Log(Z_B)$, every point $z$ of $V$ over $x$ is a regular point of $F_B$, and lies over $Bx$. By Lemma \ref{N:lem:transB}, the corresponding points $(z,\bar z)$ are distinct isolated points of $\psi_-^{-1}(e^{2Bx})$, and Bernstein's bound gives $N(x)\le\mu_B^-(V)$. So $N_{\max}\le\mu_B^-(V)$. For (\ref{N:e1.1}), by (\ref{N:e2.2'}) and (\ref{N:e2.1}) with $g=1$, $M_B(V)=\int_VJ_B\le\|B\|_{2k}\int_VJ=\|B\|_{2k}\int_{\A(V)}N\,d\Hm^{2k}\le\|B\|_{2k}N_{\max}\Hm^{2k}(\A(V))$, and $N_{\max}\le\mu_B^-(V)$. The same inequalities show that the right-hand side of (\ref{N:e1.1}) is at most $(\sum_IM_I^2)^{1/2}/N_{\max}$. Indeed $M_B\le\sum_I|\det B^I|M_I\le\|B\|_{2k}(\sum_IM_I^2)^{1/2}$ and $N_{\max}\le\mu_B^-$. So (\ref{N:e1.1}) is implied by the first inequality of (\ref{N:e1.2}). It can only improve the computable form of (\ref{N:e1.2}), with $\min_I\mu_I^-$ in the denominator, when $\mu_B^-$ is smaller than all the $\mu_I^-$ with $M_I>0$. For (\ref{N:e1.2}), apply (\ref{N:e2.1}) with $g=1/N$ on $\A(V)$, where $N$ is finite and positive almost everywhere. Then $\Hm^{2k}(\A(V))=\int_VJ(z)/N(\Log z)\,d\mathrm{vol}\ge\frac1{N_{\max}}\int_VJ$. By (\ref{N:e2.2}) and the triangle inequality for integrals of vector-valued functions,
$$\int_VJ=\int_V\big\|(J_I)_I\big\|_2\ \ge\ \Big\|\Big(\int_VJ_I\Big)_I\Big\|_2=\Big(\sum_IM_I(V)^2\Big)^{1/2},$$
using $\int_VJ_I=M_I(V)$ from Proposition \ref{N:prop:MB}. For (\ref{N:e1.3}), the same identity with $g=1/N$ gives $\Hm^{2k}(\A(V))\le\frac1{N_{\min}}\int_VJ\le\frac1{N_{\min}}\sum_I\int_VJ_I=\frac1{N_{\min}}\sum_IM_I(V)$, and Proposition \ref{N:prop:fibreB} gives $M_I(V)\le\pi^{2k}\mu_I^+(V)$. If $\int_VJ=\infty$, then (\ref{N:e2.2}) shows that $\sum_IM_I=\infty$, which is excluded by Proposition \ref{N:prop:fibreB}. So all quantities are finite.
\end{proof}

\section{Coamoeba lower bounds for projections}\label{N:s5}

\begin{proof}[Proof of Theorem \ref{N:thm:coamoeba}]
The formula for $M_B(V)$ is Proposition \ref{N:prop:MB}. Assume $V$ irreducible, so that $V\times\bar V$ is irreducible of dimension $2k$. If $\psi_+$ is dominant, it is generically finite between irreducible varieties of the same dimension. In characteristic zero, there is a dense Zariski open set $U\subset(\C^*)^{2k}$ over which $\psi_+$ is finite and \'etale, with fibres of exactly $d_B^+(V)$ points. These points are isolated zeros of the system of Proposition \ref{N:prop:fibreB}, so $d_B^+(V)\le\mu_B^+(V)$ by Bernstein's bound. The unit torus is Zariski dense in $(\C^*)^{2k}$ (the proof of Proposition~\ref{L:prop:parity}), so the set of $\theta$ with $e^{2i\theta}\notin U$ is null. For the other $\theta$, the antiholomorphic involution $\iota(z,w)=(\bar w,\bar z)$ of $V\times\bar V$ preserves $\psi_+^{-1}(e^{2i\theta})$. Indeed $\psi_+(\bar w,\bar z)=\overline{\phi_B(w)}/\overline{\phi_B(z)}=\overline{\psi_+(z,w)}^{\,-1}$, and this equals $e^{2i\theta}$ when $\psi_+(z,w)=e^{2i\theta}$. Its fixed points are the pairs $(z,\bar z)$, which by Lemma \ref{N:lem:transB} correspond to the points counted by $N^\pi_B(\theta)$. The orbits of $\iota$ on the finite fibre have one or two elements, so $N^\pi_B(\theta)\equiv d_B^+(V)$ modulo $2$. If $d_B^+(V)$ is odd, then $N^\pi_B\ge1$ almost everywhere, and $M_B(V)=\int_{[0,\pi)^{2k}}N^\pi_B\ge\pi^{2k}$. For $k=1$ and $\phi_B|_V$ generically finite, $W_B$ is a plane curve. If a reduced equation of $W_B$ has a Newton polygon which is not centrally symmetric, Proposition~\ref{L:prop:rays} gives $\mathcal M(W_B)\ge\pi^2$, and $M_B(V)\ge\deg(\phi_B|_V)\pi^2\ge\pi^2$.
\end{proof}

\begin{proof}[Proof of Corollary \ref{N:cor:explicit}]
By (\ref{N:e1.2}),
$$\Hm^{2k}(\A(V))\ge\frac{(\sum_IM_I^2)^{1/2}}{\min_I\mu_I^-}\ge\frac{(\#G\,\pi^{4k})^{1/2}}{\min_I\mu_I^-}.$$
The upper bound is (\ref{N:e1.3}).
\end{proof}

\section{The graph formula}\label{N:s6}

\begin{proof}[Proof of Proposition \ref{N:prop:graph}]
Let $I$ have $2k$ elements, with complement $J$, and treat $\mu_I^-$ first. The linear map $G(\alpha,\beta)=(\alpha,\beta-\alpha)$ of $\R^{2n}$ is unimodular. It maps the segment $[0,(e_m,e_m)]$, $m\in I$, to $[0,(e_m,0)]$, the polytope $\Delta_j\times0$ to $\{(a,-a):a\in\Delta_j\}$, and $0\times\Delta_j$ to itself. The $2k$ segments $[0,(e_m,0)]$, $m\in I$, lie in the coordinate subspace of the coordinates $\alpha_I$, and their normalised mixed volume there is $1$, the coefficient of $\prod\lambda_m$ in $\Vol(\sum\lambda_m[0,e_m])=\prod\lambda_m$. By the block factorisation, $\mu_I^-$ is the mixed volume, in $\R^J\times\R^n$, of the projections forgetting $\alpha_I$, namely $\{(a_J,-a):a\in\Delta_j\}$ and $0\times\Delta_j$. Composing with the unimodular map $(\alpha_J,\beta)\mapsto(\alpha_J,-\beta)$ gives $\Gamma_J(\Delta_j)$ and $0\times(-\Delta_j)$. For $\mu_I^+$, the map $(\alpha,\beta)\mapsto(\alpha,\beta+\alpha)$ sends $[0,(e_m,-e_m)]$ to $[0,(e_m,0)]$, $\Delta_j\times0$ to $\{(a,a)\}$ and $0\times\Delta_j$ to itself. Then the same argument gives $\Gamma_J(\Delta_j)$ and $0\times\Delta_j$. For $n=2k$, $J=\emptyset$, $\Gamma_J(\Delta)=\Delta$, and $\mu^\pm=\MV(\Delta_1,\dots,\Delta_k,\pm\Delta_1,\dots,\pm\Delta_k)$.
\end{proof}

\section{Proofs of the corollaries}\label{N:s7}

\begin{proof}[Proof of Corollary \ref{N:cor:planes}]
The direction space $D_V$ of $V$ projects injectively to $\C^I$. Otherwise a nonzero vector $y$ of $D_V$ with $y_I=0$ would give a point at infinity of $\bar V\simeq\PP^k$ lying on the $2k$ hyperplanes $\lambda_l$, $l\in I$, and on the hyperplane at infinity. These are at least $k+1$ hyperplanes, which contradicts general position. So $\pi_I|_V$ is an isomorphism onto the $k$-plane $W_I=\pi_I(V)\subset\C^{2k}$. Its hyperplane arrangement consists of the hyperplanes $\lambda_l$, $l\in I$, and $\lambda_0$, and it is in general position. By the remark following Theorem~\ref{K2:thm:nondeg}, $W_I$ can be defined by $k$ affine linear equations with all coefficients nonzero, and then by Theorem~\ref{K2:thm:nondeg} $W_I$ is nondegenerate with $M_+(W_I)=1$ and $M_-(W_I)=\binom{2k}k$, so $\mathcal M(W_I)=\pi^{2k}$ by Corollary~\ref{L:cor:Mone}. By Theorem A, the fibres of $\Log|_{W_I}$ have at most $\binom{2k}k$ points almost everywhere. Since $\pi_I|_V$ is an isomorphism, $M_I(V)=\mathcal M(W_I)=\pi^{2k}$ and $N_I\le\binom{2k}k$ almost everywhere. The proof of the bound $N_{\max}\le\mu_B^-$ in Theorem \ref{N:thm:main} applies with $\mu_I^-$ replaced by the bound $\binom{2k}k$ for the fibres of $\Log|_{W_I}$. Indeed, for $x$ outside a null set, every point of $V$ over $x$ corresponds to a regular point $w$ of $\Log|_{W_I}$ over $y=\pi_Ix$, by the argument with the set $Z_I$ in the proof of Theorem \ref{N:thm:main}. Such points $w$ are transverse, hence isolated, points of $W_I\cap\sigma_y(W_I)$ by Lemma~\ref{L:lem:transverse}, and Bernstein's bound for isolated zeros bounds their number by $M_-(W_I)$. The same argument is used in the proof of Corollary \ref{N:cor:curves}. Hence $N_{\max}\le\binom{2k}k$. The upper bound in (\ref{N:e1.3}) is $\sum_IM_I/N_{\min}\le\binom n{2k}\pi^{2k}$, and (\ref{N:e1.2}) gives the lower bound. If $V$ is real up to multiplicative translation, a twisted conjugation preserves $V$ and $\Log$, and its fixed points form a real $k$-dimensional subset whose image is $\Hm^{2k}$-null. So almost every fibre consists of pairs, and $N_{\min}\ge2$. This argument does not use the graph mixed volumes.
\end{proof}

\begin{proof}[Proof of Corollary \ref{N:cor:curves}]
Fix a pair $I$, and let $g$ be a reduced equation of the plane curve $W_I$, with Newton polygon $P=D\cdot T_0$, $T_0=\conv(0,e_1,e_2)$. Since $\pi_I|_V$ is birational, it is an isomorphism outside finite sets, so $M_I(V)=\mathcal M(W_I)$ and $N_I(y)=N_{W_I}(y)$ for almost every $y$. The triangle $P$ is not centrally symmetric, so $\mathcal M(W_I)\ge\pi^2$ by Proposition~\ref{L:prop:rays}, which needs only that $g$ is reduced. At smooth points of $W_I$, the transversality lemma of Section~\ref{L:s3} applies to $W_I$ and $\sigma_y(W_I)$, defined by $g$ and by the polynomial $\overline{g(\sigma_y(w))}$ with Newton polygon $-P$. The singular points of $W_I$ are finitely many, so for almost every $y$ Bernstein's bound gives $N_{W_I}(y)\le\MV(P,-P)=4\Area(P)=2D^2$. In the same way the rolled argument count of $W_I$ is at most $\MV(P,P)=2\Area(P)=D^2$ almost everywhere, so $\mathcal M(W_I)\le D^2\pi^2$. As in the proof of Corollary \ref{N:cor:planes}, the bound $N_I\le2D^2$ gives $N_{\max}\le2D^2$. The bounds then follow from (\ref{N:e1.2}) and (\ref{N:e1.3}), with $N_{\min}\ge1$.
\end{proof}

\section{Remarks and open questions}\label{N:s9}

The lower bounds (\ref{N:e1.1}) and (\ref{N:e1.2}) are of two different kinds. The first uses a single projection and the fibre bound for that projection. The second uses all coordinate projections at once, through Cauchy--Binet, and the fibre count of $\Log|_V$ itself. As shown in the proof of Theorem \ref{N:thm:main}, the first is always implied by the second, and it is useful only through the fibre bound $\mu_B^-$. For generic $V$ with $n\ge2k+1$ this fibre count is $|H|$ almost everywhere by Theorem~\ref{F:thm:injective}, so that, when the differences of exponents span $\Z^n$, (\ref{N:e1.2}) gives the bound $(\sum_IM_I^2)^{1/2}$ with no division at all. The real varieties, with fibres of even cardinality, show that it can fail for special $V$. For $k$-planes it would be natural to look for the minimum of $\Hm^{2k}(\A(V))$ over planes in general position, as in the half-dimensional case. Real planes are natural candidates, by analogy with Theorem~\ref{K:thm:main} and Corollary~\ref{K2:cor:strict}, in which the fibre count $2^k$ of real planes is the largest possible near a corner. The two-sided bounds of Corollary \ref{N:cor:planes} grow with $n$ like $\binom n{2k}^{1/2}$ and $\binom n{2k}$. For lines the exact growth rate is determined in Sections~\ref{G:s4} and \ref{D:s4}: for random lines the area is asymptotic to $\frac\pi2n\log n$. For planes of dimension $k\ge2$ it remains open. Finally, (\ref{N:e1.1}) can be applied to non-coordinate matrices $B$. For example, sums of coordinates give monomial maps whose images have larger Newton polytopes, and optimising over $B$ is an open problem.

\section{Generic injectivity, lines and space curves}\label{F:s1}

Let $V\subset(\C^*)^n$ be a smooth complete intersection of dimension $k$, $n\ge2k+1$, defined by Laurent polynomials $f_1,\dots,f_{n-k}$ with supports $A_1,\dots,A_{n-k}\subset\Z^n$ and Newton polytopes $\Delta_j=\conv A_j$. Theorem~\ref{N:thm:main} gives, for irreducible $V$, the bounds
$$\frac1{N_{\max}}\Big(\sum_{|I|=2k}M_I(V)^2\Big)^{1/2}\ \le\ \Hm^{2k}(\A(V))\ \le\ \frac1{N_{\min}}\sum_{|I|=2k}M_I(V).$$
Here $M_I(V)$ is the coamoeba multivolume of the coordinate projection $\pi_I|_V$, counted with multiplicity, and $N_{\max}$, $N_{\min}$ are the essential supremum and infimum of the number $N(x)$ of points of $V$ over $x\in\A(V)$. In this section and in Sections~\ref{F:s2} to \ref{F:s4} we answer four questions left open there and in Section~\ref{N:s9}: whether $N_{\max}$ can be removed from the bounds for generic varieties, how the bounds compare for real and nonreal lines in $(\C^*)^3$, whether the birationality hypothesis of Corollary~\ref{N:cor:curves} holds for generic curves, and how the area of the amoeba of a line grows with $n$. The notation is that of Sections~\ref{N:s1} and \ref{N:s2}, to which we refer for the area formula and the definitions.

Let $\Lambda_j=\Z(A_j-A_j)$ be the lattice spanned by differences of exponents of $f_j$, $\Lambda=\sum_j\Lambda_j$, and
$$H=\{t\in(\C^*)^n:\ t^\lambda=1\ \text{for all}\ \lambda\in\Lambda\}.$$
When $\Lambda$ has rank $n$, $H$ is a finite group of order $[\Z^n:\Lambda]$, consisting of points whose coordinates are roots of unity, and it lies in the unit torus $\T^n$. For $t\in H$ and $\alpha_j\in A_j$ we have $f_j(tz)=t^{\alpha_j}f_j(z)$, so $H$ acts freely on $V$ and preserves $\Log$. Hence $N\ge|H|$ on the amoeba. The coefficients of $f_1,\dots,f_{n-k}$ form a vector $c$ in a complex vector space $\mathcal C$, and we let $\mathcal U\subset\mathcal C$ be the open set of $c$ for which all coefficients indexed by the $A_j$ are nonzero and $V_c$ is a smooth complete intersection of dimension $k$.

\begin{theorem}\label{F:thm:injective}
Suppose that every $\Delta_j$ is $n$-dimensional and $n\ge2k+1$. For Lebesgue-almost every $c\in\mathcal U$, the following holds for $\Hm^{2k}$-almost every $x\in\A(V_c)$: the points of $V_c$ over $x$ form a single orbit of $H$, so that $N(x)=|H|$. In particular, if $\Lambda=\Z^n$, then $\Log|_{V_c}$ is injective outside a set whose image is $\Hm^{2k}$-null.
\end{theorem}

The proof shows more: the conclusion holds for supports that are not full-dimensional, provided $\operatorname{rank}\Lambda=n$ and $\operatorname{rank}\Lambda_S>2k+2|S|-n$ for every nonempty proper subset $S$ of the equations, where $\Lambda_S=\sum_{j\in S}\Lambda_j$ (Remark \ref{F:rem:general}).

\begin{corollary}\label{F:cor:noNmax}
Under the hypotheses of Theorem \ref{F:thm:injective}, for almost every $c$ with $V_c$ irreducible,
$$\frac1{|H|}\Big(\sum_IM_I(V_c)^2\Big)^{1/2}\ \le\ \Hm^{2k}(\A(V_c))\ \le\ \frac1{|H|}\sum_IM_I(V_c).$$
For almost every $k$-plane in $(\C^*)^n$, $n\ge2k+1$, one has $H=\{1\}$, $M_I=\pi^{2k}$ for all $I$, and $\binom n{2k}^{1/2}\pi^{2k}\le\Hm^{2k}(\A(V))\le\binom n{2k}\pi^{2k}$.
\end{corollary}

For lines the structure is completely explicit. A line in general position is $\ell_q=\{(\rho-q_1,\dots,\rho-q_n):\rho\in\C\setminus\{q_j\}\}$ up to multiplicative translation, with distinct $q_j$. It is real up to multiplicative translation exactly when the $q_j$ are collinear.

\begin{theorem}\label{F:thm:lines3}
Every real line in general position in $(\C^*)^3$ has $\Hm^2(\A(\ell))\le\frac32\pi^2$. Every nonreal line in general position in $(\C^*)^3$ has $\Hm^2(\A(\ell))\ge\sqrt3\,\pi^2$. More generally, in $(\C^*)^n$, $n\ge3$, every nonreal line in general position has $\Hm^2(\A(\ell))\ge\binom n2^{1/2}\pi^2$, and every real line has $\Hm^2(\A(\ell))\le\frac12\binom n2\pi^2$.
\end{theorem}

So in $(\C^*)^3$ the infimum of the amoeba area over lines in general position is an infimum over real lines. Conjecture~\ref{F:conj:lines3} below predicts where it is attained.

\begin{theorem}\label{F:thm:birational}
Let $d_1,\dots,d_{n-1}\ge1$, $D=\prod_jd_j$, and let $V\subset(\C^*)^n$ be the zero set of $n-1$ polynomials with Newton polytopes $d_jS_n$ and generic coefficients, with $n\ge3$. Then $V$ is a smooth irreducible curve. For every pair $I$ of coordinates, $\pi_I|_V$ is birational onto a plane curve whose reduced equation has Newton polygon $D\cdot\conv(0,e_1,e_2)$. Consequently, for generic coefficients, the bounds of Corollary~\ref{N:cor:curves} hold without further hypothesis, $\binom n2^{1/2}\pi^2/(2D^2)\le\Hm^2(\A(V))\le\binom n2D^2\pi^2$. Moreover, for almost every choice of coefficients, $N=1$ almost everywhere, so that in fact $\binom n2^{1/2}\pi^2\le\Hm^2(\A(V))\le\binom n2D^2\pi^2$.
\end{theorem}

\section{Generic injectivity}\label{F:s2}\label{F:sec:injective}

For a subset $S$ of $\{1,\dots,n-k\}$ put $\Lambda_S=\sum_{j\in S}\Lambda_j$, with $\Lambda_\emptyset=0$, and $H_S=\{t\in\T^n:t^\lambda=1\ \forall\lambda\in\Lambda_S\}$, a closed subgroup of $\T^n$ of real dimension $n-\operatorname{rank}\Lambda_S$. For $t\in\T^n$ let $S(t)=\{j:t^\lambda=1\ \forall\lambda\in\Lambda_j\}$, and put $T'_S=\{t\in\T^n:S(t)=S\}$. The sets $T'_S$ partition $\T^n$. Each $T'_S$ is an open subset of $H_S$, namely $H_S$ minus the closed subgroups $H_{S\cup\{j\}}$, $j\notin S$, so it is a manifold of dimension at most $n-\operatorname{rank}\Lambda_S$. Finally $T'_{\{1,\dots,n-k\}}=H\cap\T^n=H$ when $\Lambda$ has rank $n$.

\begin{lemma}\label{F:lem:samemoduli}
Two points $z,z'\in(\C^*)^n$ satisfy $\Log z=\Log z'$ if and only if $z'=tz$ for a unique $t\in\T^n$. If $z\in V_c$ and $t\in T'_S$, then $tz\in V_c$ if and only if $f_{c,j}(tz)=0$ for all $j\notin S$.
\end{lemma}

\begin{proof}
The first assertion is clear, with $t=z'/z$ coordinatewise. If $j\in S$, then $t^{\alpha}=t^{\alpha'}$ for all $\alpha,\alpha'\in A_j$, so $f_{c,j}(tz)=t^{\alpha_j}f_{c,j}(z)=0$ for any $\alpha_j\in A_j$.
\end{proof}

Let $X=\{(z,c)\in(\C^*)^n\times\mathcal U:f_{c,j}(z)=0,\ j=1,\dots,n-k\}$. The map $(z,c)\mapsto(f_{c,j}(z))_j$ is a submersion, since varying the coefficient of a monomial $z^\alpha$, $\alpha\in A_j$, changes $f_{c,j}(z)$ by a nonzero multiple of $z^\alpha$. So $X$ is a smooth manifold of real dimension $\dim_\R\mathcal U+2k$, and the projection $X\to\mathcal U$ has fibres $V_c$. For each $S$ consider
$$\Psi_S:X\times T'_S\to\C^{S^c},\qquad\Psi_S(z,c,t)=\big(f_{c,j}(tz)\big)_{j\notin S}.$$

\begin{lemma}\label{F:lem:submersion}
$\Psi_S$ is a submersion.
\end{lemma}

\begin{proof}
Let $(z,c,t)\in X\times T'_S$ and $j\notin S$. Since $j\notin S(t)$, there are $\alpha,\alpha'\in A_j$ with $t^\alpha\neq t^{\alpha'}$. For $x\in\C$, consider the variation of $c$ which adds $x\,z^{-\alpha}$ to the coefficient of $z^\alpha$ in $f_j$, subtracts $x\,z^{-\alpha'}$ from the coefficient of $z^{\alpha'}$, and leaves everything else unchanged. It changes $f_{c,j}(z)$ by $x-x=0$, so together with $\delta z=0$ and $\delta t=0$ it is tangent to $X\times T'_S$. It changes $f_{c,j}(tz)$ by $x(t^\alpha-t^{\alpha'})$, which takes every complex value as $x$ varies. It does not change $f_{c,j'}$ for $j'\neq j$. So the differential of $\Psi_S$ is surjective.
\end{proof}

\begin{proof}[Proof of Theorem \ref{F:thm:injective}]
By Lemma \ref{F:lem:submersion}, $Y_S=\Psi_S^{-1}(0)$ is a submanifold of $X\times T'_S$ of real codimension $2|S^c|$, of dimension $\dim_\R\mathcal U+2k+\dim T'_S-2(n-k-|S|)$. Apply Sard's theorem to the projection $Y_S\to\mathcal U$. For $c$ outside a null set $\mathcal N_S$, the fibre $Z_S(c)=\{(z,t):z\in V_c,\ t\in T'_S,\ \Psi_S(z,c,t)=0\}$ is a submanifold of dimension
$$\delta_S=2k+\dim T'_S-2(n-k-|S|).$$
By Lemma \ref{F:lem:samemoduli}, $Z_S(c)=\{(z,t):z\in V_c,\ t\in T'_S,\ tz\in V_c\}$. We claim that $\delta_S<2k$ for every $S$ other than $\{1,\dots,n-k\}$. For $S=\emptyset$, $\dim T'_\emptyset=n$ and $\delta_\emptyset=4k-n<2k$ since $n>2k$. For $S\neq\emptyset$, each $\Lambda_j$ has rank $n$ because $\Delta_j$ is $n$-dimensional, so $H_S$ is finite, $\dim T'_S=0$ and $\delta_S=2k-2(n-k-|S|)<2k$ when $|S|<n-k$.

Let $c\notin\bigcup_S\mathcal N_S$, and let $E_c\subset V_c$ be the union of the projections of the manifolds $Z_S(c)$, $S\neq\{1,\dots,n-k\}$, to $V_c$. Each is a smooth image of a manifold of dimension less than $2k$, so $\Hm^{2k}(E_c)=0$, and $\Hm^{2k}(\Log(E_c))=0$ since $\Log$ is locally Lipschitz. Let $x\in\A(V_c)\setminus\Log(E_c)$ and $z,z'\in V_c$ over $x$. Then $z'=tz$ with $t\in T'_{S(t)}$, and $(z,t)\in Z_{S(t)}(c)$. Since $z\notin E_c$, $S(t)=\{1,\dots,n-k\}$ and $t\in H$. Conversely $Hz\subset V_c$ lies over $x$, and $H$ acts freely. So the fibre is $Hz$ and $N(x)=|H|$.
\end{proof}

\begin{remark}\label{F:rem:general}
Without full-dimensionality, and assuming $\operatorname{rank}\Lambda=n$, the same proof works if $\delta_S<2k$ for all $S\neq\{1,\dots,n-k\}$, that is, if $n-\operatorname{rank}\Lambda_S<2(n-k-|S|)$, or $\operatorname{rank}\Lambda_S>2k+2|S|-n$, for every nonempty proper $S$. If $\Lambda$ has rank less than $n$, then $H$ is not compact, the subgroup $H\cap\T^n$ preserving $\Log$ is a compact torus of positive dimension, and it preserves $V_c$. The fibres of $\Log|_{V_c}$ are then infinite, and the amoeba has dimension less than $2k$. The dimension of amoebas is determined in general by Draisma, Rau and Yuen [DRY].
\end{remark}

\begin{proof}[Proof of Corollary \ref{F:cor:noNmax}]
By Theorem \ref{F:thm:injective}, $N_{\max}=N_{\min}=|H|$, and the bounds of Theorem~\ref{N:thm:main} give the first statement. For $k$-planes, take $n-k$ affine linear equations with full support $\{0,e_1,\dots,e_n\}$, whose convex hull is the $n$-dimensional simplex, and $\Lambda=\Z^n$. Every $k$-plane in general position can be defined by such equations, by the argument of Section~\ref{K2:s3}, which carries over to $n>2k$. Indeed, if the row space $R\subset\C^{n+1}$ of the coefficient matrix lay in $\{y_0=0\}$, the plane would pass through the origin, whose closure lies on the $n\ge k+1$ hyperplanes $\lambda_1,\dots,\lambda_n$. If $R$ lay in $\{y_i=0\}$ with $i\ge1$, the point at infinity in the direction $e_i$ would lie on $\lambda_0$ and on the $n-1\ge k+1$ hyperplanes $\lambda_l$, $l\ne i$. Either way general position would fail, so a generic basis of $R$ has all entries nonzero. The set of coefficient vectors defining planes in general position is open and dense. The map from coefficient matrices of rank $n-k$ to the Grassmannian of affine $k$-planes is a smooth submersion, and the bad set of coefficients is saturated, since the fibre count depends only on the plane. So a null set of coefficients corresponds to a null set of planes for the smooth measure class of the Grassmannian, which is the meaning of almost every $k$-plane. So Theorem \ref{F:thm:injective} applies, $H=\{1\}$, and Corollary~\ref{N:cor:planes} gives $M_I=\pi^{2k}$ for all $I$.
\end{proof}

\section{Lines}\label{F:s3}\label{F:sec:lines}

For the line $\ell_q$ put $w_j=1/(\rho-q_j)$. The differential of $\rho\mapsto\Log\ell_q(\rho)$, in the real coordinates $(\operatorname{Re}\rho,\operatorname{Im}\rho)$, has columns $a=\operatorname{Re}w$ and $-b$, where $b=\operatorname{Im}w$. By the Lagrange identity, $J^2=|a|^2|b|^2-(a\cdot b)^2$. Since $\sum_jw_j^2=|a|^2-|b|^2+2i\,a\cdot b$ and $\sum_j|w_j|^2=|a|^2+|b|^2$,
$$J(\rho)=\frac12\Big(\big(\textstyle\sum_j|\rho-q_j|^{-2}\big)^2-\big|\sum_j(\rho-q_j)^{-2}\big|^2\Big)^{1/2}.\eqtag{F:e3.1}$$
The Jacobian of a coordinate projection is $J_{\{i,j\}}=|a_ib_j-a_jb_i|=|\operatorname{Im}(\bar w_iw_j)|$.

\begin{lemma}\label{F:lem:lineN}
Let $n\ge3$. If the $q_j$ are not collinear, $\Log|_{\ell_q}$ is injective. If they are collinear, almost every fibre over the amoeba has exactly two points, exchanged by the reflection of $\rho$ in the line through the $q_j$. For every pair $I$, $M_I(\ell_q)=\pi^2$.
\end{lemma}

The dichotomy between real and nonreal lines is also in [MN2, Remark 7.1] and in [NP, Section 4.1].

\begin{proof}
If $|\rho-q_j|=|\rho'-q_j|$ for all $j$ with $\rho\neq\rho'$, then all $q_j$ lie on the perpendicular bisector of $[\rho,\rho']$. This proves the first two assertions, the real points $\rho$ on the line of the $q_j$ having a one-dimensional image. The projection $\pi_I(\ell_q)$ is the line $\{(\rho-q_i,\rho-q_j)\}$ of $(\C^*)^2$ with $q_i\neq q_j$, a multiplicative translate of $\{1+u+v=0\}$ after a monomial change of coordinates, and $\pi_I$ is injective on $\ell_q$. By Lemma~\ref{L:lem:line} its coamoeba multivolume is $\pi^2$.
\end{proof}

\begin{proof}[Proof of Theorem \ref{F:thm:lines3}]
For a nonreal line, $N=1$ almost everywhere by Lemma \ref{F:lem:lineN}, so Theorem~\ref{N:thm:main} gives $\Hm^2(\A(\ell))=\int J\ge(\sum_IM_I^2)^{1/2}=\binom n2^{1/2}\pi^2$. For a real line, $N=2$ almost everywhere, so $\Hm^2(\A(\ell))=\frac12\int J\le\frac12\sum_IM_I=\frac12\binom n2\pi^2$. This upper bound is also the case $p=2$, $n_j=1$ of [MN1, Theorem 4.1], and it appears for a particular real line among the examples following that theorem. For $n=3$ these are $\sqrt3\,\pi^2$ and $\frac32\pi^2<\sqrt3\,\pi^2$.
\end{proof}

For $n=3$ and real lines, an affine change of $\rho$ normalises $q=(0,1,t)$ with $t\in\R\setminus\{0,1\}$. The three values $t=-1,\frac12,2$ are the harmonic configurations: they are exchanged by the affine symmetries permuting $0,1,t$, and in each of them one of the three points is the midpoint of the other two. They are the fixed points of the involutions of the configuration space induced by these symmetries, and therefore critical points of the area along real lines, whenever the area is differentiable there. We propose the following conjecture, studied further in Section~\ref{G:sec:three}.

\begin{conjecture}\label{F:conj:lines3}
The minimum of $\Hm^2(\A(\ell))$ over lines in general position in $(\C^*)^3$ is attained exactly at the real harmonic lines.
\end{conjecture}

The growth of the area with $n$ for lines in $(\C^*)^n$ is studied in Sections~\ref{G:s4} and \ref{D:s4}.

\section{Space curves}\label{F:s4}\label{F:sec:curves}

\begin{proof}[Proof of Theorem \ref{F:thm:birational}]
Let $\bar C\subset\PP^n$ be the common zero set of the homogenisations $F_1,\dots,F_{n-1}$ of the $f_j$, forms of degrees $d_j$. For generic coefficients, $\bar C$ is a smooth complete intersection curve of degree $D$ by Bertini's theorem [Ha, Theorem II.8.18]. It is connected, as every complete intersection of positive dimension is, hence irreducible. Moreover it meets every coordinate linear subspace of codimension at least $2$ in the empty set, and every coordinate hyperplane transversally. Indeed, each of these finitely many conditions defines a proper closed subset of the space of coefficients, since it fails for a generic curve by a count of dimensions. So $V=\bar C\cap(\C^*)^n$ is a smooth irreducible curve.

Fix $I=\{a,b\}$ and let $L_I=\{x_0=x_a=x_b=0\}$, a linear subspace of dimension $n-3$. Since $\bar C\cap L_I=\emptyset$, the projection $p_I:\bar C\to\PP^2$, $x\mapsto[x_0:x_a:x_b]$, is a finite morphism, which extends $\pi_I|_V$. We show that it is birational for generic coefficients. Let $G$ be the Grassmannian of $(n-3)$-planes in $\PP^n$, and let $\mathcal B\subset\mathcal C\times G$ be the set of pairs $(c,L)$ such that $\bar C_c$ is a smooth irreducible complete intersection curve, $\bar C_c\cap L=\emptyset$, and the projection from $L$ is not birational onto its image. It is a constructible set, and by definition $\bar C_c$ is a smooth irreducible curve for every $(c,L)\in\mathcal B$. In particular the coefficients for which $\bar C_c$ is smooth but not an irreducible curve, for instance $F_1=F_2$ when $d_1=d_2$, do not occur in $\mathcal B$. For each $c$ for which $\bar C_c$ is a smooth irreducible curve, the fibre $\mathcal B_c$ is a constructible subset of $G$ with $\dim\mathcal B_c<\dim G$. This is the classical fact that a projection from a generic centre of an irreducible curve to $\PP^2$ is birational onto its image. It is obtained by projecting from generic points successively: a generic point projection $\PP^m\dashrightarrow\PP^{m-1}$ of a curve is birational onto its image for $m\ge3$ [Ha, Section IV.3]. Projecting successively from $n-2$ points amounts to projecting from their linear span, and the map from $(n-2)$-tuples of points to $G$ given by the span is dominant. The set of centres for which the projection is birational is constructible and contains the image of a dense set of tuples, so it contains a nonempty open subset of $G$. Hence $\dim\mathcal B\le\dim\mathcal C+\dim G-1$. The group $GL_{n+1}(\C)$ acts on $\mathcal C$, by linear changes of the homogeneous coordinates in the forms, and on $G$, transitively, and $\mathcal B$ is invariant. So the projection $\mathcal B\to G$ is an equivariant map onto a homogeneous space, all its fibres $\mathcal B^L=\{c:(c,L)\in\mathcal B\}$ are isomorphic, and the theorem on the dimension of fibres gives $\dim\mathcal B=\dim G+\dim\mathcal B^{L_I}$. Therefore $\dim\mathcal B^{L_I}\le\dim\mathcal C-1$, and for generic $c$ the projection $p_I$ is birational.

The image $p_I(\bar C)$ is then an irreducible plane curve of degree $D$, since $p_I$ is finite and birational and $\bar C$ has degree $D$. Let $g$ be its reduced equation, a form of degree $D$ in $(x_0,x_a,x_b)$. The coefficient of $x_0^D$ in $g$ is nonzero if and only if the point $[1:0:0]$ is not on the image, that is, if $\bar C$ does not meet $\{x_a=x_b=0\}$. This is a coordinate linear subspace of codimension $2$, so it holds. In the same way the coefficients of $x_a^D$ and $x_b^D$ are nonzero, because $\bar C$ does not meet $\{x_0=x_b=0\}$ and $\{x_0=x_a=0\}$. In the affine coordinates $(z_a,z_b)=(x_a/x_0,x_b/x_0)$, the Newton polygon of $g$ is therefore $D\cdot\conv(0,e_1,e_2)$. Its zero set in $(\C^*)^2$ is the closure of $\pi_I(V)$, since $\bar C$ meets the coordinate hyperplanes in finitely many points. This proves the first part, and the bounds of Corollary~\ref{N:cor:curves} follow.

For the last statement, the supports of the $f_j$ are all lattice points of $d_jS_n$, so each $\Lambda_j=\Z^n$. Theorem \ref{F:thm:injective} applies with $k=1$, $n\ge3$, and $H=\{1\}$, so $N=1$ almost everywhere for almost every choice of coefficients. Then (\ref{N:e1.2}) with $N_{\max}=1$ and $M_I\ge\pi^2$ gives $\Hm^2(\A(V))\ge\binom n2^{1/2}\pi^2$.
\end{proof}

\section{Amoebas of lines: statements}\label{G:s1}\label{G:sec:intro}

We now study lines in detail, continuing Section~\ref{F:sec:lines}. Let $q=(q_1,\dots,q_n)$ be $n\ge3$ distinct points of $\C$. As there, the line $\ell_q=\{(\rho-q_1,\dots,\rho-q_n):\rho\in\C\setminus\{q_1,\dots,q_n\}\}$ of $(\C^*)^n$ is in general position, and up to multiplicative translation, which only translates the amoeba, every line of $(\C^*)^n$ in general position is of this form, as recalled before Theorem~\ref{F:thm:lines3}. We write $w_j=(\rho-q_j)^{-1}$, $A=\sum_j|w_j|^2$ and $B=\sum_jw_j^2$. By formula~(\ref{F:e3.1}), the Jacobian of $\Log|_{\ell_q}$ with respect to Lebesgue measure in the parameter $\rho$ is $J=\frac12(A^2-|B|^2)^{1/2}$. For a pair $a<b$ we put $J_{ab}=|\operatorname{Im}(\bar w_aw_b)|$, the Jacobian of the projection of $\Log|_{\ell_q}$ to the coordinates $a,b$, and we recall that $\int_\C J_{ab}\,dA=\pi^2$, by the proof of Theorem~\ref{F:thm:lines3} and Lemma~\ref{L:lem:line}. The line $\ell_q$ is called real if the points $q_j$ lie on a real affine line of $\C$, and nonreal otherwise. By Lemma~\ref{F:lem:lineN}, the map $\Log|_{\ell_q}$ is injective almost everywhere for nonreal lines and two-to-one almost everywhere for real lines, so that $\Hm^2(\A(\ell_q))=\int_\C J\,dA$ in the first case and $\Hm^2(\A(\ell_q))=\frac12\int_\C J\,dA$ in the second. All areas below are areas in $\R^n$, and all integrals over $\C$ are with respect to Lebesgue measure $dA$ in the variable $\rho$.

Theorem~\ref{F:thm:lines3} shows that in $(\C^*)^3$ the amoeba of every real line has area at most $\frac32\pi^2$, while the amoeba of every nonreal line has area at least $\sqrt3\pi^2$, and Conjecture~\ref{F:conj:lines3} predicts that the smallest area is attained at the harmonic configurations. The present section and Sections~\ref{G:s2} to \ref{G:s5} address three questions left open there: whether the harmonic real line is the minimiser in $(\C^*)^3$, whether the growth for random lines is of order $n\log n$, and whether the lower bound $\binom n2^{1/2}\pi^2$ for nonreal lines can be improved when $n$ is large.

For three collinear distinct points, an affine map of $\C$ sends the two outer points to $0$ and $1$ and the middle one to a point $s\in(0,1)$. Put $g(s)=\pi^{-2}\Hm^2(\A(\ell_{(0,s,1)}))$. Our first result is the following.

\begin{theorem}\label{G:thm:A}
The function $g$ is continuous on $(0,1)$ and satisfies $g(s)=g(1-s)$. It tends to $\frac{1+\sqrt2}2$ as $s\to0$ and as $s\to1$. Every real line in general position in $(\C^*)^3$ has area $\pi^2g(s)$ for some $s\in(0,\frac12]$, every nonreal line has area at least $\sqrt3\pi^2$ by Theorem~\ref{F:thm:lines3}, and $\frac{\sqrt3}2\le g\le\frac32$ on $(0,1)$. Consequently the infimum of $\Hm^2(\A(\ell))$ over all lines in general position in $(\C^*)^3$ equals $\pi^2\inf_{(0,1/2]}g$, and it lies between $\frac{\sqrt3}2\pi^2$ and $\frac{1+\sqrt2}2\pi^2$. If $g(s_0)<\frac{1+\sqrt2}2$ for some $s_0$, then this infimum is attained, and it is attained only by real lines.
\end{theorem}

The equispaced configuration $s=\frac12$ is the fixed point of the symmetry $s\mapsto1-s$, and so it is a critical point of $g$ whenever $g$ is differentiable there. Its orbit under the affine symmetries consists of the harmonic configurations of Section~\ref{F:sec:lines}. The conjecture that the harmonic line is the minimiser is therefore equivalent to the statement that $\min_{(0,1/2]}g=g(\frac12)$, and it would follow from the monotonicity of $g$ on $(0,\frac12]$ (Conjecture~\ref{G:conj:mono}). Even the attainment of the infimum depends on the inequality $\inf g<\frac{1+\sqrt2}2$, which we have not proved; it would follow, for instance, from any upper bound for $g(\frac12)$ smaller than $\frac{1+\sqrt2}2$.

Our second result concerns random lines. Let $q_1,\dots,q_n$ be independent random points with the uniform distribution on the unit disc $D=\{|z|<1\}$. With probability one the points are distinct and not collinear, so that $\ell_q$ is a nonreal line in general position and $\Hm^2(\A(\ell_q))=\int J$.

\begin{theorem}\label{G:thm:B}
For $n\ge6$, the expected area satisfies $\frac\pi2n\log n-\frac\pi2\sqrt n\log n-\pi(\frac52-\log2)n\le\E\,\Hm^2(\A(\ell_q))\le\frac\pi2n\log n+9\pi n\sqrt{1+\log n}+\pi(\log6+\frac12+\frac{11}{\sqrt2})n$. Moreover $\Hm^2(\A(\ell_q))/(\frac\pi2n\log n)$ tends to $1$ in probability and in $L^1$ as $n\to\infty$.
\end{theorem}

In units of $\pi^2$, the area of a random line is thus asymptotic to $\frac1{2\pi}n\log n$. It lies far above the lower bound $\binom n2^{1/2}\approx n/\sqrt2$ of Theorem~\ref{F:thm:lines3} and far below the upper bound $\binom n2$. This settles the growth question of Section~\ref{F:sec:lines} for this model, with the exact constant, and it shows that for almost all nonreal lines in this model the lower bound of Theorem~\ref{F:thm:lines3} is off by a factor of order $\log n$. Our third result concerns all nonreal lines.

\begin{theorem}\label{G:thm:C}
For every nonreal line $\ell$ in general position in $(\C^*)^n$, $n\ge3$, the inequality $\Hm^2(\A(\ell))>\binom n2^{1/2}\pi^2$ is strict. For every $\varepsilon>0$, the probability that a line through $n$ independent uniform points of the unit disc satisfies $\Hm^2(\A(\ell))\ge(\frac\pi2-\varepsilon)n\log n$ tends to $1$ as $n\to\infty$.
\end{theorem}

We conjecture that a lower bound of order $n\log n$ holds for every line, real or not, and Section~\ref{G:s5} gives three pieces of evidence: the random model, the equispaced real lines, for which a heuristic computation gives the same leading term $\frac\pi2n\log n$, and a degeneration formula, proved here for $n=3$ and in general in Theorem~\ref{D:thm:D}, under which every limit of real configurations along a binary tree of clusters has area at least $\frac{\pi^2}{4\log2}n\log n$, with equality for balanced binary trees when $n$ is a power of $2$.

\section{Identities and inequalities for the Jacobian}\label{G:s2}\label{G:sec:identities}

\begin{lemma}\label{G:lem:J}
For every $\rho\in\C\setminus\{q_1,\dots,q_n\}$ the following hold. First, $J^2=\sum_{a<b}J_{ab}^2$. Second, $J_{ab}=|\sin\theta_{ab}|/(|\rho-q_a|\,|\rho-q_b|)$, where $\theta_{ab}$ is the angle at $\rho$ between the segments $[\rho,q_a]$ and $[\rho,q_b]$. Third, if $J_S$ denotes the Jacobian of the line $\ell_{q_S}$ for a subset $S$ of the indices with $|S|\ge2$, then $J_S\le J$. Fourth, $\frac12(A-|B|)\le J\le\frac12A$. Fifth, $J_{ab}\le|w_a|\,|w_b|$ and $J_{ab}\le(|\rho|\,|q_a-q_b|+|q_a|\,|q_b|)/(|\rho-q_a|^2|\rho-q_b|^2)$.
\end{lemma}

\begin{proof}
For complex numbers $u,v$ we have $|u|^2|v|^2-\operatorname{Re}(u^2\bar v^2)=|u\bar v|^2-\operatorname{Re}((u\bar v)^2)=2(\operatorname{Im}(u\bar v))^2$. Hence $A^2-|B|^2=\sum_{a,b}(|w_a|^2|w_b|^2-\operatorname{Re}(w_a^2\bar w_b^2))=2\sum_{a,b}(\operatorname{Im}(w_a\bar w_b))^2=4\sum_{a<b}J_{ab}^2$, the terms with $a=b$ being zero. This gives the first identity, which is also the Cauchy--Binet formula for the $2\times2$ minors of the real $n\times2$ Jacobian matrix of $\Log|_{\ell_q}$. Writing $\rho-q_a=r_ae^{i\phi_a}$, we get $\bar w_aw_b=e^{i(\phi_a-\phi_b)}/(r_ar_b)$, whence the second formula with $\theta_{ab}=\phi_a-\phi_b$. The third assertion follows from the first, since $J_S^2$ is the sum of the terms $J_{ab}^2$ with $a,b\in S$. For the fourth, $|B|\le A$, and then $A-|B|\le(A^2-|B|^2)^{1/2}=((A-|B|)(A+|B|))^{1/2}\le A$. For the fifth, the first bound is clear, and $\operatorname{Im}(\bar w_aw_b)=\operatorname{Im}((\rho-q_a)\overline{(\rho-q_b)})/(|\rho-q_a|^2|\rho-q_b|^2)$, with $\operatorname{Im}((\rho-q_a)\overline{(\rho-q_b)})=\operatorname{Im}(\bar\rho(q_b-q_a))+\operatorname{Im}(q_a\bar q_b)$, whose modulus is at most $|\rho|\,|q_a-q_b|+|q_a|\,|q_b|$.
\end{proof}

\begin{lemma}\label{G:lem:inv}
The area $\Hm^2(\A(\ell_q))$ is unchanged when the points $q_j$ are permuted, when $q$ is replaced by $\alpha q+\beta=(\alpha q_1+\beta,\dots,\alpha q_n+\beta)$ with $\alpha\in\C^*$ and $\beta\in\C$, and when $q$ is replaced by $\bar q$.
\end{lemma}

\begin{proof}
A permutation of the $q_j$ permutes the coordinates of $(\C^*)^n$ and of $\R^n$, which is an isometry of $\R^n$. The substitution $\rho=\alpha\rho'+\beta$ shows that $\ell_{\alpha q+\beta}=\alpha\cdot\ell_q$ is the multiplicative translate of $\ell_q$ by $(\alpha,\dots,\alpha)$, whose amoeba is the translate of $\A(\ell_q)$ by $\log|\alpha|(1,\dots,1)$. Finally $\ell_{\bar q}$ is the complex conjugate of $\ell_q$, which has the same amoeba.
\end{proof}

\begin{lemma}\label{G:lem:cont}
The functions $q\mapsto\int_\C J\,dA$ and $q\mapsto\int_\Omega J\,dA$, for a fixed measurable set $\Omega$, are continuous on the open set of configurations of $n$ distinct points of $\C$.
\end{lemma}

\begin{proof}
Let $q^{(m)}\to q$, with all configurations in a compact set $K$ of configurations of distinct points, and let $\eta>0$ be smaller than one third of the minimal distance between two points of any configuration in $K$. Then $J(\rho;q^{(m)})\to J(\rho;q)$ for every $\rho\notin\{q_1,\dots,q_n\}$. We show that the family $J(\cdot;q')$, $q'\in K$, is uniformly integrable, and continuity follows from Vitali's convergence theorem. By Lemma \ref{G:lem:J}, $J\le\sum_{a<b}J_{ab}$. Outside the union of the discs $D(q'_a,\eta)$, the fifth part of Lemma \ref{G:lem:J} gives $J_{ab}\le C\min(1,|\rho|^{-3})$ with $C$ depending only on $K$ and $\eta$, which is integrable and independent of $q'$. Inside $D(q'_a,\eta)$, every $J_{ab}$ with $b\ne a$ is at most $|w_a|\,|w_b|\le\eta^{-1}|\rho-q'_a|^{-1}$, every $J_{bc}$ with $a\notin\{b,c\}$ is at most $\eta^{-2}$, and for every measurable $E\subset D(q'_a,\eta)$ we have $\int_E|\rho-q'_a|^{-1}\,dA\le\int_{D(q'_a,r)}|\rho-q'_a|^{-1}dA=2\pi r$ with $\pi r^2=\operatorname{area}(E)$. Hence $\int_EJ\le C'(\operatorname{area}(E)^{1/2}+\operatorname{area}(E))$ for every $E\subset\C$ with a constant $C'$ depending only on $K$, $\eta$ and $n$, together with the uniform tail bound outside a large disc. This is uniform integrability.
\end{proof}

\section{Real lines in $(\C^*)^3$ and the equispaced line}\label{G:s3}\label{G:sec:three}

\subsection{Reduction to one parameter}\label{G:s3.1}
Let $\ell$ be a real line in general position in $(\C^*)^3$, so $\ell$ is a multiplicative translate of $\ell_q$ with $q_1,q_2,q_3$ distinct and collinear. After a permutation, which does not change the area by Lemma \ref{G:lem:inv}, $q_2$ lies on the segment $[q_1,q_3]$, and the affine map sending $q_1$ to $0$ and $q_3$ to $1$ sends $q_2$ to some $s\in(0,1)$. By Lemma \ref{G:lem:inv} again, $\Hm^2(\A(\ell))=\pi^2g(s)$. The affine map $\rho\mapsto1-\rho$, followed by the permutation exchanging the first and third points, sends $(0,s,1)$ to $(0,1-s,1)$, so $g(s)=g(1-s)$, and every real line has area $\pi^2g(s)$ with $s\in(0,\frac12]$. Continuity of $g$ follows from Lemma \ref{G:lem:cont}. The bounds $\frac{\sqrt3}2\le g\le\frac32$ follow from Theorem~\ref{F:thm:lines3} and its proof, since for a real line $\Hm^2=\frac12\int J$ and $\binom32^{1/2}\pi^2\le\int J\le\sum_{a<b}\int J_{ab}=3\pi^2$, the first inequality being the Minkowski inequality $\int(\sum_{a<b}J_{ab}^2)^{1/2}\ge(\sum_{a<b}(\int J_{ab})^2)^{1/2}$.

The configuration $s=\frac12$ is the equispaced one, $q=(0,\frac12,1)$, in which the middle point is the midpoint of the outer ones. In the normalisation $q=(0,1,t)$ of Section~\ref{F:sec:lines}, it corresponds to the three harmonic values $t=-1,\frac12,2$, for which one of the three points is the midpoint of the other two. The equispaced configuration is the unique fixed point of $s\mapsto1-s$, so if $g$ is differentiable at $\frac12$, then $g'(\frac12)=0$.

\subsection{The degenerate limit}\label{G:s3.2}

\begin{proposition}\label{G:prop:limit}
As $t\to+\infty$, $\int_\C J\,dA$ for $q=(0,1,t)$ tends to $(1+\sqrt2)\pi^2$. Consequently $g(s)\to\frac{1+\sqrt2}2$ as $s\to0$ and as $s\to1$.
\end{proposition}

\begin{proof}
Let $t>4$ and split $\C$ into $\Omega_1=\{|\rho|\le\sqrt t\}$ and $\Omega_2=\{|\rho|>\sqrt t\}$. On $\Omega_1$ we have $J_{12}\le J\le J_{12}+J_{13}+J_{23}$ by the first part of Lemma \ref{G:lem:J}. The function $J_{12}$ does not depend on $t$, and $\int_{\Omega_1}J_{12}\to\int_\C J_{12}=\pi^2$. For $a\in\{1,2\}$ and $\rho\in\Omega_1$ we have $|\rho-t|\ge t-\sqrt t\ge t/2$, so $J_{a3}\le|w_a|\,|w_3|\le2t^{-1}|\rho-q_a|^{-1}$, and since $\Omega_1\subset D(q_a,\sqrt t+1)$, we get $\int_{\Omega_1}J_{a3}\le2t^{-1}\cdot2\pi(\sqrt t+1)$, which tends to $0$. Hence $\int_{\Omega_1}J\to\pi^2$.

On $\Omega_2$ we substitute $\rho=tu$. Then $w_j(\rho)=t^{-1}\hat w_j(u)$ with $\hat w_j(u)=(u-p_j)^{-1}$ and $p=(0,t^{-1},1)$, so that $J_{ab}(\rho)=t^{-2}\hat J_{ab}(u)$, where $\hat J_{ab}$ is the Jacobian of the pair $(p_a,p_b)$, and $dA_\rho=t^2dA_u$. Hence $\int_{\Omega_2}J\,dA=\int_{|u|>t^{-1/2}}\hat J\,dA$, with $\hat J=(\hat J_{12}^2+\hat J_{13}^2+\hat J_{23}^2)^{1/2}$. For $|u|\ge t^{-1/2}$ we have $|u|\ge2/t$ and so $|u-t^{-1}|\ge|u|/2$. The fifth part of Lemma \ref{G:lem:J} gives $\hat J_{12}\le(|u|t^{-1})/(|u|^2|u-t^{-1}|^2)\le4t^{-1}|u|^{-3}$, and therefore $\int_{|u|>t^{-1/2}}\hat J_{12}\le4t^{-1}\cdot2\pi\int_{t^{-1/2}}^\infty r^{-2}dr=8\pi t^{-1/2}$, which tends to $0$. Let $K$ denote the Jacobian of the pair $(0,1)$, so that $\hat J_{13}=K$ and $\int_\C K=\pi^2$. The function $\hat J_{23}$, the Jacobian of the pair $(t^{-1},1)$, tends to $K$ at every $u\notin\{0,1\}$. For $|u|\ge t^{-1/2}$, Lemma \ref{G:lem:J} gives $\hat J_{23}\le|\hat w_2||\hat w_3|\le2|u|^{-1}|u-1|^{-1}$ and, since by the identity in the proof of Lemma \ref{G:lem:J} the numerator of $\hat J_{23}$ is $|\operatorname{Im}(\bar u(1-t^{-1}))|\le|u|$, the term $\operatorname{Im}(p_2\bar p_3)$ vanishing for the real points $p_2,p_3$, also $\hat J_{23}\le|u|/(|u-t^{-1}|^2|u-1|^2)\le4|u|^{-1}|u-1|^{-2}$. The same two bounds hold for $K$. The function $G(u)=\min(2|u|^{-1}|u-1|^{-1},4|u|^{-1}|u-1|^{-2})$ is integrable on $\C$, since it is $O(|u|^{-1})$ near $0$, $O(|u-1|^{-1})$ near $1$, and $O(|u|^{-3})$ at infinity. By dominated convergence, $\int_{|u|>t^{-1/2}}(\hat J_{13}^2+\hat J_{23}^2)^{1/2}dA\to\int_\C\sqrt2K\,dA=\sqrt2\pi^2$. Since $|\hat J-(\hat J_{13}^2+\hat J_{23}^2)^{1/2}|\le\hat J_{12}$, we get $\int_{\Omega_2}J\to\sqrt2\pi^2$, and so $\int_\C J\to(1+\sqrt2)\pi^2$.

For the consequence, the configuration $(0,s,1)$ is, by Lemma \ref{G:lem:inv}, equivalent to $(0,1,t)$ with $t=1/s$, through $\rho\mapsto\rho/s$. The line is real, so $g(s)=\frac12\pi^{-2}\int J$ tends to $\frac{1+\sqrt2}2$ as $s\to0$, and by the symmetry also as $s\to1$.
\end{proof}

The limit has a transparent meaning. In the region $|\rho|\lesssim1$ the third coordinate $\rho-t$ is nearly constant, and the amoeba looks like the amoeba of the line $\{(\rho,\rho-1)\}$ of $(\C^*)^2$ times a point, whose coamoeba multivolume is $\pi^2$. In the region $|\rho|\sim t$ the first two coordinates nearly coincide, and the amoeba looks like the amoeba of a line of $(\C^*)^2$ placed in the plane $x_1=x_2$ of $\R^3$, through the map $(x,y)\mapsto(x,x,y)$, which multiplies areas by $\sqrt2$.

\begin{proof}[Proof of Theorem \ref{G:thm:A}]
The first three sentences are proved above. A nonreal line has area at least $\sqrt3\pi^2>\frac32\pi^2\ge\pi^2\sup g$ by Theorem~\ref{F:thm:lines3}, so the infimum over all lines is $\pi^2\inf_{(0,1/2]}g$, and nonreal lines are never minimisers. If $g(s_0)<\frac{1+\sqrt2}2$ for some $s_0\in(0,\frac12]$, then by Proposition \ref{G:prop:limit} there is $s_1\in(0,s_0)$ with $g>g(s_0)$ on $(0,s_1]$, and the continuous function $g$ attains its infimum over the compact interval $[s_1,\frac12]$, which is its infimum over $(0,\frac12]$. The upper bound $\frac{1+\sqrt2}2\pi^2$ for the infimum follows from the limit of $g$ at $0$, and the lower bound from $g\ge\frac{\sqrt3}2$.
\end{proof}

\subsection{The monotonicity conjecture}\label{G:s3.3}
We state the following conjecture, which implies Conjecture~\ref{F:conj:lines3}.

\begin{conjecture}\label{G:conj:mono}
The function $g$ is strictly decreasing on $(0,\frac12]$. In particular, among all lines in general position in $(\C^*)^3$, the amoeba has minimal area exactly for the real lines through three equispaced points.
\end{conjecture}

A proof would need either a closed formula for $g$, which we have not found, or a sign estimate for the derivative $g'$ on $(0,\frac12)$, combined with a rigorous version of the rate of convergence in Proposition \ref{G:prop:limit}, or a convexity argument for $\int J$ along the family. The derivative cannot be computed by differentiating under the integral sign, since $\partial_sJ$ has a nonintegrable singularity of order $|\rho-s|^{-2}$ at the moving point; one must first move the singularity to a fixed position by a diffeomorphism of the plane that equals a translation near the moving point and the identity near the other points. The proof of Proposition \ref{G:prop:limit} gives explicit error terms, of size $O(t^{-1/2})$, for all pieces except the dominated convergence of $\hat J_{23}$, which can also be made explicit.

\section{Random lines: growth of order $n\log n$}\label{G:s4}\label{G:sec:random}

In this section $q_1,\dots,q_n$ are independent and uniform on the unit disc $D$, and $n\ge6$. Expected multivolumes of amoebas of random half-dimensional complete intersections are computed by Bayraktar and Ki\c{s}isel [BK]; here we study the area of the amoeba itself, for random lines. We write $\log$ for the natural logarithm. Since $\Hm^2(\A(\ell_q))=\int_\C J$ almost surely, Tonelli's theorem gives $\E\,\Hm^2(\A(\ell_q))=\int_\C\E J(\rho)\,dA(\rho)$.

\begin{lemma}\label{G:lem:beurling}
For $\rho\in D$ and $0<\delta<1-|\rho|$, $\int_{D\setminus D(\rho,\delta)}(\rho-q)^{-2}dA(q)=0$.
\end{lemma}

This is the classical fact that the Beurling transform of the indicator function of a disc vanishes inside the disc, see [AIM, Chapter 4]; we include the short proof.

\begin{proof}
On the region $\Omega=D\setminus\overline{D(\rho,\delta)}$ the function $F(q)=(\rho-q)^{-1}$ is smooth up to the boundary and $\partial F/\partial q=(\rho-q)^{-2}$. Since $d(F\,d\bar q)=\frac{\partial F}{\partial q}dq\wedge d\bar q$ and $dq\wedge d\bar q=-2i\,dA$, Stokes' theorem gives $\int_\Omega\frac{\partial F}{\partial q}\,dA=\frac i2\int_{\partial\Omega}F\,d\bar q$, with the standard orientation of $\partial\Omega$. On the unit circle $\bar q=1/q$ and $d\bar q=-q^{-2}dq$, so the integral is $\oint_{|q|=1}\frac{dq}{q^2(q-\rho)}$, whose residues at $q=0$ and $q=\rho$ are $-\rho^{-2}$ and $\rho^{-2}$ if $\rho\ne0$, with sum $0$. If $\rho=0$ the integrand is $q^{-3}$, whose residue at $0$ is $0$. On the circle $|q-\rho|=\delta$ we have $\bar q-\bar\rho=\delta^2/(q-\rho)$, so $d\bar q=-\delta^2(q-\rho)^{-2}dq$ and the integral is $\oint\delta^2(q-\rho)^{-3}dq=0$.
\end{proof}

\subsection{The lower bound}\label{G:s4.1}
Fix $\delta=\frac1{2\sqrt n}$ and $\rho$ with $s=|\rho|<1-\delta$. Let $S(\rho)$ be the set of indices $l$ with $|q_l-\rho|\ge\delta$, and let $A_S$ and $B_S$ be the sums defining $A$ and $B$ restricted to $l\in S(\rho)$. If $|S(\rho)|\ge2$, Lemma \ref{G:lem:J} gives $J\ge J_{S(\rho)}\ge\frac12(A_S-|B_S|)$, and if $|S(\rho)|\le1$ then $A_S-|B_S|=0\le J$. So $J\ge\frac12(A_S-|B_S|)$ always. We have $A_S=\sum_lX_l$ and $B_S=\sum_lY_l$ with $X_l=|\rho-q_l|^{-2}\mathbf 1_{|\rho-q_l|\ge\delta}$ and $Y_l=(\rho-q_l)^{-2}\mathbf 1_{|\rho-q_l|\ge\delta}$ independent and identically distributed. The disc $D(\rho,1-s)$ is contained in $D$, so $\E X_l\ge\frac1\pi\int_{\delta\le|x|\le1-s}|x|^{-2}dA(x)=2\log\frac{1-s}\delta$. By Lemma \ref{G:lem:beurling}, $\E Y_l=0$, and $\E|Y_l|^2\le\frac1\pi\int_{|x|\ge\delta}|x|^{-4}dA=\delta^{-2}$, so $\E|B_S|\le(\E|B_S|^2)^{1/2}=(n\E|Y_1|^2)^{1/2}\le\sqrt n/\delta$. Hence $\E J(\rho)\ge n\log\frac{1-s}\delta-\frac{\sqrt n}{2\delta}=\frac12n\log n+n\log2-n+n\log(1-s)$. Integrating over the disc $|\rho|<1-\delta$, whose area is at least $\pi(1-2\delta)=\pi-\pi/\sqrt n$, and using $\int_D\log(1-|\rho|)\,dA=2\pi\int_0^1s\log(1-s)\,ds=-\frac{3\pi}2$ together with the fact that $\log(1-s)\le0$, we obtain $\E\int J\ge\frac\pi2n\log n-\frac\pi2\sqrt n\log n-\pi(1-\log2)n-\frac{3\pi}2n$, which is the lower bound of Theorem \ref{G:thm:B}.

\subsection{The upper bound}\label{G:s4.2}
For $|\rho|>2$, every pair satisfies $J_{ab}\le(2|\rho|+1)/(|\rho|-1)^4$ by the fifth part of Lemma \ref{G:lem:J}, so $J\le\binom n2^{1/2}(2|\rho|+1)/(|\rho|-1)^4\le\frac n{\sqrt2}(2|\rho|+1)/(|\rho|-1)^4$, and $\int_{|\rho|>2}J\le\frac n{\sqrt2}2\pi\int_2^\infty\frac{r(2r+1)}{(r-1)^4}dr=\frac{11\pi}{\sqrt2}n$, because $\int_1^\infty(2x^2+5x+3)x^{-4}dx=2+\frac52+1=\frac{11}2$.

Inside $D(0,2)$, let $d_i$ be the distance from $q_i$ to the nearest other point, and $D_i=D(q_i,d_i/2)$. These discs are pairwise disjoint and contained in $D(0,2)$. Outside their union, $J\le\frac12A$ by Lemma \ref{G:lem:J}. Inside $D_i$, put $A'_i=\sum_{l\ne i}|w_l|^2$. Lemma \ref{G:lem:J} gives $J^2=\sum_{l\ne i}J_{il}^2+J_{(i)}^2\le|w_i|^2A'_i+\frac14A_i'^2$, where $J_{(i)}$ is the Jacobian of the configuration without $q_i$, so $J\le|w_i|\sqrt{A'_i}+\frac12A'_i$. Each term $|w_l|^2$ of $A'_i$ is integrated over a region disjoint from $D_l$. Therefore $\int_{D(0,2)}J\le\frac12\sum_l\int_{D(0,2)\setminus D_l}|\rho-q_l|^{-2}dA+\sum_i\int_{D_i}|w_i|\sqrt{A'_i}\,dA$. The region $D(0,2)\setminus D_l$ lies in the annulus $\frac{d_l}2\le|\rho-q_l|\le3$, so the first sum is at most $\pi\sum_l\log(6/d_l)$. For $\rho\in D_i$ and $l\ne i$ we have $|\rho-q_l|\ge|q_i-q_l|-\frac{d_i}2\ge\frac12|q_i-q_l|$, so $A'_i\le4T_i$ with $T_i=\sum_{l\ne i}|q_i-q_l|^{-2}$, and $\int_{D_i}|w_i|\,dA=\pi d_i$. The second sum is thus at most $2\pi\sum_id_i\sqrt{T_i}$.

It remains to estimate the expectations. Since $\PP(|q_m-q_l|<x)\le x^2$, we get $\PP(d_l<x)\le(n-1)x^2$, and $\E\log(1/d_l)\le\int_0^\infty\min(1,(n-1)e^{-2u})\,du=\frac12\log(n-1)+\frac12$. So the expectation of the first sum is at most $\frac\pi2n\log n+\pi(\log6+\frac12)n$. For the second, fix $i$ and $l\ne i$, and let $e$ be the distance from $q_i$ to the nearest of the other $n-2$ points. Then $d_i\le\min(e,|q_i-q_l|)$, and conditionally on $q_i$ the variables $e$ and $x=|q_i-q_l|$ are independent. For $y\le1$ the disc $D(q_i,y)$ meets $D$ in a set containing the disc $D(q_i(1-\frac y2),\frac y2)$, since $|q_i|\frac y2+\frac y2\le y$ and $|q_i|(1-\frac y2)+\frac y2\le1$, so $\PP(e>y\mid q_i)\le(1-\frac{y^2}4)^{n-2}\le e^{-(n-2)y^2/4}$, while $e\le2$. Hence $\E(e^2\mid q_i)\le\int_0^1e^{-(n-2)\sigma/4}d\sigma+3e^{-(n-2)/4}\le\frac{16}{n-2}=:c$. Since $\sigma\mapsto\min(1,\sigma/x^2)$ is concave, Jensen's inequality gives $\E(d_i^2/x^2\mid q_i,q_l)\le\min(1,c/x^2)$, and since $q_l$ has density at most $\frac1\pi$ and $x\le2$, $\E\min(1,c/x^2)\le c+2c\log(2/\sqrt c)=c(1+\log(4/c))$. For $n\ge6$ we have $c\le4$, where $c(1+\log(4/c))$ is increasing in $c$, so $\E(d_i^2T_i)\le(n-1)c(1+\log\frac4c)\le20(1+\log n)$. Then $\E(d_i\sqrt{T_i})\le(\E\,d_i^2T_i)^{1/2}\le4.5\sqrt{1+\log n}$, and the expectation of the second sum is at most $9\pi n\sqrt{1+\log n}$. Adding the three contributions gives the upper bound of Theorem \ref{G:thm:B}.

\subsection{Convergence in probability and in mean}\label{G:s4.3}
Let $L_n=\frac\pi2n\log n$ and $X_n=\Hm^2(\A(\ell_q))/L_n$. With $\delta$ and $S(\rho)$ as in the lower bound, $\Hm^2(\A(\ell_q))\ge\frac12\sum_lZ_l-\frac12\int_{|\rho|<1-\delta}|B_S|\,dA$, where $Z_l=\int_{|\rho|<1-\delta,\,|\rho-q_l|\ge\delta}|\rho-q_l|^{-2}dA$. The $Z_l$ are independent, identically distributed, with values in $[0,2\pi\log(2/\delta)]$, and $\frac12\sum_l\E Z_l\ge L_n-O(\sqrt n\log n+n)$ by the computation of the lower bound. So $\operatorname{Var}(\frac12\sum Z_l)\le n\pi^2\log^2(4\sqrt n)$, and Chebyshev's inequality gives $\PP(\frac12\sum Z_l<(1-\frac\varepsilon2)L_n)=O(\varepsilon^{-2}n^{-1})$. By the lower bound computation $\E\int|B_S|\le\pi\sqrt n/\delta=2\pi n$, and Markov's inequality gives $\PP(\frac12\int|B_S|>\frac\varepsilon2L_n)\le\frac{2\pi n}{\varepsilon L_n}=\frac{4}{\varepsilon\log n}$. Hence $\PP(X_n<1-\varepsilon)\to0$ for every $\varepsilon>0$. Together with $\E X_n\le1+o(1)$, from the upper bound, this gives $\E|X_n-1|=\E(X_n-1)+2\E(1-X_n)^+\le o(1)+2\varepsilon+2\PP(X_n<1-\varepsilon)$, so $\E|X_n-1|\to0$, and $X_n\to1$ in probability and in $L^1$. This completes the proof of Theorem \ref{G:thm:B}.

\begin{remark}\label{G:rem:heur}
The constant $\frac\pi2$ has a simple origin. At a point $\rho$ of the bulk, the points within distance $r$ of $\rho$ number about $nr^2$, so $A(\rho)\approx\int_{n^{-1/2}}^1\frac{2n\,r\,dr}{r^2}=n\log n$, while $B(\rho)$ is a sum of terms with uniformly distributed phases and is of size $O(n)$, dominated by the nearest point. So $J\approx\frac12A\approx\frac12n\log n$ on a disc of area $\pi$. The same heuristic applies to any bounded density $p$ with compact support, with $A(\rho)\approx\pi np(\rho)\log n$, and gives the same constant $\frac\pi2$, since $\int p=1$. For the standard Gaussian density the same asymptotics holds, as proved in Theorem~\ref{D:thm:E}, which also treats real random lines. The vanishing of Lemma \ref{G:lem:beurling} is convenient but not essential: the proof only needs $|\E Y_l|=O(1)$, which holds whenever the truncated Beurling transform of the density is bounded, as it is for compactly supported H\"older continuous densities, by the classical estimates for the Beurling transform [AIM, Chapter 4].
\end{remark}

\section{Lower bounds for nonreal lines when $n$ is large}\label{G:s5}\label{G:sec:nonreal}

\begin{proof}[Proof of Theorem \ref{G:thm:C}]
A line in general position in $(\C^*)^n$ is not contained in a translate of a two-dimensional subtorus, so the first statement is the case $k=1$ of Corollary~\ref{T:cor:strict}, with $N=1$ for nonreal lines (Lemma~\ref{F:lem:lineN}) and $M_I=\pi^2$ for every pair $I$. The proof of Theorem~\ref{T:thm:lower} also shows that $\int J>\binom n2^{1/2}\pi^2$ for every line in general position, real or not. The second statement is the convergence in probability of Theorem \ref{G:thm:B}, since a line through independent uniform points is nonreal with probability one.
\end{proof}

The strictness is not quantitative, and a deterministic lower bound of order $n\log n$ for all lines seems to need new ideas. The proof of Theorem \ref{G:thm:B} uses randomness only once, to show that $|B_S|$ is small compared with $A_S$ on average, through the cancellation of the phases of the terms $(\rho-q_l)^{-2}$. The deterministic part of the argument gives, for every configuration and every measurable choice of the sets $S(\rho)$, the lower bound $\Hm^2(\A(\ell_q))\ge\frac1{2N}\int(A_S-|B_S|)\,dA$, with $N=1$ or $2$, so that a lower bound of order $n\log n$ follows for every configuration on which the phases cancel in a quantitative way. We collect three pieces of evidence for the following conjecture.

\begin{conjecture}\label{G:conj:nlogn}
There is a constant $c>0$ such that every line in general position in $(\C^*)^n$ satisfies $\Hm^2(\A(\ell))\ge c\,n\log n$.
\end{conjecture}

The first piece of evidence is Theorem \ref{G:thm:B}, which shows that the conjecture holds with any $c<\frac\pi2$ for all lines outside a set of small probability in the random model.

The second is the family of equispaced real lines, $q=(1,2,\dots,n)$, which are the natural candidates for small amoebas since real lines have $N=2$. At a point $\rho=x+iy$ with $1\ll|y|\ll n$ and $0<x<n$, the sum $A$ is approximately $\int_\R\frac{dq}{(x-q)^2+y^2}=\frac\pi{|y|}$, while $B$ is approximately $\int_\R\frac{dq}{(x+iy-q)^2}=0$. So $J\approx\frac\pi{2|y|}$, and integrating over this region gives $\int J\approx\frac\pi2\cdot n\cdot2\log n=\pi n\log n$, hence $\Hm^2(\A(\ell_q))=\frac12\int J\approx\frac\pi2n\log n$, with the same constant as for random lines. By the heuristic, the area of the equispaced real line is of order $n\log n$, and not of order $n$, so that its ratio to the lower bound $\binom n2^{1/2}\pi^2$ valid for nonreal lines would tend to infinity.

The third piece of evidence comes from degenerations. Proposition \ref{G:prop:limit} is the simplest case of a degeneration formula. Suppose that the points split into clusters at very different scales, organised as a rooted tree, with the point at infinity of $\PP^1$ playing the role of the root, each internal node $\nu$ having children which are clusters of sizes $m_1,\dots,m_c$. Near the scale of the node $\nu$, the children are seen as single points with multiplicities $m_1,\dots,m_c$, and the coordinates of the points of one child nearly coincide. A cluster of $m$ coinciding coordinates maps by $x\mapsto(x,\dots,x)$, so the Jacobian of a pair of children of sizes $m_a,m_b$ is multiplied by $\sqrt{m_am_b}$. Arguing as in the proof of Proposition \ref{G:prop:limit}, one expects $\int J$ to tend to the sum over internal nodes of $\int(\sum_{a<b}m_am_bJ_{ab}^2)^{1/2}$, computed for the configuration of the children. We prove this here for $n=3$, and in general in Theorem~\ref{D:thm:D}. For $n=3$ it gives $(1+\sqrt2)\pi^2=(\sqrt{1\cdot1}+\sqrt{2\cdot1})\pi^2$. For instance, for the real configuration $q=(0,\varepsilon,1,1+\varepsilon)$, with two clusters of two points, it gives $\int J\to(\sqrt{1\cdot1}+\sqrt{1\cdot1}+\sqrt{2\cdot2})\pi^2=4\pi^2$. For a binary tree, every node contributes $\sqrt{m_1m_2}\pi^2$. If $f(m)$ denotes the minimum over binary trees with $m$ leaves of the sum of these contributions, divided by $\pi^2$, then $f(1)=0$ and $f(m)=\min_{0<a<m}(\sqrt{a(m-a)}+f(a)+f(m-a))$. The inequality $-p\log_2p-(1-p)\log_2(1-p)\le2\sqrt{p(1-p)}$ for the binary entropy, which follows from Topsøe's bound $H(p)\le(4p(1-p))^{1/\ln4}$ [To] since $1/\ln4>\frac12$, gives by induction $f(m)\ge\frac m2\log_2m$, with equality when $m$ is a power of $2$ and the tree is balanced. Since for a real line $\Hm^2=\frac12\int J$, limits of real configurations along binary trees have area at least $\frac{\pi^2}4n\log_2n=\frac{\pi^2}{4\log2}n\log n$ in the limit, again of order $n\log n$. A caterpillar tree, in which one point is added at a time, gives $\sum_{m<n}\sqrt m\,\pi^2$, of order $n^{3/2}$. Thus no degeneration seems to produce an area of order $n$.

These considerations suggest that the bound $\binom n2^{1/2}\pi^2$ of Theorem~\ref{F:thm:lines3}, of order $n$, is not of the right order for any line when $n$ is large, and that the optimal constant in Conjecture \ref{G:conj:nlogn} is close to $\frac\pi2$.

\section{Degenerations and other random models: statements}\label{D:s1}\label{D:sec:intro}

In this section and in Sections~\ref{D:s2} to \ref{D:s5} we complete the study of lines begun in Section~\ref{F:sec:lines} and continued in Sections~\ref{G:s1} to \ref{G:s5}, and we use the notation of Section~\ref{G:s1}, which we recall. For $n\ge2$ distinct points $q=(q_1,\dots,q_n)$ of $\C$, the line $\ell_q=\{(\rho-q_1,\dots,\rho-q_n):\rho\in\C\setminus\{q_j\}\}$ is a line in general position in $(\C^*)^n$. With $w_j=(\rho-q_j)^{-1}$, $A=\sum_j|w_j|^2$ and $B=\sum_jw_j^2$, the Jacobian of $\Log|_{\ell_q}$ with respect to Lebesgue measure in $\rho$ is $J=\frac12(A^2-|B|^2)^{1/2}$, and $J^2=\sum_{a<b}J_{ab}^2$ with $J_{ab}=|\operatorname{Im}(\bar w_aw_b)|$ (Lemma~\ref{G:lem:J}). Each pair contributes $\int_\C J_{ab}\,dA=\pi^2$. We write $I(q)=\int_\C J\,dA$. By Lemma~\ref{F:lem:lineN}, $\Hm^2(\A(\ell_q))=I(q)$ if the $q_j$ are not collinear and $\Hm^2(\A(\ell_q))=\frac12I(q)$ if they are collinear, for $n\ge3$.

Proposition~\ref{G:prop:limit} shows, for $n=3$, that $I(0,1,t)\to(1+\sqrt2)\pi^2$ as $t\to\infty$, and Section~\ref{G:sec:nonreal} proposes a general degeneration formula on heuristic grounds. Theorem~\ref{G:thm:B} shows that for $n$ independent uniform points of the unit disc the area is asymptotic to $\frac\pi2n\log n$, and Remark~\ref{G:rem:heur} suggests that the same holds for other distributions. The present section and Sections~\ref{D:s2} to \ref{D:s5} prove both.

\subsection{The degeneration formula}\label{D:s1.1}
To describe nested clusters we use rooted trees. Let $T$ be a finite rooted tree whose leaves are labelled by $\{1,\dots,n\}$ and whose internal vertices all have at least two children. For a vertex $v$ let $m(v)$ be the number of leaves below $v$, so that $m(v)=1$ for a leaf. For an internal vertex $v$ with children $c_1,\dots,c_k$ and a configuration $p=(p(c_1),\dots,p(c_k))$ of $k$ distinct points of $\C$, put $J^{(v)}_p=(\sum_{a<b}m(c_a)m(c_b)J_{ab}(p)^2)^{1/2}$, where $J_{ab}(p)$ is the Jacobian of the pair $(p(c_a),p(c_b))$, and $I_v(p)=\int_\C J^{(v)}_p\,dA$.

\begin{definition}\label{D:def:datum}
A degeneration datum on $T$ consists of real numbers $h(v)$ for the internal vertices, with $h(\text{root})=0$ and $h(c)>h(v)$ whenever $c$ is an internal child of $v$, and of continuous maps $\varepsilon\mapsto p_v(\varepsilon)\in\C^{k_v}$, $\varepsilon\in[0,\varepsilon_0)$, one for each internal vertex $v$ with $k_v$ children, such that the $k_v$ points of $p_v(0)$ are distinct. The associated family of configurations is $q_j(\varepsilon)=\sum_v\varepsilon^{h(v)}p_v(\varepsilon)(c_v(j))$, the sum over the internal vertices $v$ on the path from the root to the leaf $j$, where $c_v(j)$ is the child of $v$ on that path.
\end{definition}

For small $\varepsilon>0$ the points $q_j(\varepsilon)$ are distinct. The leaves below a child $c$ of the root form a cluster of diameter $O(\varepsilon^{h(c)})$ around the point $p_{\rm root}(\varepsilon)(c)$, and inside it the same structure repeats at smaller scales.

\begin{theorem}\label{D:thm:D}
For every degeneration datum, $\lim_{\varepsilon\to0}I(q(\varepsilon))=\sum_vI_v(p_v(0))$, the sum over the internal vertices of $T$. If every $p_v(0)$ consists of collinear points and the family is real, so that all the lines $\ell_{q(\varepsilon)}$ are real, the areas converge to half of this sum. If the family is nonreal for small $\varepsilon$, the areas converge to the sum itself.
\end{theorem}

The integral $I_v(p)$ has a direct meaning: $J^{(v)}_p$ is the Jacobian of the map $\rho\mapsto(\log|\rho-p(c_a)|)_a$ composed with the linear map $\R^k\to\R^{\sum m(c_a)}$ that repeats the $a$-th coordinate $m(c_a)$ times, which multiplies the Jacobian of the pair $(a,b)$ by $\sqrt{m(c_a)m(c_b)}$. When $v$ has two children of sizes $m,m'$, $I_v=\sqrt{mm'}\,\pi^2$, whatever the two points. When all the children of $v$ are leaves, $I_v(p)=I(p)$.

\begin{corollary}\label{D:cor:binary}
If $T$ is binary, the limit in Theorem \ref{D:thm:D} is $\pi^2\sum_v\sqrt{m(c_v)m(c'_v)}$, where $c_v,c'_v$ are the children of $v$, and it is at least $\frac{\pi^2}2n\log_2n$. Equality holds when $n$ is a power of $2$ and $T$ is the balanced binary tree. For the caterpillar tree, in which every internal vertex has a leaf among its children, the limit is $\pi^2\sum_{m=1}^{n-1}\sqrt m$, of order $\frac23n^{3/2}\pi^2$. For an arbitrary tree, the limit is at least $\pi^2\sum_ve_2(v)^{1/2}$, where $e_2(v)=\sum_{a<b}m(c_a)m(c_b)$.
\end{corollary}

For $n=3$ and the tree with root children $\{1,2\}$ and $3$, the formula gives $(\sqrt{1\cdot1}+\sqrt{2\cdot1})\pi^2=(1+\sqrt2)\pi^2$, which is Proposition~\ref{G:prop:limit}.

\subsection{Gaussian and real random lines}\label{D:s1.2}
\begin{theorem}\label{D:thm:E}
Let $q_1,\dots,q_n$ be independent with the standard complex Gaussian density $p(z)=\pi^{-1}e^{-|z|^2}$, and $n\ge16$. Then $\frac\pi2n\log n-6\pi n\log\log n\le\E\,\Hm^2(\A(\ell_q))\le\frac\pi2n\log n+50\sqrt2\pi n\sqrt{3+\log n}+\frac\pi2n\log(1+\log n)+\pi(\log6+\frac12+\frac{11}{\sqrt2})n$, and $\Hm^2(\A(\ell_q))/(\frac\pi2n\log n)\to1$ in probability and in $L^1$.
\end{theorem}

\begin{theorem}\label{D:thm:F}
Let $q_1,\dots,q_n$ be independent and uniform on $[0,1]\subset\R\subset\C$, so that $\ell_q$ is a real line, and $n\ge16$. Then $\frac\pi2n\log n-2n\log\log n-(\pi+\sqrt{\pi/2})n\le\E\,\Hm^2(\A(\ell_q))\le\frac\pi2n\log n+\frac\pi2(\log12+1+4\sqrt2+\frac{11}{\sqrt2})n$, and $\Hm^2(\A(\ell_q))/(\frac\pi2n\log n)\to1$ in probability and in $L^1$.
\end{theorem}

The constant $\frac\pi2$ is thus the same in all three models, although its origin differs. For points in the plane, $I(q)\approx\frac12\int A$ and the logarithm comes from the nearest-neighbour distance $n^{-1/2}$, which gives $\frac12\log n$ per point and a total $\frac\pi2n\log n$ for the complex lines, whose area is $I$. For points on a line, the nearest-neighbour distance is $n^{-1}$, which gives $\log n$ per point and $I\approx\pi n\log n$, but the area of a real line is $\frac12I$. The factor two in the logarithm and the factor $\frac12$ from the fibre count cancel exactly. This explains the observation of Section~\ref{G:sec:nonreal} that equispaced real lines and random complex lines have the same leading term.

\section{Preliminary estimates}\label{D:s2}\label{D:sec:prelim}

We use freely the identities of Lemma~\ref{G:lem:J}: $J^2=\sum_{a<b}J_{ab}^2$; $J_S\le J$ for the Jacobian $J_S$ of a subconfiguration; $\frac12(A-|B|)\le J\le\frac12A$; $J_{ab}\le|w_a||w_b|$. We also need three simple facts.

\begin{lemma}\label{D:lem:pairs}
For distinct $q_a,q_b$ and $\rho\notin\{q_a,q_b\}$, $J_{ab}(\rho)\le|q_a-q_b|/(|\rho-q_a|\,|\rho-q_b|^2)$.
\end{lemma}

\begin{proof}
We have $J_{ab}=|\operatorname{Im}((\rho-q_a)\overline{(\rho-q_b)})|/(|\rho-q_a|^2|\rho-q_b|^2)$, and $(\rho-q_a)\overline{(\rho-q_b)}=|\rho-q_a|^2+(\rho-q_a)\overline{(q_a-q_b)}$, whose imaginary part has modulus at most $|\rho-q_a|\,|q_a-q_b|$.
\end{proof}

\begin{lemma}\label{D:lem:scaling}
For $\alpha\in\C^*$ and $\beta\in\C$, the Jacobian of the configuration $\alpha q+\beta$ at $\alpha u+\beta$ equals $|\alpha|^{-2}$ times the Jacobian of $q$ at $u$, so that $\int_{\alpha\Omega+\beta}J(\cdot;\alpha q+\beta)\,dA=\int_\Omega J(\cdot;q)\,dA$ for every measurable $\Omega$. The same holds for the weighted Jacobians $J^{(v)}_p$.
\end{lemma}

\begin{proof}
Each $w_j$ is multiplied by $\alpha^{-1}$, so each $J_{ab}$ by $|\alpha|^{-2}$, and $dA$ by $|\alpha|^2$.
\end{proof}

\begin{lemma}\label{D:lem:far}
If $|q_j|\le M$ for all $j$ and $R\ge2M>0$, then $\int_{|\rho|>R}J\,dA\le36\sqrt2\,\pi\,nM/R$ and $\int_{|\rho|>2M}J\,dA\le\frac{11\pi}{\sqrt2}n$.
\end{lemma}

\begin{proof}
By Lemma~\ref{G:lem:J}, $J_{ab}\le(|\rho||q_a-q_b|+|q_a||q_b|)/(|\rho-q_a|^2|\rho-q_b|^2)\le(2M|\rho|+M^2)/(|\rho|-M)^4$ for $|\rho|>M$, so $J\le\binom n2^{1/2}(2Mr+M^2)/(r-M)^4\le\frac n{\sqrt2}(2Mr+M^2)/(r-M)^4$ with $r=|\rho|$. For $r\ge2M$, $r-M\ge r/2$, and $\int_R^\infty2\pi r\cdot16(2Mr+M^2)r^{-4}dr=32\pi(2M/R+M^2/(2R^2))\le32\pi\cdot\frac94M/R=72\pi M/R$, using $M/R\le\frac12$. Multiplying by $\frac n{\sqrt2}$ gives the first bound. For the second, the substitution $\rho=Mv$ and the computation of Section~\ref{G:s4.2} give $\frac n{\sqrt2}2\pi\int_2^\infty\frac{r(2r+1)}{(r-1)^4}dr=\frac{11\pi}{\sqrt2}n$.
\end{proof}

\section{Proof of the degeneration formula}\label{D:s3}\label{D:sec:degen}

We argue by induction on the number of internal vertices of $T$. If the root is the only internal vertex, all its children are leaves, $q(\varepsilon)=p_{\rm root}(\varepsilon)\to p_{\rm root}(0)$, and the statement is the continuity of $I$ on configurations of distinct points (Lemma~\ref{G:lem:cont}).

In general, let $c_1,\dots,c_k$ be the children of the root and $P_a(\varepsilon)=p_{\rm root}(\varepsilon)(c_a)$. Let $K_a$ be the set of leaves below $c_a$. If $c_a$ is a leaf, $K_a=\{j\}$ and $q_j(\varepsilon)=P_a(\varepsilon)$. If $c_a$ is internal, let $T_a$ be the subtree rooted at $c_a$, and put $r_j(\varepsilon)=\sum_v\varepsilon^{h(v)-h(c_a)}p_v(\varepsilon)(c_v(j))$ for $j\in K_a$, the sum over the internal vertices of $T_a$ on the path from $c_a$ to $j$. Then $q_j(\varepsilon)=P_a(\varepsilon)+\varepsilon^{h(c_a)}r_j(\varepsilon)$, and $r(\varepsilon)=(r_j(\varepsilon))_{j\in K_a}$ is the family associated with the degeneration datum on $T_a$ given by the heights $h(v)-h(c_a)$ and the same maps $p_v$. It has fewer internal vertices, so by induction $I(r(\varepsilon))\to\sum_{v\in T_a}I_v(p_v(0))$. Moreover there is $C\ge1$ with $|r_j(\varepsilon)|\le C$ for all $j$, all internal $c_a$ and all small $\varepsilon$.

Let $d>0$ with $|P_a(\varepsilon)-P_b(\varepsilon)|\ge4d$ for $a\ne b$ and all small $\varepsilon$, which exists since the points $P_a(0)$ are distinct and $p_{\rm root}$ is continuous, and let $M$ bound all $|q_j(\varepsilon)|$. For an internal child $c_a$ put $\delta_a=\varepsilon^{h(c_a)/2}$ and $\Omega_a=D(P_a(\varepsilon),\delta_a)$, and let $\Omega_0$ be the complement of the union of these discs. For $\varepsilon$ small, $\delta_a<d$, $2C\varepsilon^{h(c_a)}\le\delta_a$, and the discs $\Omega_a$ are disjoint. We study the integral of $J=J(\cdot;q(\varepsilon))$ over each region.

\subsection{The inner regions}\label{D:s3.1}
Fix an internal child $c_a$, and let $J_{K_a}$ be the Jacobian of the subconfiguration $(q_j)_{j\in K_a}$. By Lemma~\ref{G:lem:J}, $J^2=J_{K_a}^2+\sum'J_{jl}^2$, the sum $\sum'$ over pairs not contained in $K_a$, so that $0\le J-J_{K_a}\le\sum'J_{jl}$. On $\Omega_a$, every point $q_l$ with $l\notin K_a$ satisfies $|\rho-q_l|\ge4d-\delta_a-C\varepsilon^{h(c_l)}\ge2d$, where the last term is absent if $c_l$ is a leaf. Hence for $j\in K_a$ and $l\notin K_a$, $J_{jl}\le|w_j||w_l|\le(2d)^{-1}|\rho-q_j|^{-1}$, and since $\Omega_a\subset D(q_j,2\delta_a)$, $\int_{\Omega_a}J_{jl}\le(2d)^{-1}\cdot2\pi\cdot2\delta_a$. For $j,l\notin K_a$, $J_{jl}\le(2d)^{-2}$ on $\Omega_a$, and $\int_{\Omega_a}J_{jl}\le(2d)^{-2}\pi\delta_a^2$. Both tend to $0$, so $\int_{\Omega_a}(J-J_{K_a})\to0$.

By Lemma \ref{D:lem:scaling}, with $\rho=P_a(\varepsilon)+\varepsilon^{h(c_a)}u$, $\int_{\Omega_a}J_{K_a}\,dA=\int_{|u|<\varepsilon^{-h(c_a)/2}}J(u;r(\varepsilon))\,dA(u)=I(r(\varepsilon))-\int_{|u|\ge\varepsilon^{-h(c_a)/2}}J(u;r(\varepsilon))\,dA(u)$. By Lemma \ref{D:lem:far} with the bound $C$ and $R=\varepsilon^{-h(c_a)/2}\ge2C$, the last integral is at most $36\sqrt2\pi\,m(c_a)\,C\,\varepsilon^{h(c_a)/2}$, which tends to $0$. Hence $\int_{\Omega_a}J\to\sum_{v\in T_a}I_v(p_v(0))$.

\subsection{The outer region}\label{D:s3.2}
On $\Omega_0$, consider first a pair $j,l$ in the same cluster $K_a$, with $c_a$ internal. Then $|q_j-q_l|\le2C\varepsilon^{h(c_a)}$, and for $\rho\in\Omega_0$, $|\rho-q_j|\ge|\rho-P_a|-C\varepsilon^{h(c_a)}\ge\frac12|\rho-P_a|$, since $|\rho-P_a|\ge\delta_a\ge2C\varepsilon^{h(c_a)}$, and likewise for $q_l$. Lemma \ref{D:lem:pairs} gives $J_{jl}\le2C\varepsilon^{h(c_a)}\cdot8|\rho-P_a|^{-3}$, and $\int_{|\rho-P_a|\ge\delta_a}16C\varepsilon^{h(c_a)}|\rho-P_a|^{-3}dA=32\pi C\varepsilon^{h(c_a)}/\delta_a=32\pi C\varepsilon^{h(c_a)/2}$, which tends to $0$.

Now let $F_\varepsilon=(\sum''J_{jl}^2)^{1/2}$, the sum $\sum''$ over pairs $j\in K_a$, $l\in K_b$ with $a<b$. Then $|J-F_\varepsilon|\le\sum_a\sum_{j<l\in K_a}J_{jl}$, so $\int_{\Omega_0}|J-F_\varepsilon|\to0$ by the previous paragraph. For $\rho\notin\{P_1(0),\dots,P_k(0)\}$, $q_j(\varepsilon)\to P_a(0)$ for $j\in K_a$, so $J_{jl}(\rho;q(\varepsilon))\to J_{ab}(\rho;p_{\rm root}(0))$ for $j\in K_a$, $l\in K_b$, and since there are $m(c_a)m(c_b)$ such pairs, $F_\varepsilon(\rho)\to J^{(\rm root)}_{p_{\rm root}(0)}(\rho)$. Also $\mathbf 1_{\Omega_0}\to1$ almost everywhere, since $\delta_a\to0$. It remains to show that the functions $F_\varepsilon\mathbf 1_{\Omega_0}$ are uniformly integrable, and then Vitali's theorem gives $\int_{\Omega_0}F_\varepsilon\to I_{\rm root}(p_{\rm root}(0))$. For $|\rho|\ge2M$, Lemma \ref{D:lem:far} gives a bound $\frac n{\sqrt2}(2M|\rho|+M^2)(|\rho|-M)^{-4}$, independent of $\varepsilon$ and integrable. For $|\rho|<2M$ and $\rho\in\Omega_0$, each $J_{jl}$ with $j\in K_a$, $l\in K_b$, $a\ne b$ satisfies $J_{jl}\le|w_j||w_l|$, where $|w_j|\le2|\rho-P_a(\varepsilon)|^{-1}$ on $\Omega_0$ by the above, and $|w_l|\le2|\rho-P_b(\varepsilon)|^{-1}$. Near $P_a(\varepsilon)$, at distance less than $2d$, the second factor is at most $(1/d)$ and the first is $2|\rho-P_a(\varepsilon)|^{-1}$, and for every measurable set $E$ of area $\pi s^2$ we have $\int_E|\rho-P_a(\varepsilon)|^{-1}dA\le2\pi s$, uniformly in $\varepsilon$. Away from all the $P_a$ at distance at least $2d$, the integrand is bounded by $4d^{-2}$. Hence $\int_EF_\varepsilon\le C'(\operatorname{area}(E)^{1/2}+\operatorname{area}(E))$ for $E\subset D(0,2M)$, with $C'$ independent of $\varepsilon$, which together with the uniform tail bound is uniform integrability.

\subsection{Conclusion of the proof}\label{D:s3.3}
Adding the regions, $I(q(\varepsilon))=\int_{\Omega_0}J+\sum_a\int_{\Omega_a}J\to I_{\rm root}(p_{\rm root}(0))+\sum_a\sum_{v\in T_a}I_v(p_v(0))$, which is the sum over all internal vertices. The assertion on areas follows from Lemma~\ref{F:lem:lineN}, since the fibre count is $2$ for real lines and $1$ for nonreal lines. This proves Theorem \ref{D:thm:D}.

\begin{proof}[Proof of Corollary \ref{D:cor:binary}]
For two children, $J^{(v)}_p=\sqrt{m(c_v)m(c'_v)}J_{12}(p)$ and $\int J_{12}=\pi^2$. Let $f(m)$ be the minimum of $\sum_v\sqrt{m(c_v)m(c'_v)}$ over binary trees with $m$ leaves, so that $f(1)=0$ and $f(m)=\min_{0<a<m}(\sqrt{a(m-a)}+f(a)+f(m-a))$. We show $f(m)\ge\frac m2\log_2m$ by induction. If it holds below $m$, then for $a=pm$ we get $\sqrt{a(m-a)}+f(a)+f(m-a)\ge m\sqrt{p(1-p)}+\frac a2\log_2a+\frac{m-a}2\log_2(m-a)=\frac m2\log_2m+m(\sqrt{p(1-p)}-\frac12H(p))$, with $H(p)=-p\log_2p-(1-p)\log_2(1-p)$. The inequality $H(p)\le(4p(1-p))^{1/\ln4}$ of Topsøe [To], together with $4p(1-p)\le1$ and $1/\ln4>\frac12$, gives $H(p)\le2\sqrt{p(1-p)}$, so the bracket is nonnegative. Equality holds for $p=\frac12$ at every vertex, that is for the balanced tree when $m$ is a power of $2$, since then $f(2^s)\le2^{s-1}+2f(2^{s-1})$ gives $f(2^s)\le s2^{s-1}$. For the caterpillar, the internal vertices have children of sizes $(m,1)$ for $m=1,\dots,n-1$, and $\sum_{m<n}\sqrt m\sim\frac23n^{3/2}$. For an arbitrary tree, Minkowski's integral inequality gives $I_v(p)\ge(\sum_{a<b}m(c_a)m(c_b)(\int J_{ab})^2)^{1/2}=\pi^2e_2(v)^{1/2}$.
\end{proof}

\begin{remark}\label{D:rem:DM}
The configurations of $n$ distinct points of $\C$ modulo affine maps form the moduli space $M_{0,n+1}$ of $n+1$ distinct points of $\PP^1$, the point at infinity being the $(n+1)$-st. Its Deligne--Mumford compactification [DM], [Kn] has boundary strata indexed by stable trees, and the degeneration data of Definition \ref{D:def:datum} are one-parameter families approaching such a stratum, with the point at infinity on the root. Theorem \ref{D:thm:D} says that $I$ has a limit along these families, which depends only on the limit point in the stratum and not on the heights. It is natural to expect that $I$ extends continuously to the whole compactification, with the value given by Theorem \ref{D:thm:D} on each boundary stratum. We have not proved this, since an arbitrary sequence approaching the boundary need not be of the form of Definition \ref{D:def:datum}, although after passing to a subsequence it is close to one.
\end{remark}

\section{Gaussian random lines}\label{D:s4}\label{D:sec:gauss}

Let $p(z)=\pi^{-1}e^{-|z|^2}$, so that $|q|^2$ is exponential with mean $1$ when $q$ has density $p$. Three independent points are collinear with probability zero, so $\Hm^2(\A(\ell_q))=I(q)$ almost surely, and $\E\,\Hm^2=\int_\C\E J(\rho)\,dA(\rho)$ by Tonelli's theorem. We write $\|D^2p\|$ for the supremum over $\C$ of the operator norm of the real Hessian of $p$, which is $2/\pi$, attained at $0$, since the second derivatives of $e^{-x^2}$ are bounded by $2$ in absolute value and the Hessian of $e^{-|z|^2}$ is $e^{-|z|^2}(4zz^{\mathsf T}-2\,\mathrm{Id})$, of norm $\max(2,|4|z|^2-2|)e^{-|z|^2}\le2$.

\begin{lemma}\label{D:lem:beurlingp}
For $\rho\in\C$ and $0<\delta\le1$, the integral $\Phi_\delta(\rho)=\int_{|x|\ge\delta}p(\rho-x)x^{-2}dA(x)$ satisfies $|\Phi_\delta(\rho)|\le1+\frac\pi2\|D^2p\|\le2$.
\end{lemma}

\begin{proof}
Split the domain into $\delta\le|x|<1$ and $|x|\ge1$. On the second part $|x^{-2}|\le1$, and the integral is at most $\int p=1$. On the annulus $\delta\le|x|<1$, the functions $x^{-2}$ and $x^{-2}\,\ell(x)$, for any real linear function $\ell$, have zero integral, because $x^{-2}=r^{-2}e^{-2i\theta}$ and $\ell(x)=r(\alpha\cos\theta+\beta\sin\theta)$, and $\int_0^{2\pi}e^{-2i\theta}d\theta=\int_0^{2\pi}e^{-2i\theta}\cos\theta\,d\theta=\int_0^{2\pi}e^{-2i\theta}\sin\theta\,d\theta=0$. So we may replace $p(\rho-x)$ by $p(\rho-x)-p(\rho)+\nabla p(\rho)\cdot x$, which is at most $\frac12\|D^2p\||x|^2$ by Taylor's formula. The integral over the annulus is then at most $\frac12\|D^2p\|\pi$.
\end{proof}

\subsection{The lower bound}\label{D:s4.1}
Let $n\ge16$, so that $\log\log n\ge1$, and put $R=\sqrt{\log\log n}$, $r=1/\log n$ and $\delta=1/(2\sqrt n)$. Then $\delta<r\le1$. For $|\rho|<R$, let $S(\rho)$ be the set of $l$ with $|q_l-\rho|\ge\delta$, and $A_S$, $B_S$ the corresponding sums. As in Section~\ref{G:s4.1}, $J\ge\frac12(A_S-|B_S|)$. Since $A_S$ is a sum of $n$ independent copies of $X=|\rho-q|^{-2}\mathbf 1_{|\rho-q|\ge\delta}$, we have $\E A_S\ge n\int_{\delta\le|x|\le r}p(\rho-x)|x|^{-2}dA\ge2\pi n\,m_r(\rho)\log\frac r\delta$, where $m_r(\rho)=\min_{D(\rho,r)}p\ge p(\rho)e^{-2|\rho|r-r^2}\ge p(\rho)(1-3Rr)$, because $|y|^2\le|\rho|^2+2|\rho|r+r^2$ on $D(\rho,r)$ and $r\le R$. Next, $\E B_S=n\Phi_\delta(\rho)$, so $|\E B_S|\le2n$ by Lemma \ref{D:lem:beurlingp}, and $\operatorname{Var}B_S\le n\int_{|x|\ge\delta}p(\rho-x)|x|^{-4}dA\le n\pi^{-1}\cdot\pi\delta^{-2}=n/\delta^2$, so $\E|B_S|\le2n+\sqrt n/\delta=4n$. Therefore $\E J(\rho)\ge\pi n\,p(\rho)(1-3Rr)\log(2r\sqrt n)-2n$. Since $\int_{|\rho|<R}p=1-e^{-R^2}=1-1/\log n$, and $\log(2r\sqrt n)=\frac12\log n+\log2-\log\log n$, we obtain
$\E I\ge\pi n(1-\frac1{\log n})(1-\frac{3R}{\log n})(\frac12\log n-\log\log n)-2\pi R^2n$. Expanding, $\E I\ge\frac\pi2n\log n-\frac\pi2n-\frac{3\pi}2nR-\pi n\log\log n-2\pi n\log\log n$, and since $R\le\log\log n$ and $1\le\log\log n$, this is at least $\frac\pi2n\log n-6\pi n\log\log n$.

\subsection{The upper bound}\label{D:s4.2}
Let $M=\max_j|q_j|$. By Lemma \ref{D:lem:far}, $\int_{|\rho|>2M}J\le\frac{11\pi}{\sqrt2}n$. Inside $D(0,2M)$ we argue as in Section~\ref{G:s4.2}. With $d_i$ the distance from $q_i$ to its nearest neighbour and $D_i=D(q_i,d_i/2)$, these discs are disjoint and contained in $D(0,2M)$, and $\int_{D(0,2M)}J\le\frac12\sum_l\int_{D(0,2M)\setminus D_l}|\rho-q_l|^{-2}dA+2\pi\sum_id_i\sqrt{T_i}$, where $T_i=\sum_{l\ne i}|q_i-q_l|^{-2}$. The region $D(0,2M)\setminus D_l$ lies in the annulus $\frac{d_l}2\le|\rho-q_l|\le3M$, and $d_l\le2M$, so the first sum is at most $\pi\sum_l\log(6M/d_l)$. Since $p\le\pi^{-1}$, $\PP(d_l<x)\le(n-1)x^2$, and $\E\log(1/d_l)\le\frac12\log(n-1)+\frac12$ as in Section~\ref{G:s4.2}. By Jensen's inequality, $\E\log M\le\frac12\log\E M^2$, and $\E M^2=\E\max_j|q_j|^2=\sum_{j=1}^n\frac1j\le1+\log n$, the expected maximum of $n$ independent exponential variables. So the expectation of the first sum is at most $\frac\pi2n\log n+\pi n(\log6+\frac12)+\frac\pi2n\log(1+\log n)$.

For the second sum, fix $i$, let $s=|q_i|$ and $\kappa=e^{-(s+1)^2}$, so that $\pi p\ge\kappa$ on $D(q_i,1)$. Fix $l\ne i$, let $e$ be the distance from $q_i$ to the nearest of the remaining $n-2$ points, $x=|q_i-q_l|$, and $\bar e=\min(e,1)$. Conditionally on $q_i$, the variables $e$ and $x$ are independent, and $\PP(e>y\mid q_i)\le(1-\kappa y^2)^{n-2}\le e^{-(n-2)\kappa y^2}$ for $y\le1$. So $\E(\bar e^2\mid q_i)\le\int_0^1e^{-(n-2)\kappa\sigma}d\sigma\le c$ and $\PP(e>1\mid q_i)\le e^{-(n-2)\kappa}\le c$, with $c=1/((n-2)\kappa)$. Since $d_i\le\min(e,x)$, we have $d_i^2/x^2\le\min(1,\bar e^2/x^2)+\mathbf 1_{e>1}$, and by Jensen's inequality for the concave function $\sigma\mapsto\min(1,\sigma/x^2)$, $\E(d_i^2/x^2\mid q_i,q_l)\le\min(1,c/x^2)+c$. Integrating over $q_l$, whose density is at most $\pi^{-1}$, gives at most $c+\int_{\sqrt c\le|z|\le1}c|z|^{-2}\pi^{-1}dA+c+c=c(3+\log^+(1/c))$ when $c\le1$, and at most $2$ in any case. Summing over $l$, $\E(d_i^2T_i\mid q_i)\le(n-1)\min(2,c(3+\log^+(1/c)))$. If $c\le1$, this is at most $\frac{n-1}{n-2}e^{(s+1)^2}(3+\log(n-2))\le2e^{(s+1)^2}(3+\log n)$, using $\log(1/c)\le\log(n-2)$. If $c>1$, then $e^{(s+1)^2}>n-2$, and $2(n-1)\le4e^{(s+1)^2}\le2e^{(s+1)^2}(3+\log n)$. So in all cases $\E(d_i^2T_i\mid q_i)\le2e^{(|q_i|+1)^2}(3+\log n)$, and by the conditional Cauchy--Schwarz inequality $\E(d_i\sqrt{T_i})\le\sqrt{2(3+\log n)}\,\E\,e^{(|q_i|+1)^2/2}$. Finally, $\E\,e^{(|q|+1)^2/2}=\int_0^\infty2se^{-s^2}e^{(s+1)^2/2}ds=2e\int_0^\infty se^{-(s-1)^2/2}ds\le2e\int_{-\infty}^\infty(|u|+1)e^{-u^2/2}du=2e(2+\sqrt{2\pi})<25$, where we put $u=s-1$ and used $s\le|u|+1$. The expectation of the second sum is therefore at most $2\pi n\cdot25\sqrt{2(3+\log n)}=50\sqrt2\pi n\sqrt{3+\log n}$. Adding the three contributions gives the upper bound of Theorem \ref{D:thm:E}.

\subsection{Convergence in probability and in mean}\label{D:s4.3}
With $R,\delta,S(\rho)$ as in the lower bound, $I\ge\frac12\sum_lZ_l-\frac12\int_{|\rho|<R}|B_S|\,dA$, where $Z_l=\int_{|\rho|<R,\,|\rho-q_l|\ge\delta}|\rho-q_l|^{-2}dA$. The $Z_l$ are independent and identically distributed, with $0\le Z_l\le2\pi\log(R/\delta)+\pi$, since the part of the domain with $|\rho-q_l|\le R$ contributes at most $2\pi\log(R/\delta)$ and the rest has area at most $\pi R^2$ and integrand at most $R^{-2}$. Hence $\operatorname{Var}(\frac12\sum Z_l)=O(n\log^2n)$, while $\E\frac12\sum Z_l\ge L_n-O(n\log\log n)$ with $L_n=\frac\pi2n\log n$, by the computation of the lower bound. Chebyshev's inequality gives $\PP(\frac12\sum Z_l<(1-\frac\varepsilon2)L_n)\to0$. Also $\E\int_{|\rho|<R}|B_S|\le4n\pi R^2=4\pi n\log\log n$, and Markov's inequality gives $\PP(\frac12\int|B_S|>\frac\varepsilon2L_n)\le\frac{8\log\log n}{\varepsilon\log n}\to0$. So $\PP(X_n<1-\varepsilon)\to0$ for $X_n=I/L_n$, and the upper bound gives $\E X_n\le1+o(1)$. As in Section~\ref{G:s4.3}, $\E|X_n-1|\le\E(X_n-1)+2\E(1-X_n)^+\to0$. This proves Theorem \ref{D:thm:E}.

\section{Real random lines}\label{D:s5}\label{D:sec:real}

Let now $q_1,\dots,q_n$ be independent and uniform on $[0,1]$, $n\ge16$. The line $\ell_q$ is real, so $\Hm^2(\A(\ell_q))=\frac12I(q)$. We write $\rho=x+iy$.

\subsection{The lower bound}\label{D:s5.1}
Let $\eta=1/\log n$, $y_0=1/n$, and let $Q$ be the region $\eta<x<1-\eta$, $y_0<|y|<1$. For $\rho\in Q$, all points satisfy $|\rho-q_l|\ge|y|$, and we use $J\ge\frac12(A-|B|)$ with the full sums. First, $\E A=n\int_0^1\frac{dq}{(x-q)^2+y^2}=\frac n{|y|}(\arctan\frac x{|y|}+\arctan\frac{1-x}{|y|})\ge\frac{n\pi}{|y|}-\frac nx-\frac n{1-x}$, by the inequality $\arctan a\ge\frac\pi2-\frac1a$ for $a>0$. Next, $\E B=n\int_0^1(\rho-q)^{-2}dq=n(\frac1{\rho-1}-\frac1\rho)=\frac n{\rho(\rho-1)}$, whose modulus is at most $\frac n{x(1-x)}$. Finally $\operatorname{Var}B\le n\int_\R\frac{dq}{(q^2+y^2)^2}=\frac{n\pi}{2|y|^3}$, so $\E|B|\le\frac n{x(1-x)}+\sqrt{\frac{n\pi}2}|y|^{-3/2}$. Hence $\E J(\rho)\ge\frac{n\pi}{2|y|}-\frac n{2x}-\frac n{2(1-x)}-\frac n{2x(1-x)}-\frac12\sqrt{\frac{n\pi}2}|y|^{-3/2}$ on $Q$.

We integrate over $Q$. The first term gives $\frac{n\pi}2(1-2\eta)\cdot2\log\frac1{y_0}=n\pi(1-2\eta)\log n\ge\pi n\log n-2\pi n$. Since $\frac1{x(1-x)}=\frac1x+\frac1{1-x}$ and $\int_\eta^{1-\eta}\frac{dx}x\le\log\frac1\eta$, the next three terms give at least $-2\cdot2n\log\frac1\eta=-4n\log\log n$, the factor $2$ being the total length of the $y$-range. The last term gives at least $-\frac12\sqrt{\frac{n\pi}2}\cdot2\cdot2y_0^{-1/2}=-\sqrt{2\pi}\,n$. So $\E I\ge\pi n\log n-4n\log\log n-(2\pi+\sqrt{2\pi})n$, and halving gives the lower bound of Theorem \ref{D:thm:F}.

\subsection{The upper bound}\label{D:s5.2}
The points lie in the closed unit disc, so the argument of Section~\ref{G:s4.2} applies verbatim with two changes in the expectations. First, $\PP(d_l<x)\le2(n-1)x$, so $\E\log(1/d_l)\le\int_0^\infty\min(1,2(n-1)e^{-u})du=\log(2(n-1))+1$, and the first sum has expectation at most $\pi n(\log6+\log2+1+\log n)$. Second, for $i\ne l$ and $e$ as before, $\PP(e>y\mid q_i)\le(1-y)^{n-2}\le e^{-(n-2)y}$ for $y\le1$, since $[q_i-y,q_i+y]\cap[0,1]$ has length at least $y$, and $e\le1$. So $\E(e^2\mid q_i)\le\int_0^1e^{-(n-2)\sqrt\sigma}d\sigma\le\frac2{(n-2)^2}=:c$, and, the density of $q_l$ on $\R$ being at most $1$, $\E\min(1,c/x^2)\le\int_\R\min(1,c/t^2)dt=4\sqrt c$. Hence $\E(d_i^2T_i)\le(n-1)4\sqrt2/(n-2)\le8$ for $n\ge16$, and $\E(d_i\sqrt{T_i})\le\sqrt8$, so the second sum has expectation at most $2\pi n\sqrt8=4\sqrt2\pi n$. With the bound $\frac{11\pi}{\sqrt2}n$ for $|\rho|>2$, we get $\E I\le\pi n\log n+\pi(\log12+1+4\sqrt2+\frac{11}{\sqrt2})n$, and halving gives the upper bound of Theorem \ref{D:thm:F}.

\subsection{Convergence}\label{D:s5.3}
On $Q$, $I\ge\frac12\sum_lZ_l-\frac12\int_Q|B|$, with $Z_l=\int_Q|\rho-q_l|^{-2}dA\le2\pi\log(2/y_0)$, since $y_0\le|\rho-q_l|\le2$ on $Q$. Then $\operatorname{Var}(\frac12\sum Z_l)=O(n\log^2n)$, and $\E\int_Q|B|=O(n\log\log n)$ by the estimates of the lower bound. Chebyshev's and Markov's inequalities and the upper bound give the convergence as in Section \ref{D:sec:gauss}, with $L_n=\pi n\log n$ for $I$ and therefore $\frac\pi2n\log n$ for the area.

\begin{remark}\label{D:rem:general}
The proofs use the specific densities only through four properties: a bound on the density, a lower bound for it on small discs, the behaviour of the mean of $B_S$, controlled through Lemma \ref{D:lem:beurlingp} or its analogue on the line, and the tail of $\max_j|q_j|$. We expect that they extend, with the same leading term $\frac\pi2n\log n$, to any bounded density $p$ on $\C$ with bounded second derivatives and a finite moment $\E\,\kappa(q)^{-1/2}<\infty$, where $\kappa(z)=\pi\min_{D(z,1)}p$, which in the Gaussian case is $\E\,e^{(|q|+1)^2/2}$, together with $\E\log\max_j|q_j|=O(\log\log n)$, and to any density on a compact interval of $\R$ that is bounded above and below. We have not carried out these extensions, and we have not attempted to determine the minimal hypotheses.
\end{remark}

\section{Topology of the volume-extremal varieties: the case $n=2k$}\label{T:s1}\label{T:sec:half}

In this section and the next one we describe the varieties at which the bounds of the paper are attained, in terms of the topology of the amoeba map, of the symmetries of the variety and, for curves, of the amoeba itself. In this section $V\subset(\C^*)^{2k}$ is a smooth complete intersection of dimension $k$, and we use the notation of Section~\ref{L:sec:log}. So $N(x)$ is the number of points of $V$ over $x$, $S$ is the critical set of $\Log|_V$, $U=\R^{2k}\setminus\Log(S)$, $N_\pi(\theta)$ is the number of points of $V$ whose arguments are congruent to $\theta$ modulo $\pi$, and $Q_\pi=[0,\pi)^{2k}$. By Lemma~\ref{L:lem:area}, $\Log(S)$ is a closed null set, $N$ is finite and locally constant on $U$, $N\ge1$ on $\A(V)\cap U$, and $\Vol(\A(V))=\Vol(\A(V)\cap U)\le\int N=\mathcal M(V)$.

We say that $V$ is \emph{volume-maximal} if $\mathcal M(V)>0$ and $\Vol(\A(V))=\frac{\pi^{2k}}2M_+(V)$, which is the upper bound of [Mik, Theorem 1] (Corollary~\ref{T:cor:bounds} below). We say that $V$ is \emph{volume-minimal} if $\mathcal M(V)\ge\pi^{2k}$ and $\Vol(\A(V))=\pi^{2k}/M_-(V)$, that is, if $V$ attains the lower bound which Theorem~A(b) gives under the hypotheses of Theorems B, C and D. Let $\mathbb S=\{t\in(\C^*)^{2k}:|t_1|=\dots=|t_{2k}|=1\}$ be the compact torus. We say that $V$ is \emph{real up to multiplicative translation} if it is invariant under an antiholomorphic involution $c_t(z)=t\bar z=(t_j\bar z_j)_j$ with $t\in\mathbb S$. This is the case when $sV$ is defined by polynomials with real coefficients for some $s\in(\C^*)^{2k}$, with $t=\bar s/s$. Since $|t_j|=1$, the involution $c_t$ preserves $\Log$, so it permutes the points of every fibre of $\Log|_V$.

\subsection{Parity of the fibres}\label{T:s1.1}

\begin{theorem}\label{T:thm:parity}
Let $V$ be a smooth complete intersection of dimension $k$ in $(\C^*)^{2k}$. Then $N(x)$ is even for every $x\in U$. If moreover $k\ge2$ and $V$ is real up to multiplicative translation, then $N(x)$ is divisible by $4$ for every $x\in U$.
\end{theorem}

\begin{proof}
We follow the cobordism argument of Milnor for the degree modulo $2$ [Mil, \S4], which applies to the proper map $f=\Log|_V$. Let $y_0,y_1\in U$. There is a smooth family $h_t$, $t\in[0,1]$, of diffeomorphisms of $\R^{2k}$, all equal to the identity outside a fixed ball $B$, with $h_0=\mathrm{id}$ and $h_1(y_0)=y_1$. The map $G:V\times[0,1]\to\R^{2k}$, $G(z,t)=h_t^{-1}(\Log z)$, is smooth and proper, since $G^{-1}(K)\subset\Log^{-1}(K\cup\bar B)\times[0,1]$ for every compact set $K$. Since $N$ is locally constant on the open set $U$, there is a ball $D$ around $y_0$ with $D\subset U$, $h_1(D)\subset U$, $N=N(y_0)$ on $D$ and $N=N(y_1)$ on $h_1(D)$. By Sard's theorem we may choose $y\in D$ which is a regular value of $G$ and of its restrictions to $V\times\{0\}$ and $V\times\{1\}$. Then $X=G^{-1}(y)$ is a compact one-dimensional submanifold of $V\times[0,1]$ whose boundary is $X\cap(V\times\{0,1\})=F_0\sqcup F_1$, where $F_0=f^{-1}(y)\times\{0\}$ and $F_1=f^{-1}(h_1(y))\times\{1\}$. So $\#F_0=N(y_0)$ and $\#F_1=N(y_1)$. A compact one-manifold is a finite union of circles and arcs, and every arc has two end points. Hence $N(y_0)\equiv N(y_1)$ modulo $2$.

The amoeba is closed and $\Vol(\A(V))\le\mathcal M(V)\le\pi^{2k}M_+(V)<\infty$ by Theorem~A(a), so $\A(V)\neq\R^{2k}$. A point $y_1\notin\A(V)$ lies in $U$ and has $N(y_1)=0$. So $N(x)$ is even for every $x\in U$.

Now let $k\ge2$ and let $c=c_t$ be an antiholomorphic involution preserving $V$, with $t\in\mathbb S$. Its fixed set $\operatorname{Fix}(c)$ in $V$ is empty or a real-analytic submanifold of dimension $k$, since near a fixed point $c$ is conjugate to the complex conjugation of $\C^k$. Because $\Log\circ c=\Log$, we have $G(c(z),t)=G(z,t)$, so the involution $c\times\mathrm{id}$ maps $X$ to itself. The set $G(\operatorname{Fix}(c)\times[0,1])$ is the image of a manifold of dimension $k+1<2k$ under a smooth map, hence a null set, and we may choose $y$ outside it. Then $c\times\mathrm{id}$ acts on $X$ without fixed points. A continuous involution of a compact interval has a fixed point, so no arc of $X$ is mapped to itself, and the arcs come in pairs $\{a,c(a)\}$. Count the pairs according to the positions of their end points. Let $\alpha_0$, $\alpha_1$ and $\beta$ be the numbers of pairs of arcs with both end points in $F_0$, with both end points in $F_1$, and with one end point in each, respectively. Then $\#F_0=4\alpha_0+2\beta$ and $\#F_1=4\alpha_1+2\beta$, so $N(y_0)/2\equiv N(y_1)/2$ modulo $2$. Taking $y_1\notin\A(V)$ as before gives $N(y_0)\equiv0$ modulo $4$.
\end{proof}

For $k=1$ the second statement fails: a real line has $N=2$. In the proof, the image of $\operatorname{Fix}(c)\times[0,1]$ is then two-dimensional and cannot be avoided.

\begin{corollary}\label{T:cor:bounds}
Let $V$ be a smooth complete intersection of dimension $k$ in $(\C^*)^{2k}$ with $\mathcal M(V)>0$. Then $M_-(V)\ge2$ and
$$\frac{\mathcal M(V)}{2\lfloor M_-(V)/2\rfloor}\ \le\ \Vol(\A(V))\ \le\ \frac{\mathcal M(V)}2\ \le\ \frac{\pi^{2k}}2M_+(V).$$
The second inequality is an equality if and only if $N=2$ on $\A(V)\cap U$. If moreover $k\ge2$ and $V$ is real up to multiplicative translation, then $M_-(V)\ge4$ and
$$\frac{\mathcal M(V)}{4\lfloor M_-(V)/4\rfloor}\ \le\ \Vol(\A(V))\ \le\ \frac{\mathcal M(V)}4\ \le\ \frac{\pi^{2k}}4M_+(V).$$
\end{corollary}

\begin{proof}
On $\A(V)\cap U$ we have $1\le N\le M_-(V)$ by Proposition~\ref{L:prop:fibre}, and $N$ is even by Theorem~\ref{T:thm:parity}, so $2\le N\le2\lfloor M_-(V)/2\rfloor$. The set $\A(V)\cap U$ is not empty, since $\int N=\mathcal M(V)>0$, so $M_-(V)\ge2$. Integrating over $\A(V)\cap U$, whose volume is $\Vol(\A(V))$, gives $2\Vol(\A(V))\le\mathcal M(V)\le2\lfloor M_-(V)/2\rfloor\Vol(\A(V))$, and $\mathcal M(V)\le\pi^{2k}M_+(V)$ is Theorem~A(a). Equality in $2\Vol(\A(V))\le\int N$ holds if and only if $N=2$ almost everywhere on $\A(V)\cap U$, and then everywhere on this open set, since $N$ is locally constant there. In the real case, $N\equiv0$ modulo $4$ on $U$ gives $4\le N\le4\lfloor M_-(V)/4\rfloor$ on $\A(V)\cap U$, and the same argument applies.
\end{proof}

The upper bound $\Vol(\A(V))\le\frac{\pi^{2k}}2M_+(V)$ is [Mik, Theorem 1], and the parity of $N$ is the reason for the factor $\frac12$ in it. The other inequalities of Corollary~\ref{T:cor:bounds} are new. When $M_-(V)$ is odd, the lower bound improves that of Theorem~A(b), so the bound $\mathcal M(V)/M_-(V)$ is never attained. For real $2$-planes in general position in $(\C^*)^4$, where $\mathcal M=\pi^4$ and $M_-=6$, the real lower bound gives $\Vol(\A(V))\ge\pi^4/4$. This is the exact value of the volume of the amoeba of a generic real plane [NP, Theorem 1.2], [MN2, Corollary 7.1]. The real upper bound is attained for $k=2$. Let $C_1,C_2$ be curves with $\Area(\A(C_i))=\pi^2\Area(\Delta_i)$, for instance simple Harnack curves [MR]. By Corollary~\ref{T:cor:bounds} for $k=1$ they have $N=2$ on the regular values of their amoebas and $\mathcal M(C_i)=2\Area(\A(C_i))=\pi^2M_+(C_i)$. By Lemma~\ref{L:lem:product}, the real surface $C_1\times C_2\subset(\C^*)^4$ has $N=4$ on the regular values of its amoeba, $\mathcal M=\pi^4M_+$, and $\Vol(\A(C_1\times C_2))=\frac{\pi^4}4M_+(C_1\times C_2)$.

\subsection{The extremal varieties and the amoeba map}\label{T:s1.2}

\begin{theorem}\label{T:thm:extremal}
Let $V$ be a smooth complete intersection of dimension $k$ in $(\C^*)^{2k}$.

\emph{(a)} $V$ is volume-maximal if and only if $M_+(V)>0$, $N=2$ at every point of $\A(V)\cap U$, and $N_\pi=M_+(V)$ almost everywhere on $Q_\pi$. In that case $\mathcal M(V)=\pi^{2k}M_+(V)$, the restriction of $\Log$ to $\Log^{-1}(\A(V)\cap U)\cap V$ is a two-sheeted covering of $\A(V)\cap U$, and, if $k\ge2$, $V$ is not real up to multiplicative translation.

\emph{(b)} Suppose that $\mathcal M(V)\ge\pi^{2k}$. Then $V$ is volume-minimal if and only if $N=M_-(V)$ at every point of $\A(V)\cap U$ and $\mathcal M(V)=\pi^{2k}$. In that case the restriction of $\Log$ to $\Log^{-1}(\A(V)\cap U)\cap V$ is a covering of $\A(V)\cap U$ with $M_-(V)$ sheets, $M_-(V)$ is even, and, if $k\ge2$ and $V$ is real up to multiplicative translation, $M_-(V)$ is divisible by $4$.
\end{theorem}

\begin{proof}
(a) By Corollary~\ref{T:cor:bounds}, $\Vol(\A(V))\le\mathcal M(V)/2\le\frac{\pi^{2k}}2M_+(V)$. The first inequality is an equality exactly when $N=2$ on $\A(V)\cap U$. The second is an equality exactly when $\int_{Q_\pi}N_\pi=\pi^{2k}M_+(V)$, that is, since $N_\pi\le M_+(V)$ almost everywhere by Proposition~\ref{L:prop:argbound}, when $N_\pi=M_+(V)$ almost everywhere. Under these conditions $\mathcal M(V)=\pi^{2k}M_+(V)>0$ if and only if $M_+(V)>0$. The covering statement is Lemma~\ref{L:lem:area}. If $k\ge2$ and $V$ were real up to multiplicative translation, Corollary~\ref{T:cor:bounds} would give $\Vol(\A(V))\le\frac{\pi^{2k}}4M_+(V)<\frac{\pi^{2k}}2M_+(V)$, since $M_+(V)\ge\mathcal M(V)/\pi^{2k}>0$.

(b) We have $\Vol(\A(V))\ge\mathcal M(V)/M_-(V)\ge\pi^{2k}/M_-(V)$, where the first inequality comes from $N\le M_-(V)$ on $\A(V)\cap U$. Equality in the first holds if and only if $N=M_-(V)$ almost everywhere on $\A(V)\cap U$, hence everywhere on this open set, and equality in the second if and only if $\mathcal M(V)=\pi^{2k}$. Then $M_-(V)$ is a value of $N$ on $U$, so it is even, and it is divisible by $4$ in the real case with $k\ge2$, by Theorem~\ref{T:thm:parity}.
\end{proof}

The condition $N_\pi=M_+(V)$ almost everywhere, that is $\mathcal M(V)=\pi^{2k}M_+(V)$, is Mikhalkin's definition of a multiHarnack variety [Mik, Definition 3.2]. Theorem~\ref{T:thm:extremal}(a) therefore says that a variety is volume-maximal if and only if it is multiHarnack and its amoeba map is two to one. For $k\ge2$, products of curves of maximal area are multiHarnack by Lemma~\ref{L:lem:product}, but their amoeba maps are $2^k$ to one. We do not know whether volume-maximal varieties exist for $k\ge2$; by Theorem~\ref{T:thm:extremal}(a) they cannot be real.

\subsection{Symmetries of extremal varieties}\label{T:s1.3}

Let $\Gamma_V^+=\{t\in\mathbb S:tV=V\}$. For curves this is the group $\Gamma_C^+$ of Section~\ref{L:sec:curve-equality}. Every $t\in\Gamma_V^+$ preserves $\Log$, so $\Gamma_V^+$ acts on every fibre of $\Log|_V$, and it acts freely on $V$.

\begin{proposition}\label{T:prop:group}
If $\mathcal M(V)>0$, then $\Gamma_V^+$ is finite and $|\Gamma_V^+|\le M_+(V)$.
\end{proposition}

\begin{proof}
Let $H\subset\Gamma_V^+$ be a finite subgroup. Since $\int_{Q_\pi}N_\pi=\mathcal M(V)>0$, and by the proof of Proposition~\ref{L:prop:argbound}, there is $\theta\in Q_\pi$ with $N_\pi(\theta)\ge1$ such that all points of $V$ with arguments congruent to $\theta$ modulo $\pi$ are regular points of $\Log|_V$. Let $z$ be such a point, so that $\tau_\theta(z)=z$. For $t\in H$ we have $\tau_\theta(tz)=e^{2i\theta}\bar t\bar z=t^{-1}\tau_\theta(z)=t^{-1}z\in V$, so $tz\in V\cap\tau_\theta(V)$. In the same way $t^{-1}\tau_\theta(V)=\tau_\theta(tV)=\tau_\theta(V)$, so multiplication by $t^{-1}$ preserves both $V$ and $\tau_\theta(V)$ and maps $tz$ to $z$. The two varieties meet transversally at $z$ by Lemma~\ref{L:lem:transverse}, hence also at $tz$. So the $|H|$ points of the orbit $Hz$ are isolated points of $V\cap\tau_\theta(V)$. This intersection is defined by $f_1,\dots,f_k$ and the polynomials $\overline{f_j(\tau_\theta(w))}$, whose Newton polytopes are $\Delta_1,\dots,\Delta_k,\Delta_1,\dots,\Delta_k$. Bernstein's bound for isolated zeros gives $|H|\le M_+(V)$. A closed subgroup of $\mathbb S$ of positive dimension contains finite subgroups of every order, so $\Gamma_V^+$ is finite, and the bound applies to $H=\Gamma_V^+$.
\end{proof}

For curves the bound is attained exactly by the trinomial curves, and this characterises the volume-minimal curves by their symmetries.

\begin{theorem}\label{T:thm:curvegroup}
Let $f$ be an irreducible Laurent polynomial in two variables with two-dimensional Newton polygon $\Delta$ and smooth zero set $C$. Then $|\Gamma_C^+|\le2\Area(\Delta)=M_+(C)$, with equality if and only if $f$ is a monomial times a trinomial. If $\Delta$ is not centrally symmetric, $C$ is volume-minimal if and only if $|\Gamma_C^+|=2\Area(\Delta)$.
\end{theorem}

\begin{proof}
Let $H=\Gamma_C^+$. It is finite. Otherwise it contains a circle subgroup, whose Zariski closure is a one-dimensional subtorus $T_1$ preserving $C$. Then the irreducible curve $C$ is a coset of $T_1$, and $f$ would be a binomial, with one-dimensional Newton polygon. By Lemma~\ref{L:lem:subgroup} there is an integer matrix $P$ with $|\det P|=|H|$ such that $\phi_P$ is surjective with kernel $H$. The map $\phi_P$ is finite, so $W=\phi_P(C)$ is an irreducible curve, and $C=\phi_P^{-1}(W)$ because $C$ is invariant under $H$, whose cosets are the fibres of $\phi_P$. Let $g$ be an irreducible Laurent polynomial defining $W$. Since $\phi_P$ is \'etale, $g\circ\phi_P$ is reduced, and it defines $C$. So $f$ and $g\circ\phi_P$ differ by a unit of the ring of Laurent polynomials, a constant times a monomial. By Section~\ref{L:s2.3}, the Newton polygon of $g\circ\phi_P$ is $P^{\mathsf T}\Newt(g)$, so $\Area(\Delta)=|H|\Area(\Newt(g))$. The polygon $\Newt(g)$ is two-dimensional, because $\Delta$ is, and a two-dimensional lattice polygon has area at least $\frac12$, with equality exactly for unimodular triangles. Hence $|H|\le2\Area(\Delta)$. Equality holds if and only if $\Newt(g)$ is a unimodular triangle. A unimodular triangle contains no lattice points other than its vertices, so in that case $g$ is a trinomial, and $f$ is a monomial times the trinomial $g\circ\phi_P$. Conversely, if $f$ is a monomial times a trinomial with Newton triangle $\Delta$, the proof of Proposition~\ref{L:prop:trinomial} gives $C=\phi_P^{-1}(L')$ with $|\det P|=2\Area(\Delta)$. The kernel of $\phi_P$ is a group of roots of unity contained in $\Gamma_C^+$, so $|\Gamma_C^+|\ge2\Area(\Delta)$. The last assertion follows from Theorem~\ref{L:thm:curve-equality}.
\end{proof}

\subsection{The topology of extremal curves}\label{T:s1.4}

For curves the two extremal classes can be described completely, and they turn out to be opposite in every respect. For a lattice polygon $\Delta$ let $i(\Delta)$ and $b(\Delta)$ be the numbers of interior and boundary lattice points. Recall [FPT], [Mi00, Theorem 4] that the components of the complement of the amoeba of a curve with Newton polygon $\Delta$ inject into $\Delta\cap\Z^2$, and that every vertex of $\Delta$ corresponds to a component. So the number of components is at least the number of vertices and at most $i(\Delta)+b(\Delta)$.

\begin{theorem}\label{T:thm:curves}
Let $f$ be an irreducible Laurent polynomial in two variables with two-dimensional Newton polygon $\Delta$ and smooth zero set $C$.

\emph{(a)} If $C$ is volume-maximal, that is, $\Area(\A(C))=\pi^2\Area(\Delta)$, then $N=2$ at every regular value of $\Log|_C$ in $\A(C)$. By [MR], up to multiplicative translation, $C$ is defined over $\R$ and its real locus is in cyclically maximal position in the sense of [Mi00, Definition 2], that is, $C$ is a simple Harnack curve. The topological type of its real locus in the toric surface of $\Delta$ depends only on $\Delta$ [Mi00, Theorem 3]. The components of $\R^2\setminus\A(C)$ are in bijection with $\Delta\cap\Z^2$, and the $i(\Delta)$ bounded ones are discs bounded by the images of the ovals [Mi00, Corollary 10]. So the amoeba has the maximal number $i(\Delta)$ of holes.

\emph{(b)} Suppose that $\Delta$ is not centrally symmetric and that $C$ is volume-minimal, that is, $\Area(\A(C))=\pi^2/\MV(\Delta,-\Delta)$. Then $\Delta$ is a triangle, $f$ is a monomial times a trinomial, and the following hold. First, $C$ is a nondegenerate curve of genus $i(\Delta)$ with $b(\Delta)$ punctures. Second, $\Gamma_C^+$ has order $2\Area(\Delta)$ and acts freely on $C$, and $C\to C/\Gamma_C^+$ is an unramified Galois covering of a line, which is a sphere with three punctures. Third, $\A(C)$ is the image of the amoeba of the line $\{1+z+w=0\}$ under an invertible affine map of $\R^2$. Hence $\R^2\setminus\A(C)$ has exactly three components, which are convex and unbounded, the minimal number, and the amoeba has no holes even when $i(\Delta)>0$. Fourth, every fibre of $\Log|_C$ over a regular value in the interior of $\A(C)$ is a single free orbit of the group $\Gamma_C$ of Section~\ref{L:sec:curve-equality}, with $\MV(\Delta,-\Delta)=4\Area(\Delta)$ points.
\end{theorem}

\begin{proof}
(a) The first assertion is Theorem~\ref{T:thm:extremal}(a), and the others are the cited results of [MR] and [Mi00].

(b) By Theorem~\ref{L:thm:curve-equality}, $\Delta$ is a triangle and $f$ is a monomial times a trinomial. The initial forms of $f$ are monomials and binomials $c_az^a+c_bz^b$. The zero set of such a binomial in $(\C^*)^2$ is a finite union of cosets of a one-dimensional subtorus, hence smooth, so $C$ is nondegenerate in the sense of Definition~\ref{L:def:nondeg}. By Khovanskii's formulas [Kh], a nondegenerate curve with Newton polygon $\Delta$ has genus $i(\Delta)$ and meets the boundary of the toric surface of $\Delta$ in $b(\Delta)$ points. The statements on $\Gamma_C^+$ are Theorem~\ref{T:thm:curvegroup} and its proof, where $C/\Gamma_C^+\simeq\phi_P(C)$ is a curve whose Newton polygon is a unimodular triangle, that is, a line after a monomial change of coordinates and a multiplicative translation. As in the proof of Proposition~\ref{L:prop:trinomial}, $C=\phi_P^{-1}(L')$, where $L'$ is a multiplicative translate of $L=\{1+z+w=0\}$ and $|\det P|=2\Area(\Delta)$. Since $\Log\circ\phi_P=P\circ\Log$ and $\phi_P$ is surjective, $\A(C)=P^{-1}(\A(L'))$, and $\A(L')$ is a translate of $\A(L)$. The complement of $\A(L)$ consists of the three regions $\{e^x+e^y<1\}$, $\{1+e^x<e^y\}$ and $\{1+e^y<e^x\}$, where $1$, $e^x$ and $e^y$ violate a triangle inequality. Each of them is convex, being a sublevel set of a convex function of $(x,y)$, and unbounded. The last assertion is Proposition~\ref{L:prop:single-orbit}, with $\mu=\MV(\Delta,-\Delta)=4\Area(\Delta)$ for a triangle.
\end{proof}

So a volume-maximal curve is real, its amoeba map is two to one, and its amoeba has all the $i(\Delta)$ holes allowed by its Newton polygon. A volume-minimal curve has the largest possible group of symmetries, its amoeba map is $4\Area(\Delta)$ to one, and its amoeba has no holes at all. In both cases the curve itself has genus $i(\Delta)$ and $b(\Delta)$ punctures, like every nondegenerate curve with Newton polygon $\Delta$, so it is the position of the curve with respect to the maps $\Log$ and $\Arg$, and not its intrinsic topology, which distinguishes the extremal curves. In higher dimension the volume-minimal towers are the iterated twisted products of trinomial curves (Theorem E). By Theorem F, such a product need not be a toric product of curves, although an open dense subset of the surface of Theorem F is isomorphic to an open subset of a product of two lines.

\section{Topology of the volume-extremal varieties: the case $n>2k$}\label{T:s2}\label{T:sec:ngt}

Let $V\subset(\C^*)^n$, $n>2k$, be an irreducible smooth complete intersection of dimension $k$ with $\Hm^{2k}(\A(V))>0$, and use the notation of Sections~\ref{N:s1} and \ref{N:s2}. By Theorem~\ref{N:thm:main},
$$\frac1{N_{\max}}\Big(\sum_{|I|=2k}M_I(V)^2\Big)^{1/2}\ \le\ \Hm^{2k}(\A(V))\ \le\ \frac1{N_{\min}}\sum_{|I|=2k}M_I(V).\eqtag{T:e1}$$
We say that $V$ is volume-minimal, respectively volume-maximal, in $(\C^*)^n$ if the first, respectively the second, inequality in (\ref{T:e1}) is an equality. We show that this can only happen when $V$ lies in a translate of a subtorus of dimension $2k$, so that the question reduces to the half-dimensional case of Section~\ref{T:sec:half}.

Write $x=\Log z$. For a coordinate set $I=\{i_1<\dots<i_{2k}\}$ let $\omega_I$ be the restriction to $V$ of the real $2k$-form $dx_{i_1}\wedge\dots\wedge dx_{i_{2k}}$. Then $J_I=|\omega_I|/d\mathrm{vol}$, and $J=(\sum_IJ_I^2)^{1/2}$ by (\ref{N:e2.2}). The variety $V$ is a connected complex manifold, being smooth and irreducible, and the forms $\omega_I$ are real-analytic. So each $\omega_I$ either vanishes identically on $V$ or vanishes only on a closed subset with empty interior and measure zero, and the same holds for $J$. Since $\int_VJ\ge\Hm^{2k}(\A(V))>0$, $J$ does not vanish identically.

\begin{lemma}\label{T:lem:minkowski}
Let $F$ be an integrable function on a measure space with values in $\R^m$, with $\int F\neq0$. Then $\|\int F\|\le\int\|F\|$ for the Euclidean norm, with equality if and only if $F=\|F\|u$ almost everywhere for a fixed unit vector $u$.
\end{lemma}

\begin{proof}
Let $u=\int F/\|\int F\|$. Then $\|\int F\|=\langle u,\int F\rangle=\int\langle u,F\rangle\le\int\|F\|$, with equality if and only if $\langle u,F\rangle=\|F\|$ almost everywhere, that is, $F=\|F\|u$ almost everywhere.
\end{proof}

\begin{theorem}\label{T:thm:upper}
The variety $V$ is volume-maximal in $(\C^*)^n$ if and only if $N=N_{\min}$ almost everywhere on $\A(V)$ and there is a coordinate set $I$ with $|I|=2k$ such that every coordinate $z_j$, $j\notin I$, is constant on $V$. In that case $\pi_I|_V$ is an isomorphism onto a smooth irreducible variety $W\subset(\C^*)^{2k}$ of dimension $k$, $\A(V)$ is the image of $\A(W)$ under an isometric affine embedding of $\R^{2k}$ into $\R^n$, $M_I(V)=\mathcal M(W)$, $M_{I'}(V)=0$ for $I'\neq I$, and the equality reads $\Vol(\A(W))=\mathcal M(W)/N_{\min}$.
\end{theorem}

\begin{proof}
By the area formula (\ref{N:e2.1}), $\int_VJ=\int N\,d\Hm^{2k}\ge N_{\min}\Hm^{2k}(\A(V))$, and $\int_VJ\le\sum_I\int_VJ_I=\sum_IM_I(V)$, since $(\sum_IJ_I^2)^{1/2}\le\sum_IJ_I$. So equality in (\ref{T:e1}) holds if and only if $N=N_{\min}$ almost everywhere on $\A(V)$ and $(\sum_IJ_I^2)^{1/2}=\sum_IJ_I$ almost everywhere on $V$, that is, if at almost every point of $V$ at most one $J_I$ is nonzero. If two forms $\omega_I$ and $\omega_{I'}$ were not identically zero, both would be nonzero outside a null set, which is impossible. Since $J\not\equiv0$, there is exactly one $I$ with $\omega_I\not\equiv0$. Let $z\in V$ be a point where $\omega_I\ne0$. The differentials $dx_i|_V$, $i\in I$, form a basis of the real cotangent space at $z$, and for $j\notin I$ we write $dx_j|_V=\sum_{i\in I}a_i\,dx_i|_V$ at $z$. For $i_0\in I$ and $I'=(I\setminus\{i_0\})\cup\{j\}$ we get $\omega_{I'}=\pm a_{i_0}\omega_I$ at $z$, and $\omega_{I'}\equiv0$, so $a_{i_0}=0$. Hence $dx_j|_V=0$ on the dense set where $\omega_I\neq0$, and so everywhere. Since $V$ is connected, $\log|z_j|$ is constant on $V$, and the holomorphic function $z_j$, having constant modulus on a connected complex manifold, is constant. Conversely, if $z_j$ is constant on $V$ for all $j\notin I$, then every $\omega_{I'}$ with $I'\neq I$ contains a factor $dx_j|_V=0$, so $J=J_I=\sum_{I'}J_{I'}$. In that case $V=\{(w,c):w\in W\}$ after reordering the coordinates, with $c\in(\C^*)^{n-2k}$ and $W=\pi_I(V)$. Then $\A(V)=\A(W)\times\{\Log c\}$, the numbers of points over corresponding points agree, $M_I(V)=\int_VJ_I=\mathcal M(W)$ by Lemma~\ref{L:lem:area} applied to $W$, and $M_{I'}(V)=0$ for $I'\neq I$.
\end{proof}

\begin{theorem}\label{T:thm:lower}
The variety $V$ is volume-minimal in $(\C^*)^n$ if and only if $N=N_{\max}$ almost everywhere on $\A(V)$ and $V$ is contained in a translate $pT'$ of an algebraic subtorus $T'\subset(\C^*)^n$ of dimension $2k$. In that case, let $\psi:(\C^*)^{2k}\to T'$ be an isomorphism, $\psi(w)=(w^{a_1},\dots,w^{a_n})$, let $\Lambda$ be the $n\times2k$ integer matrix with rows $a_1,\dots,a_n$, and let $W=\psi^{-1}(p^{-1}V)$. Then $W$ is a smooth irreducible variety of dimension $k$ in $(\C^*)^{2k}$, $\A(V)=\Log p+\Lambda(\A(W))$, $\Hm^{2k}(\A(V))=\|\Lambda\|_{2k}\Vol(\A(W))$ and $\int_VJ=\|\Lambda\|_{2k}\mathcal M(W)$, where $\|\Lambda\|_{2k}=\sqrt{\det\Lambda^{\mathsf T}\Lambda}$, and the equality reads $\Vol(\A(W))=\mathcal M(W)/N_{\max}$.
\end{theorem}

\begin{proof}
By the area formula, $\Hm^{2k}(\A(V))\ge N_{\max}^{-1}\int N\,d\Hm^{2k}=N_{\max}^{-1}\int_VJ$, with equality if and only if $N=N_{\max}$ almost everywhere on $\A(V)$. By Lemma~\ref{T:lem:minkowski} applied to $F=(J_I)_I$, $\int_VJ=\int_V\|F\|\ge\|\int_VF\|=(\sum_IM_I(V)^2)^{1/2}$. Here $\int F\neq0$, since otherwise all $M_I(V)$ vanish and $\Hm^{2k}(\A(V))=0$ by Theorem~\ref{N:thm:main}. Equality holds if and only if $(J_I)_I=Ju$ almost everywhere for a fixed unit vector $u$.

Suppose that this holds. Let $V^*\subset V$ be the open dense set where $J\neq0$. At $z\in V^*$ the image $\Pi_z=d\Log(T_zV)$ is a real subspace of dimension $2k$ of $\R^n$, and its Plücker coordinates are the numbers $p_I(z)=\omega_I/d\mathrm{vol}$, with $|p_I(z)|=J_I(z)=J(z)u_I$. Each $p_I$ is continuous and does not vanish on $V^*$ when $u_I>0$, so its sign is constant on each connected component $\Omega$ of $V^*$, and $p_I\equiv0$ when $u_I=0$. Hence the point of the Grassmannian with Plücker coordinates $(p_I(z))_I$ is constant on $\Omega$, and $\Pi_z=\Pi$ is constant there. For every $\nu\in\Pi^\perp$ the real-analytic function $\langle\nu,\Log z\rangle$ has zero differential on $\Omega$, so it is constant on $\Omega$, hence on the connected manifold $V$. Thus $\A(V)\subset x_0+\Pi$ for some $x_0\in\R^n$. Fix $p\in V$ and put $V'=p^{-1}V$, an irreducible variety containing $1$ and contained in the closed subgroup $G_\Pi=\{z:\Log z\in\Pi\}$ of $(\C^*)^n$. The subgroup of $(\C^*)^n$ generated by $V'$ is a closed connected algebraic subgroup [Bo, Proposition I.2.2], hence a subtorus $T'$, and $T'\subset G_\Pi$. The image $\Log(T')$ is a linear subspace of $\R^n$ of dimension $\dim_\C T'$, contained in $\Pi$, so $\dim T'\le2k$. On the other hand $\A(V)\subset\Log p+\Log(T')$ has positive $\Hm^{2k}$-measure, so $\dim T'\ge2k$. So $V\subset pT'$ with $\dim T'=2k$.

Conversely, suppose that $V\subset pT'$ with $\dim T'=2k$, and let $\psi$, $\Lambda$ and $W$ be as in the statement. Since $\log|\psi(w)_j|=\langle a_j,\Log w\rangle$, we have $\Log\circ\psi=\Lambda\circ\Log$, so $\A(V)=\Log p+\Lambda(\A(W))$, and $\Lambda$ is injective because $\psi$ is. Then $d\Log(T_zV)\subset\Lambda(\R^{2k})$ for every $z$, so all the Plücker vectors $(p_I(z))_I$ are multiples of the Plücker vector of $\Lambda(\R^{2k})$, and $(J_I)_I=Ju$ for its normalised absolute value $u$. This is the equality case of Lemma~\ref{T:lem:minkowski}. The map $z\mapsto\psi^{-1}(p^{-1}z)$ is an isomorphism of $V$ onto $W$ which identifies the fibres of $\Log|_V$ over $\Log p+\Lambda y$ with the fibres of $\Log|_W$ over $y$. An injective linear map multiplies $2k$-dimensional measures by $\|\Lambda\|_{2k}$, so $\Hm^{2k}(\A(V))=\|\Lambda\|_{2k}\Vol(\A(W))$. By the area formula and Lemma~\ref{L:lem:area}, $\int_VJ=\int N\,d\Hm^{2k}=\|\Lambda\|_{2k}\int N_W=\|\Lambda\|_{2k}\mathcal M(W)$. Dividing the equality case of (\ref{T:e1}) by $\|\Lambda\|_{2k}$ gives the last assertion.
\end{proof}

\begin{corollary}\label{T:cor:strict}
If $V$ is not contained in a translate of a subtorus of dimension $2k$, both inequalities in (\ref{T:e1}) are strict. This holds for every $k$-plane in general position in $(\C^*)^n$, $n>2k$, in the sense of Corollary~\ref{N:cor:planes}, and in particular for every line in general position in $(\C^*)^n$, $n\ge3$, real or not. Consequently the inequalities of Corollaries~\ref{N:cor:planes} and \ref{F:cor:noNmax} for such planes are strict.
\end{corollary}

\begin{proof}
A coordinate subtorus translate is a translate of a subtorus, so the first assertion follows from Theorems~\ref{T:thm:upper} and \ref{T:thm:lower}. Let $V$ be a $k$-plane in general position, and suppose that $V\subset pT'$ with $\dim T'=2k<n$. There is a nonzero $\nu\in\Z^n$ with $t^\nu=1$ on $T'$, so $z^\nu$ is constant on $V$. On $\bar V\simeq\PP^k$ write $z_j=\ell_j/\ell_0$, where $\ell_j$ is a linear form defining the hyperplane $\lambda_j$. Then $z^\nu=\prod_j\ell_j^{\nu_j}/\ell_0^{\sum_j\nu_j}$ is constant. The hyperplanes $\lambda_0,\dots,\lambda_n$ are pairwise distinct, because two equal ones together with $k-1$ others would form $k+1$ hyperplanes with nonempty intersection. By unique factorisation of homogeneous polynomials, all $\nu_j$ vanish, a contradiction. The bounds of Corollaries~\ref{N:cor:planes} and \ref{F:cor:noNmax} are obtained from (\ref{T:e1}) by bounding $N_{\max}$, $N_{\min}$ and the $M_I$, so they are strict as well.
\end{proof}

For nonreal lines this recovers the strictness in Theorem~\ref{G:thm:C}, and for real lines it shows that $\Hm^2(\A(\ell))<\frac12\sum_IM_I=\frac12\binom n2\pi^2$, in the upper bound of Theorem~\ref{F:thm:lines3}. Theorems~\ref{T:thm:upper} and \ref{T:thm:lower} show that, for $n>2k$, the extremal varieties of (\ref{T:e1}) are the half-dimensional varieties of Section~\ref{T:sec:half} placed in a translate of a subtorus of dimension $2k$, by a coordinate embedding for the upper bound and by an arbitrary monomial embedding for the lower bound. In both cases the extremality of $V$ is equivalent to that of $W\subset(\C^*)^{2k}$ for the bound $\Vol(\A(W))\le\mathcal M(W)/N_{\min}$, respectively $\Vol(\A(W))\ge\mathcal M(W)/N_{\max}$, whose equality cases are those of Corollary~\ref{T:cor:bounds} and Theorem~\ref{T:thm:extremal} when $W$ is a complete intersection. For instance, if $W$ is a complete intersection with $\mathcal M(W)\ge\pi^{2k}$, then $V$ is volume-minimal in $(\C^*)^n$ with $N_{\max}=M_-(W)$ and $\Hm^{2k}(\A(V))=\|\Lambda\|_{2k}\pi^{2k}/M_-(W)$ if and only if $W$ is volume-minimal in the sense of Section~\ref{T:sec:half}.

\section{Open problems}\label{sec:open}

We collect the questions left open by the two parts. In the half-dimensional case, the main open problem is Conjecture~\ref{L:conj:equality}, which describes the equality case of the bound $\Vol(\A(V))\ge\pi^{2k}/M_-(V)$ for all varieties with $\mathcal M(V)\ge\pi^{2k}$. It is proved for towers (Theorem E), and Corollary~\ref{K2:cor:strict} confirms it for $k$-planes in general position, which never attain the bound. The extension of the corner method sketched in Section~\ref{K:s6} would give a necessary condition for equality in general, but it is not proved, and we do not know how to verify the condition beyond planes. We also do not know a description of the tuples of polytopes for which every nondegenerate complete intersection satisfies $\mathcal M(V)\ge\pi^{2k}$, beyond the case of curves, where the answer is that the polygon is not centrally symmetric (Theorem B).

For $k\ge2$ we do not know whether volume-maximal varieties exist, that is, whether Mikhalkin's bound $\Vol(\A(V))\le\frac{\pi^{2k}}2M_+(V)$ is ever attained. By Theorem~\ref{T:thm:extremal} such a variety would be multiHarnack with a two-to-one amoeba map, and it could not be real; for real varieties the sharp bound is $\frac{\pi^{2k}}4M_+(V)$ when $k=2$ (Corollary~\ref{T:cor:bounds}), and its value for $k\ge3$ is unknown. It would also be interesting to know which varieties attain $|\Gamma_V^+|=M_+(V)$ in Proposition~\ref{T:prop:group} for $k\ge2$, and whether they are exactly the isogeny pullbacks of varieties with $M_+=1$.

In the case $n>2k$, the lower bound (\ref{N:e1.1}) uses a single projection $B$, and its optimisation over non-coordinate matrices $B$ is open. Theorem~\ref{F:thm:injective} assumes full-dimensional Newton polytopes, and Remark~\ref{F:rem:general} gives a sufficient condition otherwise; a necessary and sufficient condition for generic injectivity is not known. For $k$-planes with $k\ge2$ in $(\C^*)^n$, $n>2k$, we do not know the minimum of $\Hm^{2k}$ over planes in general position, nor whether real planes are the minimisers, nor the growth rate in $n$, for which Corollary~\ref{F:cor:noNmax} only gives the window between $\binom n{2k}^{1/2}\pi^{2k}$ and $\binom n{2k}\pi^{2k}$.

For lines, three questions stand out. The first is Conjecture~\ref{G:conj:mono}, that the function $g$ is strictly decreasing on $(0,\frac12]$, which would identify the equispaced lines as the unique minimisers in $(\C^*)^3$. Even the attainment of the minimum is open: it would follow from the inequality $\inf g<\frac{1+\sqrt2}2$, for which we have no proof. A natural approach to both questions is to differentiate $g$ after moving the singular point to a fixed position by a diffeomorphism of the plane, as in Section~\ref{G:s3.3}. The second is Conjecture~\ref{G:conj:nlogn}, a deterministic lower bound of order $n\log n$ for every line, supported by the random models, the equispaced lines and the degeneration formula. The third is the continuity of the area on the whole Deligne--Mumford compactification of the configuration space (Remark~\ref{D:rem:DM}), which together with the degeneration formula would describe the infimum of the area. Finally, the error terms in Theorems~\ref{G:thm:B} and \ref{D:thm:E} are of order $n\sqrt{\log n}$ from above, and we expect the true second-order term to be smaller.

\subsection*{Statements and declarations}
The author declares no competing interests. No data were generated or analysed in this study.

\end{document}